\documentclass{article}
\usepackage{amsfonts}
\usepackage{makeidx}
\usepackage{amssymb}
\usepackage{amsmath}
\usepackage{amstext}
\usepackage{graphicx}
\usepackage{latexsym}
\usepackage{amsthm}

\newtheorem{theorem}{Theorem}[section]

\newtheorem{axiom}[theorem]{Axiom}

\newtheorem{corollary}[theorem]{Corollary}

\newtheorem{definition}[theorem]{Definition}

\newtheorem{lemma}[theorem]{Lemma}

\newtheorem{proposition}[theorem]{Proposition}
\newtheorem{remark}[theorem]{Remark}

\input tcilatex
\begin{document}

\title{An improved setting for generalized functions: fine ultrafunctions }
\author{Vieri Benci \thanks{
Dipartimento di Matematica, Universit\`{a} degli Studi di Pisa, Via F.
Buonarroti 1/c, Pisa, ITALY}}
\maketitle

\begin{abstract}
Ultrafunctions are a particular class of functions defined on a Non
Archimedean field $\mathbb{E}\supset \mathbb{R}$. They have been introduced
and studied in some previous works (\cite{ultra},\cite{BBG},\cite{belu2012}%
,..,\cite{bls}). In this paper we develop the notion of fine ultrafunctions
which improves the older definitions in many crucial points. Some
applications are given to show how ultrafunctions can be applied in studing
Partial Differential Equations. In particular, it is possible to prove the
existence of \textit{ultrafunction solutions} to ill posed evolution poblems.

\medskip


\end{abstract}

\tableofcontents

\section{Introduction}

In many circumstances, the notion of real function is not sufficient to the
needs of a theory and it is necessary to extend it. The ultrafunctions are a
kind of generalized functions based on a field $\mathbb{E}$ containing the
field of real numbers $\mathbb{R}$. The field $\mathbb{E}$ (called field of
Euclidean numbers) is a peculiar hyperreal field which satisfies some
properties useful for the purposes of this paper.

The ultrafunctions provide generalized solutions to certain equations which
do not have any solution, not even among the distributions.

We list some of the main properties of the ultrafunctions:

\begin{itemize}
\item the ultrafunctions are defined on a set $\Gamma ,$%
\begin{equation*}
\mathbb{R}^{N}\subset \Gamma \subset \mathbb{E}^{N},
\end{equation*}
and take values in $\mathbb{E}$; actually they form an algebra $V%
{{}^\circ}%
$ over the field $\mathbb{E}$;

\item to every function $f:\mathbb{R}^{N}\rightarrow \mathbb{R}$ corresponds
a unique ultrafunction%
\begin{equation*}
f%
{{}^\circ}%
:\Gamma \rightarrow \mathbb{E}
\end{equation*}%
that extends $f$ to $\Gamma $ and satisfies suitable properties described
below;

\item there exists a linear functional 
\begin{equation*}
\doint :V%
{{}^\circ}%
\rightarrow \mathbb{E}
\end{equation*}%
called pointwise integral such that $\forall f\in C_{c}^{0,1}\left( \mathbb{R%
}^{N}\right) ,$%
\begin{equation*}
\doint f%
{{}^\circ}%
(x)dx=\int f(x)dx
\end{equation*}

\item there are $N$ operators%
\begin{equation*}
D_{i}:V%
{{}^\circ}%
\rightarrow V%
{{}^\circ}%
,\ i=1,...,N
\end{equation*}%
called generalized partial derivatives such that $\forall f\in
C_{c}^{1,1}\left( \mathbb{R}^{N}\right) ,$%
\begin{equation*}
\left( \frac{\partial f}{\partial x_{i}}\right) ^{\circ }=D_{i}f%
{{}^\circ}%
\end{equation*}

\item to every distribution $T\in \mathcal{D}^{\prime }$ corresponds an
ultrafunction $T%
{{}^\circ}%
$ such that $\forall \varphi \in \mathcal{D}$%
\begin{equation*}
\doint T%
{{}^\circ}%
(x)\varphi 
{{}^\circ}%
dx=\left\langle T,\varphi \right\rangle
\end{equation*}%
and 
\begin{equation*}
\left( \frac{\partial T}{\partial x_{i}}\right) ^{\circ }=D_{i}T%
{{}^\circ}%
\end{equation*}

\item if $u$ is the solution of a PDE, then $u%
{{}^\circ}%
$ is the solution of the same PDE "translated" in the framework of
ultrafunctions.

\item $\Gamma $ is a hyperfinite set (see section \ref{hs}) so that we have
enough compactness to prove the existence of a solution for a very large
class of equations which includes many ill posed problems.
\end{itemize}

The ultrafunctions have been recently introduced in \cite{ultra} and
developed in \cite{BBG},\cite{belu2012},..,\cite{bls}. In these papers
different models of ultrafuctions have been analyzed and several
applications have been provided. In the present paper, we introduce an
improved model: the space of \textit{fine} ultrafunctions. The fine
ultrafunctions form an algebra in which the pointwise integral and the
generalized derivative satisfy most of the familiar properties that are
consistent with the algebraic structure of $V%
{{}^\circ}%
$. In particular, these properties allow to solve many evolution problem in
the space $C^{1}(\mathbb{E},V%
{{}^\circ}%
)$ (see sections \ref{TDU} and \ref{ep}).

\bigskip

This paper is organized as follow: in the rest of this introduction, we
frame the theory of ultrafunction and expose our point of view on
Non-Archimedean Mathematics and on the notion of generalized functions.

In section \ref{PN}, we present the preliminary material necessary to the
rest of the paper. In particular we present an approach to Non Standard
Analysis (NSA) suitable for the theory of ultrafunctions. This approach is
based on the notion of $\Lambda $-limit (see also \cite{ultra} and \cite%
{BDN2018}) which leads to the field of Euclidean numbers (see also \cite%
{BL2021}). This part has been written in such a way to be understood also by
a reader who is not familiar with NSA.

In section \ref{U} we recall the notion of ultrafunction, we define the
spaces of fine ultrafunctions and of time dependent ultrafunctions.

The main properties of the fine ultrafunctions are examined in section \ref%
{BP}.

Section \ref{SA} is devoted to some applications that exemplify the use of
ultrafunctions in PDE's.

Section \ref{CSU} is devoted to the proof that the definition of
ultrafunctions is consistent. In fact, even if this definition is based on
notions which appear quite natural, the consistency of these notions is a
delicate issue. We prove this consistency by the construction of a very
involved model; we do not know if a simpler model exists. This section is
very technical and we assume the reader to be used with the techniques of
NSA.

\subsection{Few remarks on Non-Archimedean Mathematics}

The scientific community has always accepted new mathematical entities,
especially if these are useful in the modeling of natural phenomena and in
solving the problems posed by the technique. Some of these entities are the
infinitesimals that have been a carrier of the modern science since the
discovery of the infinitesimal calculus at the end of XVII century. But
despite the successes achieved with their employment, they have been opposed
and even fought by a considerable part of the scientific community (see e.g. 
\cite{bair2,Bla,benciISO}). At the end of the 19th century they were placed
on a more rigorous basis thanks to the works of Du Bois-Reymond \cite{DBR},
Veronese \cite{veronese2}, Levi-Civita \cite{LC} and others, nevertheless
they were fought (and defeated) by the likes of Russell (see e.g. \cite%
{Russell1}) and Peano \cite{peano}. Also the reception of the Non Standard
Analysis created in the '60s by Robinson has not been as good as it
deserved, even though a minority of mathematicians of the highest level has
elaborated interesting theories based on it (see e.g. \cite{ALBE}, \cite%
{nelson}, \cite{tao}).

Personally, I am convinced that the Non-Archimedean Mathematics is branch of
mathematics very rich and allows to construct models of the real world in a
more efficient way. Actually, this is the main motivation of this paper.

\subsection{Few remarks on generalized function}

The intensive use of the Laplace transform in engineering led to the
heuristic use of symbolic methods, called operational calculus. An
influential book on operational calculus was Oliver Heaviside's
Electromagnetic Theory of 1899 \cite{Hea}. During the late 1920s and 1930s
further steps were taken, very important to future work. The Dirac delta
function was boldly defined by Paul Dirac as a central aspect of his
scientific formalism. Jean Leray and Sergei Sobolev, working in partial
differential equations, defined the first adequate theory of generalized
functions and generalized derivative in order to work with weak solutions of
partial differential equations. Sobolev's work was further developed in an
extended form by Laurent Schwartz. Today, among people working in partial
differential equations, the theory of distributions of L. Schwartz is the
most commonly used, but also other notions of generalized functions have
been introduced by J.F. Colombeau \cite{col85} and M. Sato \cite{sa59}.

After the discovery of Non Standard Analysis, many models of generalized
functions based on hyperreal fields appeared. The existence of infinite and
infinitesimal numbers allows to relate the \textit{delta of Dirac }$\delta $%
\textit{\ }to a function which takes an infinite value in a neighborhood of $%
0$ and vanishes in the other points. So, in this context, expression such as 
$\sqrt{\delta _{a}}$ or $\delta _{a}^{2}$ make absolutely sense.

The literature in this context is quite large and, without the hope to be
exhaustive, we refer to the following papers and their references: Albeverio
, Fenstad, Hoegh-Krohn \cite{ALBE}, Nelson \cite{nelson}, Arkeryd, Cutland,
Henson \cite{cutland}, Bottazzi \cite{Bott}, Todorov \cite{todo2011}.\bigskip

\subsection{Notations\label{sec:not}}

For the sake of the reader, we list the main notation used in this paper. If 
$X$ is any set and $\Omega $\ is a measurable subset of $\mathbb{R}^{N}$,
then

\begin{itemize}
\item $\mathfrak{\wp }\left( X\right) $ denotes the power set of $X$ and $%
\mathfrak{\wp }_{fin}\left( X\right) $ denotes the family of finite subsets
of $X;$

\item $|X|$ will denote the cardinality of $X$;

\item $B_{r}(x_{0})=\left\{ x\in \mathbb{R}^{N}\ |\ \left\vert
x-x_{0}\right\vert <r\right\} ;$

\item $N_{\varepsilon }(\Omega )=\left\{ x\in \mathbb{R}^{N}\ |\
dist(x,\Omega )<\varepsilon \right\} $;

\item $int(\Omega )$ denotes the interior part of $\Omega $ and $\overline{%
\Omega }$ denotes the closure of $\Omega ;$

\item $\mathfrak{F}\left( X,Y\right) $ denotes the set of all functions from 
$X$ to $Y$ and $\mathfrak{F}\left( \Omega \right) =\mathfrak{F}\left( \Omega
,\mathbb{R}\right) $;

\item if $W$ is any function space, then $W_{c}$ will denote the function
space of functions in $W$ having compact support;

\item if $W\subset \mathfrak{F}\left( \mathbb{R}^{N}\right) $ is any
function space, then $W\left( \Omega \right) $ will denote the space of
their restriction to $\Omega .$

\item if $W$ is a generic function space, its topological dual will be
denoted by $W^{\prime }$ and the pairing by $\left\langle \cdot ,\cdot
\right\rangle _{W};$

\item $C_{0}^{0}\left( \overline{\Omega }\right) $ denotes the set of
continuous functions defined on $\overline{\Omega }$ which vanish on $%
\partial \Omega ;$

\item $C^{k}\left( \Omega \right) $ denotes the set of functions defined on $%
\Omega $ which have continuous derivatives up to the order $k$;

\item $C^{k,1}\left( \Omega \right) $ denotes the set of $C^{k}\left( \Omega
\right) $-functions whose $k$-derivative is Lipschitz continuous.

\item $\mathcal{L}^{p}\left( \Omega \right) $ ($\mathcal{L}_{loc}^{p}\left(
\Omega \right) $) denotes the set of the functions $u$ such that $\left\vert
u\right\vert ^{p}$ is integrable (locally integrable) functions in $\Omega ;$

\item $L^{1}\left( \Omega \right) $ ($L_{loc}^{1}\left( \Omega \right) $)
denotes the usual equivalence classes of integrable (locally integrable)
functions in $\Omega $; these classes will be denoted by:%
\begin{equation*}
\left[ f\right] _{L^{1}}=\left\{ g\in \mathcal{L}_{loc}^{1}\left( \Omega
\right) \ |\ g(x)=f(x)\ a.e.\right\} ;
\end{equation*}

\item $H^{k,p}\left( \Omega \right) $ denotes the usual Sobolev space of
functions defined on $\Omega $;

\item $BV\left( \Omega \right) $ denotes space of functions of bounded
variation defined on $\Omega $;

\item $\mathfrak{D}\left( \Omega \right) $ denotes the set of the infinitely
differentiable functions with compact support defined on $\Omega ;\ 
\mathfrak{D}^{\prime }\left( \Omega \right) $ denotes the topological dual
of $\mathfrak{D}\left( \Omega \right) $, namely the set of distributions on $%
\Omega ;$

\item if $\Omega =\mathbb{R}^{N}$, when no ambiguity is possible, we will
write $\mathfrak{D}$, $\mathfrak{D}^{\prime }$, $L^{1},$ $\mathcal{L}%
^{1},... $ instead of $\mathfrak{D}\left( \mathbb{R}^{N}\right) $, $%
\mathfrak{D}^{\prime }\left( \mathbb{R}^{N}\right) $, $L^{1}\left( \mathbb{R}%
^{N}\right) ,$ $\mathcal{L}^{1}\left( \mathbb{R}^{N}\right) ,...$

\item $\mathbb{E}$ will denote the field of Euclidean number which will be
defined in section \ref{lt};

\item for any $\xi \in \mathbb{E}^{N},\rho \in \mathbb{E}$, we set $%
\mathfrak{B}_{\rho }(\xi )=\left\{ x\in \mathbb{E}^{N}:\left\vert x-\xi
\right\vert <\rho \right\} $;

\item $supp(f)$ denotes the usual notion of support of a function or a
distribution in $\mathbb{R}^{N}$;

\item $\mathfrak{supp}(f)=\left\{ x\in \Gamma :f(x)\neq 0\right\} $ denotes
the support of a grid function

\item $\mathfrak{mon}(x)=\{y\in \mathbb{E}^{N}:x\sim y\}$ where $x\sim y$
means that $x-y$ is infinitesimal; the set $\mathfrak{mon}(x)$ is called
monad of $x$ (see Def. \ref{MG});

\item $\mathfrak{gal}(x)=\{y\in \mathbb{E}^{N}:x-y$ is finite$\}$; the set $%
\mathfrak{gal}(x)$ is called galaxy of $x$ (see Def. \ref{MG});

\item we denote by $\chi _{X}$ the indicator (or characteristic) function of 
$X$, namely 
\begin{equation*}
\chi _{X}(x)=%
\begin{cases}
1 & if\,\,\,x\in X \\ 
0 & if\,\,\,x\notin X\,;%
\end{cases}%
\end{equation*}%
If $X=\left\{ a\right\} $ and $a$ is an atom, then, in order to simplify the
notation, we will write $\chi _{a}(x)$ instead of $\chi _{\left\{ a\right\}
}(x).$

\item $\nabla =(\partial _{1},...,\partial _{N})$ denotes the usual gradient
of standard functions;

\item $D=(D_{1},...,D_{N})$ will denote the extension of the gradient in the
sense of the ultrafunctions;

\item $\nabla \cdot \phi $ will denote the usual divergence of standard
vector fields $\phi \in \left[ C^{1}\right] ^{N}$;

\item $D\cdot \phi \ $will denote the extension of the divergence in the
sense of ultrafunctions;

\item $\Delta $ denotes the usual Laplace operator of standard functions;

\item $D^{2}$ will denote the extension of the Laplace operator in the sense
of ultrafunctions.
\end{itemize}

\section{Preliminary notions\label{PN}}

As we have already remarked in the introduction, in this section, we present
the material necessary to the rest of the paper. In particular we present an
approach to NSA based on $\Lambda $-theory. This part has been written in
such a way to be understood also by a reader who is not familiar with NSA. $%
\Lambda $-theory can be considered a different approach to Nonstandard
Analysis. It can be introduced via the notion of $\Lambda $-limit, and it
can be easily used for the purposes of this paper.

\subsection{Non Archimedean Fields\label{naf}}

Here, we recall the basic definitions and some facts regarding non
Archimedean fields.

\begin{definition}
A field ${\mathbb{K}}$ is called ordered if there is a set ${\mathbb{K}}%
^{+}\subset {\mathbb{K}}$ such that

\begin{enumerate}
\item $x,y\in {\mathbb{K}}^{+}\Rightarrow x+y,xy\in {\mathbb{K}}^{+}$

\item ${\mathbb{K=K}}^{+}\cup \left\{ 0\right\} \cup {\mathbb{K}}^{-}$ where 
${\mathbb{K}}^{-}=\left\{ x\in {\mathbb{K\ }}|\ -x\in {\mathbb{K}}%
^{+}\right\} $
\end{enumerate}
\end{definition}

In an ordered field the order relation is defined as follows:%
\begin{equation*}
x<y\Leftrightarrow y-x\in {\mathbb{K}}^{+}
\end{equation*}

In the following, ${\mathbb{K}}$ will denote an ordered field. Its elements
will be called numbers. It is well known that every ordered field contains
(a copy of) the rational numbers; hence the following definitions makes
sense:

\begin{definition}
Let $\mathbb{K}$ be an ordered field. Let $\xi \in \mathbb{K}$. We say that:

\begin{itemize}
\item $\xi $ is infinitesimal if, for all positive $n\in \mathbb{N}$, $|\xi
|<\frac{1}{n}$;

\item $\xi $ is finite if there exists $n\in \mathbb{N}$ such that $|\xi |<n$%
;

\item $\xi $ is infinite if, for all $n\in \mathbb{N}$, $|\xi |>n$
(equivalently, if $\xi $ is not finite).
\end{itemize}
\end{definition}

\begin{definition}
An ordered field $\mathbb{K}$ is called Non-Archimedean if it contains an
infinite number.
\end{definition}

It's easily seen that all infinitesimal are finite, that the inverse of an
infinite number is a nonzero infinitesimal number. Infinitesimal numbers can
be used to formalize a new notion of "closeness":

\begin{definition}
\label{def infinite closeness} We say that two numbers $\xi ,\zeta \in {%
\mathbb{K}}$ are infinitely close if $\xi -\zeta $ is infinitesimal. In this
case, we write $\xi \sim \zeta $.
\end{definition}

Clearly, the relation "$\sim $" of infinite closeness is an equivalence
relation.

\begin{theorem}
\label{na}If $\mathbb{K\supseteq R}$ is an ordered field every finite number 
$\xi \in \mathbb{K}$ is infinitely close to a unique real number $r\sim \xi $%
. $r$ is called the the \textbf{standard part} of $\xi $ and denoted by $%
st(\xi )$.
\end{theorem}

\textbf{Proof}: Given a finite number $\xi \in \mathbb{K}$, we set 
\begin{equation*}
A:=\{r\in \mathbb{R}\ |\ r<\xi \},B:=\{r\in \mathbb{R}\ |\ r\geq \xi \}
\end{equation*}%
We have that $(A,B)$ is a section of $\mathbb{R}$; moreover, since $\xi $ is
finite, $A\neq \varnothing $ and $B\neq \varnothing .$ Then, by the
completeness of the reals $\exists c\in \mathbb{R}$ 
\begin{equation*}
\forall a\in A,\forall b\in B,\ a\leq c\leq b.
\end{equation*}%
Now, it is not difficult to check that $c\sim \xi $.

$\square $

\bigskip

\begin{corollary}
\label{co1}If $\mathbb{K\supset R}$, ($\mathbb{K\neq R}$) is an ordered
field, then it is Non-Archimedean.
\end{corollary}

\textbf{Proof}: Take $\xi \in \mathbb{K}\backslash \mathbb{R}$. If $\xi $ is
infinite, then $\mathbb{K}$ is Non-Archimedean by definition. If $\xi $ is
finite then 
\begin{equation*}
\zeta :=\frac{1}{\xi -st\left( \xi \right) }
\end{equation*}%
is infinite, and hence $\mathbb{K}$ is Non-Archimedean.

$\square $

\bigskip

Now let us examine some (obvious) properties of the function $st(\cdot ).$

\begin{proposition}
\label{PS}Let $\xi $ and $\zeta $ be finite numbers, then

\begin{enumerate}
\item if $\xi \in \mathbb{R}$, $st\left( \xi \right) =\xi ;$

\item \label{d}$\xi \leq \zeta \Rightarrow st\left( \xi \right) \leq
st\left( \zeta \right) ;$

\item $st\left( \xi +\zeta \right) =st\left( \xi \right) +st\left( \zeta
\right) ;$

\item $st\left( \xi \cdot \zeta \right) =st\left( \xi \right) \cdot st\left(
\zeta \right) ;$

\item if $st\left( \zeta \right) \neq 0,$ then $st\left( \frac{\xi }{\zeta }%
\right) =\frac{st\left( \xi \right) }{st\left( \zeta \right) }.$
\end{enumerate}
\end{proposition}

\textbf{Proof:\ }The first four statements can be proved easily. In order to
prove (5), we put%
\begin{eqnarray*}
\xi &=&r+\varepsilon \\
\zeta &=&s+\theta
\end{eqnarray*}%
where $r,s\in \mathbb{E}_{\kappa }$, $\varepsilon \sim \theta \sim 0.$ Then, 
\begin{equation*}
st\left( \xi \cdot \zeta \right) =st\left[ \left( r+\varepsilon \right)
\left( s+\theta \right) \right] =st\left[ rs+\left( \varepsilon s+\theta
r+\varepsilon \theta \right) \right]
\end{equation*}%
Since $\varepsilon s+\theta r+\varepsilon \theta \sim 0$, we have that $%
st\left( \xi \cdot \zeta \right) =rs=st\left( \xi \right) \cdot st\left(
\zeta \right) .$ Finally 
\begin{equation*}
st\left( \zeta \right) \cdot st\left( \frac{\xi }{\zeta }\right) =st\left(
\zeta \cdot \frac{\xi }{\zeta }\right) =st\left( \xi \right) ;
\end{equation*}%
hence%
\begin{equation*}
st\left( \frac{\xi }{\zeta }\right) =\frac{st\left( \xi \right) }{st\left(
\zeta \right) }.
\end{equation*}

$\square $

\begin{definition}
\label{MG}Let $\mathbb{K}$ be a Non-Archimedean field, and $\xi \in \mathbb{K%
}$ a number. The monad of $\xi $ is the set of all numbers that are
infinitely close to it:%
\begin{equation*}
\mathfrak{m}\mathfrak{o}\mathfrak{n}(\xi )=\{\zeta \in \mathbb{K}:\xi \sim
\zeta \},
\end{equation*}%
and the galaxy of $\xi $ is the set of all numbers that are finitely close
to it: 
\begin{equation*}
\mathfrak{gal}(\xi )=\{\zeta \in \mathbb{K}:\xi -\zeta \ \text{is\ finite}\}.
\end{equation*}
\end{definition}

By definition, it follows that the set of infinitesimal numbers is $%
\mathfrak{mon}(0)$ and that the set of finite numbers is $\mathfrak{gal}(0)$%
. Moreover, the standard part can be regarded as a function:%
\begin{equation}
st:\mathfrak{gal}(0)\rightarrow \mathbb{R}.  \label{sh}
\end{equation}

\subsection{$\Lambda $-theory \label{lt}}

In order to construct a space of ultrafunctions it is useful to take the set 
$\Lambda $ sufficiently large; for example a superstructure over $\mathbb{R}$
defined as follows: 
\begin{equation*}
\Lambda =V_{\infty }(\mathbb{R})=\bigcup_{n\in \mathbb{N}}V_{n}(\mathbb{R}),
\end{equation*}%
where the sets $V_{n}(\mathbb{R})$ are defined by induction:%
\begin{equation*}
V_{0}(\mathbb{R})=\mathbb{R}
\end{equation*}%
and, for every $n\in \mathbb{N}$, 
\begin{equation}
V_{n+1}(\mathbb{R})=V_{n}(\mathbb{R})\cup \mathfrak{\wp }\left( V_{n}(%
\mathbb{R})\right) .  \label{Vn}
\end{equation}%
Identifying the couples with the Kuratowski pairs and the functions and the
relations with their graphs, it follows that{\ }$V_{\infty }(\mathbb{R}%
\mathbb{)}$ contains every mathematical entities used in PDE's.

Let 
\begin{equation}
\mathfrak{L}=\mathfrak{\wp }_{fin}\left( \Lambda \right)  \label{elle}
\end{equation}%
be the family of finite subsets of $\Lambda .$ $\mathfrak{L}$ equipped with
the partial order structure "$\subset $" is a directed set. A function $%
\varphi :\mathfrak{L}\rightarrow E$ will be called \textit{net }(with values
in $E$). The limit of a net is well defined: for example if $\varphi :%
\mathfrak{L}\rightarrow \mathbb{R}$, we set 
\begin{equation}
L=\lim_{\lambda \rightarrow \Lambda }\varphi (\lambda )  \label{lim+}
\end{equation}%
if and only if, $\forall \varepsilon \in \mathbb{R}^{+}$, $\exists \lambda
_{0}\in \mathfrak{L},$ such that $\forall \lambda \supset \lambda _{0},\ $%
\begin{equation*}
|\varphi (\lambda )-L|\ \leq \varepsilon
\end{equation*}%
Notice that in the notation (\ref{lim+}), $\Lambda $ can be regarded as the
"point at infinity" of $\mathfrak{L}$. A typical example of a limit of a net
defined on $\mathfrak{L}$ is provided by the definition of the Cauchy
integral:%
\begin{equation*}
\int_{a}^{b}f(x)dx=\lim_{\lambda \rightarrow \Lambda }\sum_{x\in \left[ a,b%
\right] \cap \lambda }f(x)(x^{+}-x);\ \ x^{+}=\min \left\{ y\in \mathbb{R}%
\cap \lambda \ |\ y>x\right\} .
\end{equation*}%
Now we will introduce axiomatically a new notion of limit:

\begin{axiom}
\label{EU}There is a field $\mathbb{E}\supset \mathbb{R}$, called field of
Euclidean numbers, such that every net 
\begin{equation*}
\varphi :\mathfrak{L}\rightarrow V_{n}(\mathbb{R}),\ \ n\in \mathbb{N},
\end{equation*}%
has a unique $\Lambda $-limit 
\begin{equation*}
\lim_{\lambda \uparrow \Lambda }\varphi (\lambda )\in V_{n}(\mathbb{E})
\end{equation*}%
which satisfies the following properties:

\begin{enumerate}
\item \label{EU1}\textit{if eventually }$\varphi (\lambda )=\psi (\lambda )$,%
\footnote{%
We say that a relation $\varphi (\lambda )\mathcal{R}\psi (\lambda )$ 
\textit{eventually} holds \textit{if } $\exists \lambda _{0}\in \mathfrak{L}$
such that $\forall \lambda \supset \lambda _{0},\ \varphi (\lambda )\mathcal{%
R}\psi (\lambda )$.}then 
\begin{equation*}
\lim_{\lambda \uparrow \Lambda }\varphi (\lambda )=\lim_{\lambda \uparrow
\Lambda }\psi (\lambda );
\end{equation*}

\item \label{A3}if $\varphi _{1}(\lambda ),...,\varphi _{n}(\lambda )$ are
nets, then 
\begin{equation*}
\lim_{\lambda \uparrow \Lambda }\left\{ \varphi _{1}(\lambda ),...,\varphi
_{n}(\lambda )\right\} =\left\{ \lim_{\lambda \uparrow \Lambda }\varphi
_{1}(\lambda ),...,\lim_{\lambda \uparrow \Lambda }\varphi _{n}(\lambda
)\right\} ;
\end{equation*}

\item \label{A4}if $E_{\lambda }$ is a net of sets, then 
\begin{equation*}
\lim_{\lambda \uparrow \Lambda }E_{\lambda }=\left\{ \lim_{\lambda \uparrow
\Lambda }\varphi (\lambda )\ |\ \forall \lambda \in \mathfrak{L,\ }\varphi
(\lambda )\in E_{\lambda }\right\} ;
\end{equation*}

\item \label{EU2}we have that 
\begin{equation}
\mathbb{E}:=\left\{ \lim_{\lambda \uparrow \Lambda }\ x_{\lambda }\ |\
\forall \lambda \in \mathfrak{L},\mathfrak{\ }x_{\lambda }\in \mathbb{R}%
\right\}  \label{EEE}
\end{equation}%
and \textit{if }$x_{\lambda },y_{\lambda }\in \mathbb{R}$, then, 
\begin{equation*}
\lim_{\lambda \uparrow \Lambda }\ \left( x_{\lambda }+y_{\lambda }\right)
=\lim_{\lambda \uparrow \Lambda }\ x_{\lambda }+\lim_{\lambda \uparrow
\Lambda }\ y_{\lambda },
\end{equation*}%
\begin{equation*}
\lim_{\lambda \uparrow \Lambda }\ \left( x_{\lambda }\cdot y_{\lambda
}\right) =\lim_{\lambda \uparrow \Lambda }\ x_{\lambda }\cdot \lim_{\lambda
\uparrow \Lambda }\ y_{\lambda };
\end{equation*}%
\begin{equation}
x_{\lambda }\geq y_{\lambda }\Rightarrow \lim_{\lambda \uparrow \Lambda }\
x_{\lambda }\geq \lim_{\lambda \uparrow \Lambda }\ y_{\lambda }.
\label{pippa}
\end{equation}
\end{enumerate}
\end{axiom}

Notice that in order to distinguish the limit (\ref{lim}) (which we will
call \textit{Cauchy limit}) from the $\Lambda $-limit, we have used the
symbols "$\lambda \rightarrow \Lambda $" and "$\lambda \uparrow \Lambda $"
respectively.

In the rest of this section, we will make some remarks for the readers who
are not familiar with NSA. The points 1,2,4 are not surprising since we
expect them to be satisfied by any notion of limit provided the target space
be equipped with a reasonable topology. The point 3 can be considered as a
definition. In axiom \ref{EU}, the new (and, for someone, surprising) fact
is that every net has a $\Lambda $-limit. Nevertheless Axiom \ref{EU} is not
contradictory and a model for it can be constructed in ZFC (see section \ref%
{cen} or \cite{ultra} or \cite{BL2021} for further details).

Probably the first question which a newcomer to the world of NSA would ask
is the following: what is the limit of a sequence such that $\varphi
(\lambda ):=(-1)^{\left\vert \lambda \right\vert }$ since it takes the
values $+1$ if $\left\vert \lambda \right\vert $ is even or $-1$ if $%
\left\vert \lambda \right\vert $ is odd. Let us see what Axiom \ref{EU}
tells us. By \ref{EU}.\ref{A3} and \ref{EU}.\ref{EU1},

\begin{equation*}
\lim_{\lambda \uparrow \Lambda }\left\{ 1,-1\right\} =\left\{ \lim_{\lambda
\uparrow \Lambda }1,\ \lim_{\lambda \uparrow \Lambda }\left( -1\right)
\right\} =\left\{ 1,-1\right\}
\end{equation*}%
and by \ref{EU}.\ref{A4}, we have that 
\begin{equation*}
\lim_{\lambda \uparrow \Lambda }(-1)^{\left\vert \lambda \right\vert }\in
\lim_{\lambda \uparrow \Lambda }\left\{ 1,-1\right\} ;
\end{equation*}%
and hence 
\begin{equation*}
\lim_{\lambda \uparrow \Lambda }(-1)^{\left\vert \lambda \right\vert }\in
\left\{ 1,-1\right\}
\end{equation*}%
Then either $\lim_{\lambda \uparrow \Lambda }(-1)^{\left\vert \lambda
\right\vert }=1$ or $\lim_{\lambda \uparrow \Lambda }(-1)^{\left\vert
\lambda \right\vert }=-1.$ Which alternative occurs cannot be deduced by
axiom \ref{EU}; each alternative can be added as an independent axiom. In
the models constructed in \cite{BL2021} this limit is +1; however this and
similar questions are not relevant for this paper and we refer to the
mentioned references for a deeper discussion of this point (in particular
see \cite{BDN2018} and \cite{BL2021}). Axiom \ref{EU} is sufficient for our
applications.

The second question which a newcomer would ask is about the limit of a
divergent sequence such that $\varphi (\lambda ):=\left\vert \lambda \cap 
\mathbb{N}\right\vert $. Let us put%
\begin{equation}
\alpha :=\lim_{\lambda \uparrow \Lambda }\ \left\vert \lambda \cap \mathbb{N}%
\right\vert  \label{alfa}
\end{equation}%
What can we say about $\alpha $? By (\ref{pippa}), $\alpha \notin \mathbb{R}%
. $ In order to give a feeling of the "meaning" of $\alpha ,$ we will put it
in relation with other infinite numbers. If $E\in \Lambda \backslash \mathbb{%
R},$ we put%
\begin{equation}
\mathfrak{num}\left( E\right) =\lim_{\lambda \uparrow \Lambda }\ |E\cap
\lambda |.  \label{num}
\end{equation}%
If $E$ is a finite set, the sequence is eventually equal to the number of
elements of $E$; then, by axiom \ref{EU}.\ref{EU1}, 
\begin{equation*}
\mathfrak{num}\left( E\right) =\left\vert E\right\vert \in \mathbb{N}.
\end{equation*}%
If $E$ is an infinite set, $\mathfrak{num}\left( E\right) \notin \mathbb{N}$%
. Hence, the limits like (\ref{num}) give mathematical entities that extend
the notion of "\textit{number of elements of a set}" to infinite sets and it
is legitimate to call them "\textit{infinite numbers}". The infinite number $%
\mathfrak{num}\left( E\right) $ is called \textit{numerosity} of $E$. The
theory of numerosities can be considered as an extension of the Cantorian
theory of cardinal and ordinal numbers. The reader interested to the details
and the developments of this theory is referred to \cite%
{benci95b,BDN2003,BDNF1,BL2021,BF}.

If a real net $x_{\lambda }$ admits the Cauchy limit, the relation between
the two limits is given by the following identity:%
\begin{equation}
\lim_{\lambda \rightarrow \Lambda }\ x_{\lambda }=st\left( \lim_{\lambda
\uparrow \Lambda }\ x_{\lambda }\right)  \label{bn}
\end{equation}

An other important relation between the two limits is the following:

\begin{proposition}
\label{lim}If 
\begin{equation*}
\lim_{\lambda \uparrow \Lambda }\ x_{\lambda }=\xi \in \mathbb{E}
\end{equation*}%
and $\xi $ is bounded, then there exists a sequence $\lambda _{n}\in 
\mathfrak{L}$ such that%
\begin{equation*}
\lim_{n\rightarrow \infty }\ x_{\lambda _{n}}=st(\xi ).
\end{equation*}
\end{proposition}

\textbf{Proof: }Set $x_{0}=st(\xi )$ and for every $n\in \mathbb{N}$, take $%
\lambda _{n}$ such that $\ x_{\lambda _{n}}\in B_{1/n}(x).$

$\square $

\begin{remark}
As we have already remarked, the field of Euclidean numbers is a hyperreal
field in the sense of Non Standard Analysis. We do not use the name
"hyperreal numers" to emphasize the fact that $\mathbb{E}$ has been defined
by the notion of $\Lambda $-limit and hence it satisfies some properties
which are not shared by other hyperreal fields. These properties are
relevant in the definitions of ultrafunctions. The explanation of the choice
of the name "Euclidean numbers" can be found in \cite{BF}.
\end{remark}

\subsection{Extension of sets and functions}

In this section we recall some basic notions of Non Standard Analysis
presented in the framework of $\Lambda $-theory. Given a set $A\in \Lambda $%
, we define 
\begin{equation}
A^{\ast }=\left\{ \lim_{\lambda \uparrow \Lambda }\varphi (\lambda )\ |\
\forall \lambda ,\ \varphi (\lambda )\in A\right\} ;  \label{Astar}
\end{equation}%
Following Keisler \cite{keisler76}, $A^{\ast }$ will be called the \textbf{%
natural extension} of $A.$ By (\ref{EEE}), we have that $\mathbb{R}^{\ast }=%
\mathbb{E}$. If we identify a relation $\mathcal{R}$ or a function $f$ with
its graph, then, by (\ref{Astar}) $\mathcal{R}^{\ast }$ and $f^{\ast }$ are
well defined.

In particular any function 
\begin{equation*}
f:A\rightarrow B,\ \ A,B\in \Lambda ,
\end{equation*}%
can be extended to $A^{\ast }$ and we have that%
\begin{equation}
f^{\ast }\left( \lim_{\lambda \uparrow \Lambda }\ x_{\lambda }\right)
=\lim_{\lambda \uparrow \Lambda }~f\left( x_{\lambda }\right) ;  \label{42}
\end{equation}%
the function 
\begin{equation*}
f^{\ast }:A^{\ast }\rightarrow B^{\ast },
\end{equation*}%
will be called \textbf{natural extension} of $f.$ More in general, if%
\begin{equation*}
u_{\lambda }:A\rightarrow B
\end{equation*}%
is a net of functions, we have that for any $x=\lim_{\lambda \uparrow
\Lambda }x_{\lambda },\ x_{\lambda }\in B$, 
\begin{equation}
u(x)=\lim_{\lambda \uparrow \Lambda }~u_{\lambda }\left( x_{\lambda }\right)
\label{10}
\end{equation}%
is a function from $A^{\ast }\ $to $B^{\ast }.$

\textbf{Example}: Let $f:\mathbb{R}^{N}\rightarrow \mathbb{R}$ be a
differentiable function, then%
\begin{equation*}
\partial _{i}^{\ast }f^{\ast }=\left( \partial _{i}f\right) ^{\ast }
\end{equation*}%
where the operator%
\begin{equation*}
\partial _{i}=\frac{\partial }{\partial x_{i}}:C^{1}\left( \mathbb{R}%
^{N}\right) \rightarrow C^{0}\left( \mathbb{R}\right)
\end{equation*}%
is regarded as a function between functional spaces and hence 
\begin{equation*}
\partial _{i}^{\ast }:C^{1}\left( \mathbb{R}^{N}\right) ^{\ast }\rightarrow
C^{0}\left( \mathbb{R}\right) ^{\ast }.
\end{equation*}

Following the current literature in NSA, we give the following definition:

\begin{definition}
A set $E_{\Lambda }$ obtained as $\Lambda $-limit of a net of sets $%
E_{\lambda }\in \Lambda $ is called \textbf{internal}.
\end{definition}

In particular, if $E\in \Lambda ,$ if you compare Axiom \ref{EU}.\ref{A4}
with (\ref{Astar}), then you see that 
\begin{equation*}
E^{\ast }=\lim_{\lambda \uparrow \Lambda }~E_{\lambda }
\end{equation*}%
in the case in which $E_{\lambda }$ is the net identically equal to $E$.

Let us see an example of \textbf{external} set i.e. of a set which is not
internal. By axiom \ref{EU}.\ref{A4}, the set $E^{\ast }$ contains a unique
"copy" $x^{\ast }$ of every element $x\in E$. Now set 
\begin{equation}
E^{\sigma }:=\left\{ x^{\ast }\in E^{\ast }\ |\ x\in E\right\} .
\label{sigma+}
\end{equation}%
We have that $E^{\sigma }\subseteq E^{\ast }$ and the equality holds if and
only if $E$ is finite. It is easy to see that if $E$ is infinite $E$ and $%
E^{\sigma }$ are external.

\bigskip

\textbf{Example}: Let $E=C^{0}\left( \mathbb{R}\right) ;$ then $C^{0}\left( 
\mathbb{R}\right) ^{\sigma }\subset C^{0}\left( \mathbb{R}\right) ^{\ast }.$
If we take%
\begin{equation*}
\sin ^{\ast }\left( x\right) =\lim_{\lambda \uparrow \Lambda }\ \sin \left(
x_{\lambda }\right) ,\ \ x=\lim_{\lambda \uparrow \Lambda }x_{\lambda }
\end{equation*}%
and, using the notation (\ref{alfa}), 
\begin{equation*}
\sin ^{\ast }\left( \alpha x\right) =\lim_{\lambda \uparrow \Lambda }\ \sin
(\left\vert \lambda \cap \mathbb{N}\right\vert \cdot x_{\lambda }),
\end{equation*}%
we have that $\sin ^{\ast }\left( x\right) \in C^{0}\left( \mathbb{R}\right)
^{\sigma }\subset C^{0}\left( \mathbb{R}\right) ^{\ast }\ $and $\sin ^{\ast
}\left( \alpha x\right) \in C^{0}\left( \mathbb{R}\right) ^{\ast }\backslash 
$ $C^{0}\left( \mathbb{R}\right) ^{\sigma }.$

\subsection{Hyperfinite sets\label{hs}}

An other fundamental notion in NSA is the following:

\begin{definition}
We say that a set $F\in \Lambda $ is \textbf{hyperfinite} if there is a net $%
\left\{ F_{\lambda }\right\} _{\lambda \in \Lambda }$ of finite sets such
that 
\begin{equation*}
F=\lim_{\lambda \uparrow \Lambda }~F_{\lambda }=\left\{ \lim_{\lambda
\uparrow \Lambda }\ x_{\lambda }\ |\ x_{\lambda }\in F_{\lambda }\right\}
\end{equation*}
\end{definition}

The hyperfinite sets share many properties of finite sets. For example, a
hyperfinite set $F\subset \mathbb{E}$ has a maximum $x_{M}$ and a minimum $%
x_{m}$ respectively given by%
\begin{equation*}
x_{M}=\lim_{\lambda \uparrow \Lambda }\max F_{\lambda };\ \
x_{m}=\lim_{\lambda \uparrow \Lambda }\min F_{\lambda }
\end{equation*}

Moreover, it is possible to "add" the elements of an hyperfinite set of
numbers. If $F$ is an hyperfinite set of numbers, the \textbf{hyperfinite sum%
} of the elements of $F$ is defined as follows: 
\begin{equation*}
\sum_{x\in F}x=\ \lim_{\lambda \uparrow \Lambda }\sum_{x\in F_{\lambda }}x.
\end{equation*}

One of the advantage to use the field of Euclidean numbers rather than a
generic hyperreal field lies in the possibility to associate a unique
hyperfinite set $E^{\circledcirc }$ to any set $E\in V_{\infty }(\mathbb{R})$
according to the following definition:

\begin{definition}
Given a set $E\subset \Lambda ,$ the set 
\begin{equation*}
E^{\circledcirc }:=\lim_{\lambda \uparrow \Lambda }\ \left( E\cap \lambda
\right)
\end{equation*}%
is called \textbf{hyperfinite extension} of $E.$
\end{definition}

To any set $E$ we can associate the sets $E^{\sigma },$ $E^{\circledcirc }$
and $E^{\ast }$ which are ordered as follows:%
\begin{equation*}
E^{\sigma }\subseteq E^{\circledcirc }\subseteq E^{\ast };
\end{equation*}%
the inclusions are strict if and only if the set $E$ is infinite. If $%
F=\lim_{\lambda \uparrow \Lambda }~F_{\lambda }$ is a hyperfinite set, its 
\textbf{hypercardinality} is defined by%
\begin{equation*}
\left\vert F\right\vert ^{\ast }:=\lim_{\lambda \uparrow \Lambda
}~\left\vert F_{\lambda }\right\vert
\end{equation*}%
Notice that the hypercardinality of $E^{\circledcirc },$ defined by%
\begin{equation*}
\left\vert E^{\circledcirc }\right\vert ^{\ast }=\lim_{\lambda \uparrow
\Lambda }\ \left\vert E\cap \lambda \right\vert ,
\end{equation*}%
is the numerosity of $E$ as it has been defined by (\ref{num}).

\subsection{Grid functions}

\begin{definition}
A hyperfinite set $\Gamma $ such that $\mathbb{R}^{N}\subset \Gamma \subset 
\mathbb{E}^{N}$ is called \textbf{hyperfinite grid}.
\end{definition}

For example the set $\left( \mathbb{R}^{N}\right) ^{\circledcirc }$ is an
hyperfinite grid.

\begin{definition}
A space of grid functions is a family $\mathfrak{G}(\Gamma )$ of internal
functions 
\begin{equation*}
u:\Gamma \rightarrow \mathbb{R}
\end{equation*}
\end{definition}

If $w\in \mathfrak{F}(\mathbb{R}^{N})^{\ast },$ the restriction of $w$ to $%
\Gamma $ is a grid function which we will denote by $w%
{{}^\circ}%
$ namely, if $w=\ \lim_{\lambda \uparrow \Lambda }w_{\lambda }$ and $x=\
\lim_{\lambda \uparrow \Lambda }x_{\lambda }\in \Gamma $, we have that 
\begin{equation}
w%
{{}^\circ}%
(x)=\ \lim_{\lambda \uparrow \Lambda }\ w_{\lambda }(x_{\lambda }).
\label{giusi}
\end{equation}%
For every $a\in \Gamma ,$ 
\begin{equation*}
\chi _{a}(x)\in \mathfrak{G}(\Gamma )
\end{equation*}%
is a grid function, and hence every grid function can be represented by the
following sum:%
\begin{equation}
u(x)=\sum_{a\in \Gamma }u(a)\chi _{a}(x)  \label{pu}
\end{equation}%
namely $\left\{ \chi _{a}\right\} _{a\in \Gamma }$ is a basis for $\mathfrak{%
G}(\Gamma )$ considered as a vector space over $\mathbb{E}$. Given $f\in 
\mathfrak{F}(\mathbb{R}^{N}),$ we will write $f%
{{}^\circ}%
$ instead of $(f^{\ast })%
{{}^\circ}%
,$ namely 
\begin{equation}
f%
{{}^\circ}%
(x):=\lim_{\lambda \uparrow \Lambda }\ f(x_{\lambda })=\sum_{a\in \Gamma
}f^{\ast }(a)\chi _{a}(x).  \label{lina}
\end{equation}%
so, $\mathfrak{G}(\Gamma )$ contains a unique copy $f%
{{}^\circ}%
$ of every function $f\in \mathfrak{F}(\mathbb{R}^{N}).$ If a function, such
as $1/|x|$ is not defined in some point and $x\in \Gamma $, we put $(1/|x|)%
{{}^\circ}%
$ equal to $0$ for $x=0;$ in general, if $\Omega $ is a subset of $\mathbb{R}%
^{N}$ and $f$ is defined in $\Omega $, we set%
\begin{equation}
f%
{{}^\circ}%
(x)=\sum_{a\in \Omega 
{{}^\circ}%
}f^{\ast }(a)\chi _{a}(x)  \label{zero}
\end{equation}%
where for every set $E\subset \mathbb{R}^{N}$, we define%
\begin{equation}
E%
{{}^\circ}%
=E^{\ast }\cap \Gamma .  \label{omegat}
\end{equation}

\section{Ultrafunctions\label{U}}

If we have a differential equation, it is relatively easy to find an
approximated solution in a suitable space of grid functions. If this
equation has a "classic" solution, this solution, in some sense,
approximates the classic solution. Then, if we take a grid, the "grid
solutions" is almost equal to the classic solution. However the "grid
solutions" cannot be considered as generalizations of the classic solutions
since \textit{they do not coincide with them}. The theory of ultrafunctions
is based on the idea of a space of functions (defined on a hyperfinite grid)
in which the generalized derivative and the generalized integral coincide
with the usual ones for every function $f$ in $C^{1}$ and in $C_{c}^{0}$
respectively. This fact implies that a "ultrafunction solution" coincides
with the classical one if the latter exists and hence it is legitimate to be
considered a generalized solution.

\subsection{Definition of ultrafunctions}

Let $V=V\left( \mathbb{R}^{N}\right) $ be a function space such that

\begin{equation*}
\mathfrak{D}\left( \mathbb{R}^{N}\right) \subset V\subset \mathcal{L}%
_{loc}^{1}\left( \mathbb{R}^{N}\right)
\end{equation*}%
and let $\left\{ V_{\lambda }\right\} _{\lambda \in \mathfrak{L}}$ be a net
of finite dimensional subspaces of $V$ such that%
\begin{equation*}
\dbigcup\limits_{\lambda \in \mathfrak{L}}V_{\lambda }=V.
\end{equation*}%
Now we set%
\begin{equation*}
V_{\Lambda }=V_{\Lambda }\left( \mathbb{E}^{N}\right) =\lim_{\lambda
\uparrow \Lambda }\ V_{\lambda }=\left\{ \lim_{\lambda \uparrow \Lambda }\
u_{\lambda }\ |\ u_{\lambda }\in V_{\lambda }\right\} ;
\end{equation*}%
$V_{\Lambda }$ is an internal vector space of hyperfinite dimension. Clearly 
$V_{\Lambda }\subset V^{\ast }$ since%
\begin{equation*}
V^{\ast }=\left\{ \lim_{\lambda \uparrow \Lambda }\ u_{\lambda }\ |\
u_{\lambda }\in V\right\} .
\end{equation*}

The space $V_{\Lambda }$ allows to equip a space of grid functions $%
\mathfrak{G}(\Gamma )$ of a richer structure:

\begin{definition}
\label{lella}A space of ultrafunctions $V%
{{}^\circ}%
(\Gamma )$ modelled on $V_{\Lambda }$ is a family of grid functions $%
\mathfrak{G}(\Gamma )$ such that the restriction map 
\begin{equation}
{{}^\circ}%
:V_{\Lambda }\rightarrow V%
{{}^\circ}%
(\Gamma )  \label{RO}
\end{equation}%
is an internal isomorphism between hyperfinite dimensional vector spaces.
\end{definition}

So, $u\in V%
{{}^\circ}%
(\Gamma )\ $if and only if there exists$\ $a net $u_{\lambda }\in V_{\lambda
}\left( \mathbb{R}^{N}\right) $ such that%
\begin{equation*}
u=\left( \lim_{\lambda \uparrow \Lambda }\ u_{\lambda }\right) ^{\circ }
\end{equation*}
In the following of this paper, if $u\in V%
{{}^\circ}%
(\Gamma ),$ such a net will be denoted by $u_{\lambda }$. We will denote by $%
\sigma _{a}(x)$ the only function in $V_{\Lambda }\left( \mathbb{E}%
^{N}\right) $ such that 
\begin{equation}
\sigma _{a}^{\circ }=\chi _{a}  \label{sigma}
\end{equation}%
Clearly $\left\{ \sigma _{a}(x)\right\} _{a\in \Gamma }$ is a basis of $%
V_{\Lambda }\left( \mathbb{E}^{N}\right) $ which will be called $\sigma $%
-basis. The $\sigma _{a}$'s allow to write the inverse of the map (\ref{RO})%
\begin{equation}
\left( \cdot \right) _{\Lambda }:V%
{{}^\circ}%
(\Gamma )\rightarrow V_{\Lambda }\left( \mathbb{E}^{N}\right)  \label{lambda}
\end{equation}%
as follows: if $u\in V%
{{}^\circ}%
(\Gamma ),$ 
\begin{equation*}
u_{\Lambda }(x):=\sum_{a\in \Gamma }u(a)\sigma _{a}(x).
\end{equation*}%
If $f\in \mathfrak{F}\left( \mathbb{E}^{N}\right) ,$ in order to simplify
the notation, we will write $f_{\Lambda }$ instead of $\left( f%
{{}^\circ}%
\right) _{\Lambda },$ namely we have that%
\begin{equation*}
f_{\Lambda }(x)=\sum_{a\in \Gamma }f^{\ast }(a)\sigma _{a}(x).
\end{equation*}%
Notice that%
\begin{equation}
f_{\Lambda }=f^{\ast }\Leftrightarrow f\in V\left( \mathbb{R}^{N}\right)
\label{ciccia2}
\end{equation}%
More in general, if $w\in \mathfrak{F}\left( \mathbb{E}^{N}\right) ^{\ast }$%
, we will write 
\begin{equation}
w_{\Lambda }(x)=\sum_{a\in \Gamma }w(a)\sigma _{a}(x),  \label{ciccia1}
\end{equation}%
In this case, the map 
\begin{equation}
\left( \cdot \right) _{\Lambda }:\mathfrak{F}\left( \mathbb{R}^{N}\right)
^{\ast }\rightarrow V_{\Lambda }\left( \mathbb{E}^{N}\right)  \label{lamda+}
\end{equation}%
is just a projection.

\bigskip

If $u\in V%
{{}^\circ}%
\left( \Gamma \right) \ $and $u_{\Lambda }\in \mathcal{L}^{1}\left( \mathbb{R%
}^{N}\right) ^{\ast }$, the integral can be defined as follows:%
\begin{equation}
\doint u(x)dx:=\int^{\ast }u_{\Lambda }(x)dx=\lim_{\lambda \uparrow \Lambda
}\int u_{\lambda }(x)dx.  \label{int1}
\end{equation}
We will refer to 
\begin{equation*}
\doint :V%
{{}^\circ}%
\left( \Gamma \right) \rightarrow \mathbb{E}
\end{equation*}%
as to the \textbf{pointwise integral. }The reason of this name is due to the
fact that (\ref{pu}) and (\ref{int1}) imply that%
\begin{equation}
\doint u(x)dx=\sum_{a\in \Gamma }u(a)d(a)  \label{int3}
\end{equation}%
where%
\begin{equation}
d(a):=\doint \chi _{a}(x)dx=\int^{\ast }\sigma _{a}(x)dx.  \label{int4}
\end{equation}%
We may think of $d(a)$ as the "measure" of the point $a\in \Gamma $. The
pointwise integral extends the usual Lebesgue integral from $V$ to $V%
{{}^\circ}%
$, more exactly, if $f\in V\cap \mathcal{L}^{1},$ then%
\begin{equation}
\doint f%
{{}^\circ}%
(x)dx=\int f(x)dx  \label{mara}
\end{equation}%
However the equality above is not true for every Lebesgue integrable
function. In fact, if $a\in \mathbb{R}^{N}$,%
\begin{equation*}
\int \chi _{a}(x)dx=0
\end{equation*}%
but, by (\ref{int3}), we have that 
\begin{equation*}
\doint \chi _{a}^{\circ }(x)dx>0
\end{equation*}%
at least for some $a\in \mathbb{R}^{N}.$ This fact is quite natural, in fact
when we work in a non-Archimedean world infinitesimals matter and cannot be
forgotten as the Riemann and the Lebesgue integrals do. Also the above
inequality shows that it is necessary to use a different symbol to
distinguish the pointwise integral from the Lebesgue integral (here we have
used $\doint $).

Given $u\in V%
{{}^\circ}%
,\ $if $u_{\lambda }\in C^{1}\cap V_{\lambda },$ and $\partial
_{i}u_{\lambda }\in V_{\lambda },$ it is natural to define the\textbf{\
partial derivative} in a point $x=\ \lim_{\lambda \uparrow \Lambda
}x_{\lambda }\in \Gamma $, as follows%
\begin{equation}
D_{i}u(x):=\left[ \partial _{i}^{\ast }u_{\Lambda }(x)\right] 
{{}^\circ}%
=\lim_{\lambda \uparrow \Lambda }\ \partial _{i}u_{\lambda }(x_{\lambda }).
\label{der}
\end{equation}%
So, if $f\in $ $C^{1}\cap V,\ $and $\partial f\in $ $C^{0}\cap V,$ we have 
\begin{equation*}
D_{i}f%
{{}^\circ}%
(x)=\left[ \partial _{i}^{\ast }f^{\ast }(x)\right] 
{{}^\circ}%
\end{equation*}%
and hence, if $x\in \mathbb{R}$, $D_{i}f%
{{}^\circ}%
(x)=\partial _{i}f(x).$

In particular, if we choose $V=C_{c}^{1},$ we have that the integral and the
derivative are defined for every ultrafunction in $V%
{{}^\circ}%
$.

\subsection{Epilogic functions and the space $V$\label{EF}}

There are many ultrafunction spaces which depend on the choice of the space $%
V,$ the net $\left\{ V_{\lambda }\right\} $ and the grid $\Gamma .$ However
there are some basic properties which should be satisfied by a "good" space
of ultrafunctions which make the theory rich and flexible. In order to get
such a space, the first step consists in choosing an appropriate space $V$.
The simplest choice is to take $V=C_{0}^{1}$ so that the pointwise integral
and the generalized derivative be well defined for every ultrafunctions just
using (\ref{int}) and (\ref{der}). However this choice is too restrictive
for many applications. In fact it is useful to work with the characteristic
function $\chi _{\Omega }$ of an open set at least when $\partial \Omega $
is sufficiently smooth. Then, it seem reasonable to work in $BV\cap
L_{c}^{\infty }.$ Unfortunately this space is not suitable since a function $%
f\in L^{\infty }$ is not pointwise defined and hence tha map $\left( 
{{}^\circ}%
\right) $ is not well defined. However we can overcome this difficulty
taking a space isomorphic to $L^{\infty }$ by choosing one function in each
equivalence class of $L^{\infty }$. This choice must be done in a way
consistent with the main operations in $L^{\infty }$. Let us see how to do
it.

We start recalling the following standard terminology: for every function $%
f\in \mathcal{L}_{loc}^{1}$ we say that a point $x\in \mathbb{R}^{N}$ is a
Lebesgue point for $f$ if%
\begin{equation*}
f(x)=\lim_{r\rightarrow 0^{+}}\frac{1}{m(B_{r}(x))}\int_{B_{r}(x)}f(y)dy,
\end{equation*}%
where $m(B_{r}(x))$ is the Lebesgue measure of the ball $B_{r}(x);$ we
recall the very important Lebesgue theorem (see e.g. \cite{lebesgue}), that
we will need in the following:

\begin{theorem}
\label{Lebesgue}If $f\in \mathcal{L}_{loc}^{1}(\mathbb{R}^{N})$ then a.e. $%
x\in \mathbb{R}^{N}$ is a Lebesgue point for $f$.
\end{theorem}

Now, we fix once for ever an infinitesimal number $\eta >0$; for example we
can choose $\eta =\alpha ^{-1}$ (see (\ref{alfa})). Given a function $f\in 
\mathcal{L}_{loc}^{\infty },$ we set%
\begin{equation}
\overline{f}(x)=st\left( \frac{1}{m(B_{\eta }(x))}\int_{B_{\eta
}(x)}f(y)dy\right) ,  \label{due+}
\end{equation}%
where $m(B_{\eta }(x))$ is the Lebesgue measure of the ball $B_{\eta }(x).$
By this trick we can choose a unique function $\overline{f}$ in the
equivalence class $\left[ f\right] _{L^{\infty }}\in L^{\infty }.$

\begin{lemma}
\label{aa}The operator $f\mapsto \overline{f}$ satisfies the following
properties:

\begin{enumerate}
\item \label{a1+}if $x$ is a Lebesgue point for $f$ then $\overline{f}%
(x)=f(x);$

\item \label{a2}$f(x)=\overline{f}(x)\ $a.e.$;$

\item \label{a3}if $f(x)=g(x)$ a.e. then $\overline{f}(x)=\overline{g}(x);$

\item \label{a4}$\overline{\overline{f}}(x)=\overline{f}(x).$

\item \label{a5}$\overline{f+g}=\overline{f}+\overline{g}$.
\end{enumerate}
\end{lemma}

\textbf{Proof}: (\ref{a1+}) - If $x$ is a Lebesgue point for $f$ then%
\begin{equation*}
\frac{1}{m(B_{\eta }(x))}\int_{B_{\eta }(x)}f(y)dy\sim f(x),
\end{equation*}%
so $\overline{f}(x)=f(x).$

(\ref{a2}) - This follows immediately by Theorem \ref{Lebesgue} and (\ref%
{a1+}).

(\ref{a3}) - Let $x\in \mathbb{R}^{N}.$ Since $f(x)=g(x)$ a.e., we obtain
that $\int_{B_{\eta }(x)}f(y)dy=\int_{B_{\eta }(x)}g(y)dy,$ so 
\begin{eqnarray*}
\overline{f}(x) &=&st\left( \frac{1}{m(B_{\eta }(x))}\int_{B_{\eta
}(x)}f(y)dy\right) \\
&=&st\left( \frac{1}{m(B_{\eta }(x))}\int_{B_{\eta }(x)}g(y)dy\right) =%
\overline{g}(x).
\end{eqnarray*}

(\ref{a4}) - This follows from (\ref{a2}) and (\ref{a3}).

(\ref{a5}) - Follows from the linearity of the function $st.$

$\square $

We now define the space of \textbf{epilogic functions} as follows: 
\begin{equation*}
EPL:=\left\{ u\in \mathcal{L}_{loc}^{\infty }\ |\ \overline{u}%
(x)=u(x)\right\}
\end{equation*}%
This name comes from the greek $\varepsilon \pi \iota \lambda o\gamma \eta $
= "choice" since each function in $EPL$ has been chosen in an equivalence
class of $L_{loc}^{\infty }.$

We list some properties of $EPL$ that will be useful in the following:

\begin{proposition}
\label{diana}The following properties hold:

\begin{enumerate}
\item \label{b3+}if $u,v\in EPL$ then $u=v\ $a.e. if and only if $u=v$;

\item \label{b3}$EPL$ is a module over the ring $C^{0},$ namely, $\varphi
\in C^{0}$ and $f\in EPL$ implies $\varphi f\in EPL$;

\item \label{b4}the $L^{2}$ norm is a norm for $EPL\cap \mathcal{L}^{2}$
(and not a pseudonorm).
\end{enumerate}
\end{proposition}

\textbf{Proof}: (\ref{b3+}) - Let $u,v\in EPL$. If $u=v$ then clearly $\bar{u%
}=\bar{v}$ a.e.; conversely, let us suppose that $u=v$ a.e.; by Lemma \ref%
{aa}, (\ref{a3}) we deduce that $\overline{u}(x)=\overline{v}(x).$ But $%
u,v\in EPL,$ so $u(x)=\overline{u}(x)=\overline{v}(x)=v(x).$

(\ref{b3}) - if $f(x)\in EPL,$ then for every $\varphi $ in $C^{k},$ 
\begin{eqnarray*}
\overline{\varphi (x)f(x)} &=&st\left( \frac{1}{m(B_{\eta }(x))}%
\int_{B_{\eta }(x)}\varphi (y)f(y)dy\right) \\
&\sim &st\left( \frac{1}{m(B_{\eta }(x))}\int_{B_{\eta }(x)}\varphi
(x)f(y)dy\right) \\
&=&\varphi (x)st\left( \frac{1}{m(B_{\eta }(x))}\int_{B_{\eta }(x)\varphi
(x)}f(y)dy\right) \\
&=&\varphi (x)\overline{f(x)}=\varphi (x)f(x).
\end{eqnarray*}

(\ref{b4}) - Let $u\in EPL$ be such that $\left\Vert u\right\Vert
_{L^{2}}=0. $ Then $u=0$ a.e. and since $0\in EPL,$ by (\ref{b3+}), we
deduce that $u=0.$

$\square $

\textbf{Example}: If $\Omega \subset \mathbb{R}^{N}$ is a measurable set,
the \textbf{density function} of $\Omega $ is defined as follows: 
\begin{equation*}
\Theta _{\Omega }(x)=\lim_{r\rightarrow 0^{+}}\frac{m(B_{r}(x)\cap \Omega )}{%
m(B_{r}(x))},
\end{equation*}%
Hence, by the Lebesgue theorem we have that $\Theta _{\Omega }(x)$ is
defined $a.e.$ Using the operator (\ref{due+}), we can define $\Theta
_{\Omega }$ in every point by setting 
\begin{equation}
\Theta _{\Omega }(x)=\bar{\chi}_{\Omega }(x).  \label{estate}
\end{equation}%
Clearly $\bar{\chi}_{\Omega }(x)$ is a function whose value is 1 in $%
int(\Omega )$ and 0 in $\mathbb{R}^{N}\setminus \overline{\Omega }.$\ If $%
\Omega $ is an open set with smooth boundary, we have that $\forall x\in 
\mathbb{R}^{N}$%
\begin{equation}
\bar{\chi}_{\Omega }(x)=\left\{ 
\begin{array}{cc}
1 & \text{if\ }\ x\in \Omega ; \\ 
0 & \text{if\ }\ x\notin \Omega ; \\ 
\frac{1}{2} & \text{if\ }\ x=\partial \Omega .%
\end{array}%
\right.  \label{chi}
\end{equation}

The next ingredient necessary to define $V$ is the space of the function of
bounded variation $BV$. We recall that $f\in BV$ if $f\in \mathcal{L}^{1}$
and its derivative $\partial _{i}f$ (in the sense of distributions) is a
Radon measure, namely, for every continuous function $\varphi $, the
functional $\varphi \mapsto \left\langle \partial _{i}f,\varphi
\right\rangle $ is well defined; it is well known that this measure can be
estended to every Borellian function and, with some abuse of notation, we
will write%
\begin{equation}
\int g(x)\partial _{i}f\ dx  \label{nota}
\end{equation}%
rather than $\int g(x)\ d\left( \partial _{i}f\right) $, since for $\partial
_{i}f\in L^{1}$, $\partial _{i}f$ coincides with a measure density and hence%
\begin{equation*}
\int g(x)\ d\left( \partial _{i}f\right) =\int g(x)\partial _{i}f(x)\ dx
\end{equation*}

Finally, we can define the space $V$ by setting%
\begin{equation}
V=BV\cap EPL\cap \mathcal{R}  \label{VV}
\end{equation}%
where $\mathcal{R}$ is the space of Riemann integrable function (we recall
that in the usual definition the function in $\mathcal{R}$ have compact
support). From now on $V$ will denote the space (\ref{VV}) and $\partial
_{i}f$ will denote the $BV$-derivative of $f$.

The space $V$ is suitable for our purposes; we have that 
\begin{equation*}
C_{c}^{0,1}\subset V
\end{equation*}%
and that $\bar{\chi}_{\Omega }\in V$ provided that $\chi _{\Omega }\in 
\mathcal{R\ }$and $\Omega $ is a Caccioppoli set, namely $\chi _{\Omega }\in
BV$. Moreove, by Prop. \ref{diana}.\ref{b3}, $V$ is a $C^{0,1}$ module:%
\begin{equation*}
\varphi \in C^{0,1},\ f\in V\Rightarrow \varphi f\in V.
\end{equation*}

\subsection{Definition of fine ultrafunctions\label{AF}}

Roughly speaking a space of ultrafunctions is \textit{fine} if many of the
properties of standard functions are satisfied.

\begin{definition}
\label{A}A space of ultrafunctions $V%
{{}^\circ}%
,$ $V=BV\cap EPL\cap \mathcal{R},$ is called fine if

\begin{enumerate}
\item \label{uno}if $u,v\in V^{\circledcirc },$ then 
\begin{equation}
\overline{u_{\lambda }v_{\lambda }}\in V_{\lambda }  \label{F1}
\end{equation}

\item \label{due}there is a linear internal functional 
\begin{equation*}
\doint :V%
{{}^\circ}%
\rightarrow \mathbb{E}
\end{equation*}%
called \textbf{pointwise integral} which satisfies the following properties:

\begin{enumerate}
\item \label{11}if $u=\ \lim_{\lambda \uparrow \Lambda }u_{\lambda },$ $%
u_{\lambda }\in V_{\lambda },$ then%
\begin{equation}
\doint u(x)dx=\lim_{\lambda \uparrow \Lambda }\ \dint u_{\lambda
}(x)dx=\int^{\ast }u_{\Lambda }dx.  \label{int}
\end{equation}

\item \label{22}for every $a\in \Gamma ,$ 
\begin{equation}
\doint \chi _{a}(x)dx>0;  \label{lana}
\end{equation}

\item \label{YY}if $f\in \mathcal{L}^{1},$ then, 
\begin{equation}
\doint f%
{{}^\circ}%
(x)\ dx\sim \int f(x)\ dx.  \label{YYY}
\end{equation}
\end{enumerate}

\item \label{33}there are $N$ internal operators 
\begin{equation*}
D_{i}:V%
{{}^\circ}%
\rightarrow V%
{{}^\circ}%
,\ i=1,...,N
\end{equation*}%
called \textbf{generalized partial derivatives} such that the following
properties are fulfilled:

\begin{enumerate}
\item \label{3}for every $u,v\in V%
{{}^\circ}%
,$ 
\begin{equation}
\doint D_{i}u(x)v(x)dx=\lim_{\lambda \uparrow \Lambda }\int \partial
_{i}u_{\lambda }(x)v_{\lambda }(x)dx=\int^{\ast }\left( \partial _{i}^{\ast
}u_{\Lambda }\right) v_{\Lambda }dx\mathbf{;}  \label{A2}
\end{equation}

\item \label{2+}for every $a\in \Gamma ,$%
\begin{equation}
\mathfrak{supp}\left[ D_{i}\chi _{a}\right] \subset \mathfrak{mon}(a).
\label{loc}
\end{equation}
\end{enumerate}
\end{enumerate}
\end{definition}

Now some comments on Def. \ref{A}.

Assumptions (\ref{uno}) states that $V%
{{}^\circ}%
$ is sufficiently large with respect to 
\begin{equation*}
V^{\circledcirc }=\lim_{\lambda \uparrow \Lambda }V\cap \lambda .
\end{equation*}%
This is a technical assumption which simplifies the definition of regular
ultrafunctions (see section \ref{RU}).

Assumption (\ref{11}) is nothing else but the definition of the pointwise
integral. Assumption (\ref{22}) is quite natural, but it is not satisfied by
every space of ultrafunctions; it is very important and, among other things,
it implies that the bilinear form $(u,v)\mapsto \doint uv\ dx$ defines a
scalar product (see section \ref{spu}).

Assumption (\ref{3}) defines the generalized derivative. By (\ref{A2}) we
get that for every function $f\in C^{1,1}$ and every $x\in \mathbb{R}^{N},$%
(see Corollary \ref{eee+}), 
\begin{equation}
D_{i}f%
{{}^\circ}%
(x)=\partial _{i}f(x);  \label{bela}
\end{equation}%
hence the generalized derivative extends the usual derivative; by (\ref{lina}%
) we have that,$\ \forall x\in \Gamma $,%
\begin{equation*}
D_{i}f%
{{}^\circ}%
(x)=\sum_{a\in \Gamma }\partial _{i}^{\ast }f^{\ast }(a)\chi _{a}(x).
\end{equation*}%
Moreover,%
\begin{equation}
\doint D_{i}u(x)v(x)dx=-\doint D_{i}u(x)v(x)dx  \label{ip}
\end{equation}%
since%
\begin{eqnarray*}
\doint D_{i}u(x)v(x)dx &=&\lim_{\lambda \uparrow \Lambda }\int \partial
_{i}u_{\lambda }(x)v_{\lambda }(x)dx \\
&=&-\lim_{\lambda \uparrow \Lambda }\int u_{\lambda }(x)\partial
_{i}v_{\lambda }(x)dx=-\doint D_{i}u(x)v(x)dx
\end{eqnarray*}

Equation (\ref{ip}) is of primary importance in the theory of weak
derivatives, distributions, calculus of variations etc. Usually this
equality is deduced by the Leibniz rule%
\begin{equation*}
D(fg)=Dfg+fDg
\end{equation*}%
However, it is inconsistent to assume that Leibniz rule be satisfied by
every pair of ultrafunction (see \cite{Schwartz} and the discussion in
Section \ref{RU}). Nevertheless the identity (\ref{ip}) holds for the fine
ultrafunctions. In particular, by (\ref{ip}) and (\ref{bela}), we have that $%
\forall u\in V%
{{}^\circ}%
,\ \forall \varphi \in \mathcal{D},$ 
\begin{equation*}
\doint D_{i}u\varphi 
{{}^\circ}%
dx=-\doint uD_{i}\varphi 
{{}^\circ}%
dx=-\doint u~\left( \partial _{i}\varphi \right) 
{{}^\circ}%
dx;
\end{equation*}%
this equality relates the generalized derivative to the notion of weak
derivative.

Property (\ref{2+}) is a natural request and you expect that it is satisfied
by (\ref{A2}); on the contrary, it does not follows from the other
properties and it needs to be stated explicitly.

\bigskip

While the construction of a generic space of ultrafunction is a relatively
easy task, the construction of a space of fine ultrafunctions is much more
delicate. We have the following theorem which is one of the main results of
this paper:

\begin{theorem}
\label{main1}The requests of Def. \ref{A} are consistent.
\end{theorem}

The proof of this theorem is rather involved and it will be given in section %
\ref{CSU} (see Th. \ref{main}). The rest of this section and section \ref{BP}
will be devoted in showing that a fine space ultrafunctions provides a quite
rich structure and many interesting and \textit{natural }properties can be
proved in a relative simple way.

\begin{remark}
A space of ultrafunctions cannot be uniquely defined, since its existence
depends on $\mathbb{E}$ and hence on Zorn's Lemma. However, if we exclude
the choice of the space $V$, the properties required in Def. \ref{A} are
quite natural; hence, if a physical phenomenon is modelled by
ultrafunctions, the properties which can be deduced can be considered
reliable.
\end{remark}

From now on, we will treat only with fine ultrafunctions and the word "fine"
will be usually omitted.

\subsection{The pointwise scalar product of ultrafunctions\label{spu}}

By (\ref{lana}), we have that $d(a)>0$ and hence, the pointwise integral
allows to defines the following scalar product which we will call \textbf{%
pointwise scalar product}: 
\begin{equation}
\doint u(x)v(x)dx=\sum_{x\in \Gamma }u(x)v(x)d(x).  \label{psc}
\end{equation}%
If $f,g,fg\in V,$ we have that 
\begin{equation}
\doint f%
{{}^\circ}%
g%
{{}^\circ}%
dx=\int^{\ast }f^{\ast }g^{\ast }dx=\int fg\ dx;  \label{bella}
\end{equation}%
however we must be careful since, for some $u,v\in V%
{{}^\circ}%
$ such that $u_{\Lambda }v_{\Lambda }\notin V_{\Lambda },$ we might have
that 
\begin{equation}
\doint uv\ dx\neq \int^{\ast }u_{\Lambda }v_{\Lambda }dx.  \label{brutta}
\end{equation}%
even if these two quantities are not too different.

\textbf{Example:} Let $\Omega \subset \mathbb{R}^{N}$ be a bounded open set.
We have that 
\begin{equation*}
\doint \chi _{\Omega }\chi _{\Omega }\ dx\neq \doint \chi _{\overline{\Omega 
}}\chi _{\overline{\Omega }}\ dx;
\end{equation*}%
in fact 
\begin{equation*}
\doint \chi _{\overline{\Omega }}\chi _{\overline{\Omega }}\ dx-\doint \chi
_{\Omega }\chi _{\Omega }\ dx=\doint \left( \chi _{\overline{\Omega }}-\chi
_{\Omega }\right) dx=\doint \chi _{\partial \Omega }dx>0
\end{equation*}%
since $\chi _{\partial \Omega }>0.$ Hence, it is not possible that 
\begin{equation*}
\doint f%
{{}^\circ}%
g%
{{}^\circ}%
dx=\int fg\ dx
\end{equation*}%
for all measurable functions.

In any case, by (\ref{F1}), if $u,v\in \left( V\cap C^{0}\right)
^{\circledcirc }$ 
\begin{equation}
u,v\in V_{\Lambda }\Rightarrow \doint uv\ dx=\lim_{\lambda \uparrow \Lambda
}\ \dint u_{\lambda }v_{\lambda }dx=\int^{\ast }u_{\Lambda }v_{\Lambda }dx.
\label{uliva}
\end{equation}

The pointwise scalar product allows to get the ultrafunction analogous of
the Riesz representation theorem in the following form:

\begin{theorem}
\label{riz}If 
\begin{equation*}
\Phi :V%
{{}^\circ}%
\rightarrow \mathbb{E}
\end{equation*}%
is an internal linear functional, there exists an unique ultrafunction $%
u_{\Phi }$ such that, $\forall v\in V%
{{}^\circ}%
,$ 
\begin{equation*}
\Phi \left( v\right) =\doint u_{\Phi }v\ dx
\end{equation*}
\end{theorem}

\textbf{Proof}: The scalar product (\ref{psc}) is the $\Lambda $-limit of a
net of scalar products $\left\langle \cdot \ ,\ \cdot \right\rangle
_{_{\lambda }}$ defined over $V_{\lambda }$. Since $\Phi $ is an internal
functional, there exists a net of functionals $\Phi :V_{\lambda }\rightarrow 
\mathbb{R}$ such that 
\begin{equation*}
\Phi =\lim_{\lambda \uparrow \Lambda }\Phi _{_{\lambda }}
\end{equation*}%
and hence $\exists u_{\Phi ,\lambda },\ \forall v\in V_{\lambda },$ 
\begin{equation*}
\Phi _{_{\lambda }}\left( v\right) =\int u_{\Phi ,\lambda },v\ dx
\end{equation*}%
Taking 
\begin{equation*}
u_{\Phi }:=\lim_{\lambda \uparrow \Lambda }\ u_{\Phi ,\lambda }
\end{equation*}%
we get that $\forall v\in V%
{{}^\circ}%
,$%
\begin{equation*}
\Phi \left( v\right) =\lim_{\lambda \uparrow \Lambda }\Phi _{_{\lambda
}}\left( v_{_{\lambda }}\right) =\lim_{\lambda \uparrow \Lambda }\ \int
u_{\Phi ,\lambda },v\ dx=\doint u_{\Phi }v\ dx.
\end{equation*}

$\square $

\bigskip

The pointwise product provides the \textbf{pointwise (Euclidean) norm} of an
ultrafunction: 
\begin{equation*}
\left\Vert u\right\Vert =\left( \sum_{a\in \Gamma }|u(a)|^{2}d(a)\right) ^{%
\frac{1}{2}}=\left( \doint |u(x)|^{2}dx\right) ^{\frac{1}{2}}.
\end{equation*}%
Also, we can define other norms which might be useful in the applications:
for $p\in \left[ 1,\infty \right) ,$ we set%
\begin{equation*}
\left\Vert u\right\Vert _{p}=\left( \sum_{a\in \Gamma }|u(a)|^{p}\
d(a)\right) ^{\frac{1}{p}}=\left( \doint |u(x)|^{p}\ dx\right) ^{\frac{1}{p}}
\end{equation*}%
and, for $p=\infty $, obviously, we set%
\begin{equation*}
\left\Vert u\right\Vert _{\infty }=\max_{a\in \Gamma }\ |u(a)|.
\end{equation*}%
Notice that all these norms are equivalent in the sense that given any two
norms $\left\Vert u\right\Vert _{p}$ and $\left\Vert u\right\Vert _{q},$
there exists two numbers $m$ and $M$ such that $\forall u\in V%
{{}^\circ}%
,$%
\begin{equation}
m\leq \frac{\left\Vert u\right\Vert _{p}}{\left\Vert u\right\Vert _{q}}\leq M
\label{norma}
\end{equation}%
Of course, if $p\neq q$, $m$ is an infinitesimal number and $M$ is an
infinite number.

The scalar product (\ref{psc}) also allows to define the \textbf{delta (or
the Dirac) ultrafunction }as follows: for every $a\in \Gamma $,%
\begin{equation}
\delta _{a}(x)=\frac{\chi _{a}(x)}{d(a)}.  \label{dirac2}
\end{equation}%
As it is natural to expect, for every $u\in V%
{{}^\circ}%
,$ we have that%
\begin{equation*}
\doint \delta _{a}(x)u(x)dx=\dsum\limits_{x\in \Gamma }u(x)\delta
_{a}(x)d(x)=\dsum\limits_{x\in \Gamma }u(x)\frac{\chi _{a}(x)}{d(a)}d(x)=u(a)
\end{equation*}%
The delta ultrafunctions are orthogonal with each other with respect to the
pointwise scalar product; hence, if normalized, they provide an orthonormal
basis, called \textbf{delta-basis}, given by%
\begin{equation*}
\left\{ \sqrt{\delta _{a}}\right\} _{a\in \Gamma }=\left\{ \frac{\chi _{a}}{%
\sqrt{d(a)}}\right\} _{a\in \Gamma }
\end{equation*}%
So, the identity (\ref{pu}) can be rewritten as follows:%
\begin{equation*}
u(x)=\sum_{a\in \Gamma }\left( \doint u(\xi )\sqrt{\delta _{a}(\xi )}d\xi
\right) \sqrt{\delta _{a}(x)}.
\end{equation*}

\subsection{Regular and smooth ultrafunctions\label{RU}}

\begin{theorem}
\label{eee}If $f\in C_{c}^{1,1}$ 
\begin{equation*}
D_{i}f%
{{}^\circ}%
(x)=\left( \partial _{i}^{\ast }f^{\ast }\right) 
{{}^\circ}%
.
\end{equation*}
\end{theorem}

\textbf{Proof: }If $f\in C_{c}^{1,1},$ then $\partial _{i}f\in
C_{c}^{0,1}\subset V\subset V^{\circledcirc };$ then, $\forall v\in
V^{\circledcirc },$ 
\begin{equation*}
\doint D_{i}f%
{{}^\circ}%
v\ dx=\int^{\ast }\partial _{i}^{\ast }f^{\ast }v_{\Lambda }dx
\end{equation*}%
Since $\left( \partial _{i}f\right) 
{{}^\circ}%
,v\in V^{\circledcirc },$ by (\ref{F1}), $\partial _{i}fv_{\lambda }=%
\overline{\partial _{i}fv_{\lambda }}\in V_{\lambda },$ and hence, by (\ref%
{int}) 
\begin{equation*}
\doint D_{i}f%
{{}^\circ}%
vdx=\int^{\ast }\partial _{i}^{\ast }f^{\ast }v_{\Lambda }dx=\doint \left(
\partial _{i}^{\ast }f^{\ast }\right) 
{{}^\circ}%
vdx
\end{equation*}%
Since this equality holds for every $v\in V^{\circledcirc },$ the conclusion
follows.

$\square $

\begin{corollary}
\label{eee+}If $f\in C^{1,1}$ in a neighborhood $B_{\varepsilon }(x_{0})$ of 
$x_{0}\in \mathbb{R}^{N};$ then 
\begin{equation*}
D_{i}f%
{{}^\circ}%
(x_{0})=\partial _{i}f(x_{0}).
\end{equation*}
\end{corollary}

\textbf{Proof}: Let $\varphi \in C_{c}^{\infty }$ be $=1$ for $x\in
B_{\varepsilon /2}(x)$ and null for $x\notin B_{\varepsilon /2}(x).$ Then $%
\varphi f\in C_{c}^{1,1}$ and by Th. \ref{eee}, 
\begin{equation*}
D_{i}\left( \varphi f\right) 
{{}^\circ}%
(x)=\left( \partial _{i}^{\ast }\varphi ^{\ast }f^{\ast }\right) 
{{}^\circ}%
(x)
\end{equation*}%
By (\ref{loc}), $D_{i}\left( \varphi f\right) 
{{}^\circ}%
(x)$ depends only on the values in $\mathfrak{mon}(x_{0})\subset
B_{\varepsilon }(x_{0})^{\ast }$ and hence, for every $x\in \mathfrak{mon}%
(x_{0}),$ 
\begin{equation*}
D_{i}f%
{{}^\circ}%
(x)=\left( \partial _{i}^{\ast }f^{\ast }\right) 
{{}^\circ}%
(x)
\end{equation*}%
and from here the conclusion.

$\square $

\bigskip

\textbf{Example 1}. If 
\begin{equation*}
f(x)=\min \left( 0,x\right)
\end{equation*}%
then $Df(0)=\frac{1}{2}.$

\bigskip

\textbf{Example 2}. If 
\begin{equation*}
f(x)=\int_{0}^{x}t\sin \frac{1}{t^{2}}dt
\end{equation*}%
by Cor. \ref{eee+}, for $x\in \mathbb{R}\backslash \left\{ 0\right\} $, $%
Df(x)=x\sin \frac{1}{x^{2}},$ but the properties of Def. \ref{A} do not
guarantee that $Df(0)=0,$ since%
\begin{equation*}
x\sin \frac{1}{x^{2}}\notin BV.
\end{equation*}%
In any case, we have that 
\begin{equation*}
Df(0)=\doint D_{i}u(x)\delta _{0}(x)dx=\lim_{\lambda \uparrow \Lambda }\int
\partial _{i}u_{\lambda }(x)\delta _{0,\lambda }(x)dx\sim 0.
\end{equation*}

\bigskip

As we already remarked the Leibniz rule does not hold for ultrafunctions and
it is not possible to define a generalized derivative which has this
property. It is easy to check this fact for idempotent functions, in fact by
the Leibniz rule, we should have%
\begin{equation*}
D\chi _{E}^{2}=2\chi _{E}D\chi _{E}
\end{equation*}%
and since $\chi _{E}^{2}=\chi _{E},$ we deduce that 
\begin{equation*}
D\chi _{E}=2\chi _{E}D\chi _{E}
\end{equation*}%
and hence, for every $x\in E,$ 
\begin{equation}
D\chi _{E}=0  \label{schw}
\end{equation}%
and this fact contradicts any reasonable generalization of the notion of
derivative. Actually the Schwartz impossibility theorem states that the
Leibniz rule cannot be satisfied by any algebra which contains, not only the
idempotent functions, but also the continuous functions (see \cite{Schwartz},%
\cite{algebra}). So it is interesting to investigate the subspaces of
ultrafunctions for which the Leibniz rule holds and more in general to
determine spaces in which many of the usual properties of smooth functions
be satisfied. This is important for the applications when we want to study
the qualitative properties (and in particular the regularity) of the
solutions of an equation.

\begin{definition}
\label{ddd}We set%
\begin{equation*}
U^{0}:=\left\{ u\in V%
{{}^\circ}%
\ |\ u_{\Lambda }\in \left( C^{0,1}\cap V\right) ^{\circledcirc }\right\}
\end{equation*}%
and for $m>1,$ we define by induction%
\begin{equation*}
U^{m}:=\left\{ u\in U^{m-1}\ |\ \forall i=1,...,N,\ D_{i}u\in U^{m-1}\right\}
\end{equation*}%
If $u\in U^{m}$, we will say that $u$ is $m$-regular.
\end{definition}

Let us see the main properties of the spaces $U^{m}.$

\begin{theorem}
\label{TU}The spaces of $m$-regular ultrafunctions satisfy the following
properties:

\begin{enumerate}
\item \label{U0}if $u\in U^{1},$ then%
\begin{equation}
D_{i}u=\left( \partial _{i}^{\ast }u_{\Lambda }\right) 
{{}^\circ}%
;  \label{U4}
\end{equation}

\item \label{U3}if $u\in U^{m},$ then 
\begin{equation*}
D_{i}u(x)=\left[ \partial _{i}^{\ast }u_{\Lambda }(x)\right] 
{{}^\circ}%
\end{equation*}

\item \label{U1}if $f\in C_{c}^{m,1},$ then $f%
{{}^\circ}%
\in U^{m};$

\item \label{U7}if $u,v,uv\in U^{1}$ then the Leibniz rule holds:%
\begin{equation*}
D_{i}\left( uv\right) =D_{i}uv+uD_{i}v.
\end{equation*}
\end{enumerate}
\end{theorem}

\textbf{Proof}: \ref{U0} - If $u\in U^{1},\ D_{i}u\in U^{0}\subset
V^{\circledcirc };$ then, by (\ref{F1}), $\forall v\in U^{0}\subset
V^{\circledcirc },$ $\overline{\left( D_{i}u\right) _{\lambda }v_{\lambda }}%
\in V_{\lambda }.$ Since $D_{i}u\ $is continuous, $\left( D_{i}u\right)
_{\lambda }v_{\lambda }=\overline{\left( D_{i}u\right) _{\lambda }v_{\lambda
}}\in V_{\lambda }$ and hence $\left( D_{i}u\right) _{\Lambda }v_{\Lambda
}\in V_{\Lambda }$ and 
\begin{equation*}
\doint D_{i}uv\ dx=\int^{\ast }\left( D_{i}u\right) _{\Lambda }v_{\Lambda }\
dx
\end{equation*}%
Then, by (\ref{A2}) we have that%
\begin{equation*}
\int^{\ast }\left( D_{i}u\right) _{\Lambda }v_{\Lambda }\ dx=\int^{\ast
}\left( \partial _{i}^{\ast }u_{\Lambda }\right) v_{\Lambda }\ dx
\end{equation*}%
and so, 
\begin{equation*}
\left( D_{i}u\right) _{\Lambda }=\partial _{i}^{\ast }u_{\Lambda }.
\end{equation*}

\ref{U3} - it follows trivially by (\ref{U0}) since $U^{m}\subset U^{1}$

\ref{U1} - For $m=0,$ $f\in C_{c}^{0,1},$ and hence $f%
{{}^\circ}%
\in U^{0}$. For $m\geq 1,$ since $\partial _{i}f\in C_{c}^{m-1,1},$ by (\ref%
{U4}), $\left( \partial _{i}f\right) 
{{}^\circ}%
=D_{i}f%
{{}^\circ}%
;$ hence $D_{i}f%
{{}^\circ}%
\in U^{m-1}$. So $f%
{{}^\circ}%
\in U^{m}.$

\ref{U7} - If $u,v\in U^{1}$ then $u,v,D_{i}u,D_{i}v\in U^{0}$ and so, using
(\ref{U0}) and (\ref{F1}), we have that 
\begin{eqnarray*}
D_{i}\left( uv\right) &=&\left[ \partial _{i}^{\ast }\left( uv\right)
_{\Lambda }\right] 
{{}^\circ}%
=\left[ \partial _{i}^{\ast }\left( u_{\Lambda }v_{\Lambda }\right) \right] 
{{}^\circ}%
=\left[ \partial _{i}^{\ast }u_{\Lambda }v_{\Lambda }+u_{\Lambda }\partial
_{i}^{\ast }v_{\Lambda }\right] 
{{}^\circ}
\\
&=&\left[ \partial _{i}^{\ast }u_{\Lambda }v_{\Lambda }\right] 
{{}^\circ}%
+\left[ u_{\Lambda }\partial _{i}^{\ast }v_{\Lambda }\right] 
{{}^\circ}
\\
&=&\left[ \partial _{i}^{\ast }u_{\Lambda }\right] 
{{}^\circ}%
v+u\left[ \partial _{i}^{\ast }v_{\Lambda }\right] 
{{}^\circ}%
=D_{i}uv+uD_{i}v
\end{eqnarray*}

$\square $

\bigskip

Since the $U^{m}$'s have hyperfinite dimension then there exist finite
dimensional spaces $U_{\lambda }^{m}$'s such that 
\begin{equation}
U^{m}=\left\{ u%
{{}^\circ}%
\ |\ u\in U_{\Lambda }^{m}\right\} ;\ \ \ U_{\Lambda }^{m}:=\lim_{\lambda
\uparrow \Lambda }\ U_{\lambda }^{m};  \label{cicca}
\end{equation}%
the space of \textbf{smooth ultrafunctions} (or $\infty $-regular
ultrafunctions) is defined as follows:\texttt{\ }%
\begin{equation}
U^{\infty }:=\dbigcap\limits_{m\leq \alpha }U^{m}=\left[ \lim_{\lambda
\uparrow \Lambda }\left( \dbigcap\limits_{m\leq \left\vert \lambda \cap 
\mathbb{N}\right\vert }U_{\lambda }^{m}\right) \right] ^{\circ }.
\label{Uinf}
\end{equation}%
where $\alpha $ has been defined by (\ref{alfa}). Clearly $U^{\infty }\neq
\varnothing ,$ since 
\begin{equation}
f\in \mathcal{D}\Rightarrow f%
{{}^\circ}%
\in U^{\infty }.  \label{ciccia}
\end{equation}%
For $m\in \mathbb{N}\cup \left\{ \infty \right\} ,$ we set 
\begin{equation*}
\left( U^{m}\right) ^{\perp }=\left\{ u\in V%
{{}^\circ}%
\ |\ \forall \psi \in U^{m},\ \doint u\psi \ dx=0\right\}
\end{equation*}%
and we denote by $\Pi _{m}u$ and $\Pi _{m}^{\perp }u$ the relative
"orthogonal" projection of $u$ on $U^{m}\ $and $\left( U^{m}\right) ^{\perp
}.$ Then every ultrafunction $u$ can be split as follows%
\begin{equation}
u=\Pi _{m}u+\Pi _{m}^{\perp }u  \label{181}
\end{equation}

\begin{definition}
\label{sm}Given the splitting (\ref{181}), $\forall m\in \mathbb{N}\cup
\left\{ \infty \right\} $, $\Pi _{m}u$ will be called the $m$\textbf{-reguar}
\textbf{part} of $u$ and $\Pi _{m}^{\perp }u$ the $m$\textbf{-singular part}
of $u.$
\end{definition}

\begin{remark}
\label{reg}It is possible to define different types of regular ultrafunction
namely we can choose different subspaces of $V%
{{}^\circ}%
$ that satisfy suitable conditions. For example, we can set%
\begin{equation*}
C_{\Lambda }^{m}=\lim_{\lambda \uparrow \Lambda }\ C_{0}^{m}\cap \lambda
\end{equation*}%
We have that $C_{\Lambda }^{m}\supset U^{m}$ and hence the functions in $%
C_{\Lambda }^{m}$ satisfy less properties. Similarly, we can choose more
regular spaces such as%
\begin{equation}
U^{m,p}:=\left\{ u\in U^{m}\ |\ \left\vert u\right\vert ^{p-2}u\in
U^{0}\right\} ;\ p\geq 2.  \label{ganzetta}
\end{equation}%
Of course the choice of a particular space depends on the problems that we
would like to treat. We can make an analogy with the theory of
distributions; in this case the spaces $C^{m}$'s and the Sobolev spaces $%
H^{m,p}$'s can be considered as subspaces of $\mathcal{D}^{\prime }$ which
present different kinds of regularity.
\end{remark}

\subsection{Time-dependent ultrafunctions\label{TDU}}

In evolution problems the time variable plays a different role that the
space variables; then the functional spaces used in these problems (e.g. $%
C^{k}(\left[ 0,T\right] ,H_{0}^{1}(\Omega )),$ $L_{loc}^{p}(\mathbb{R}%
,H^{k}(\Omega ))$ etc.) reflect this fact. The same is true in the frame of
ultrafunctions. This section is devoted in the description of the
appropriate ultrafunction-spaces for evolution problems.

First of all we need to recall some well known facts about free modules:

\begin{definition}
Given a ring $R$ and a module $M$ over $R$, the set $B\subset M$ is a basis
for $M$ if:

\begin{itemize}
\item $B$ is a generating set for $M$; that is to say, every element of $M$
is a finite sum of elements of $B$ multiplied by coefficients in $R$;

\item $B$ is linearly independent.
\end{itemize}
\end{definition}

\begin{definition}
A free module is a module with a basis.
\end{definition}

The following is a well known theorem:

\begin{theorem}
If $R$ is a commutative ring and $M$ is a free $R$-module, then all the
bases of $E$ have the same cardinality. The cardinality of a basis is called 
\textbf{rank} of $M$.
\end{theorem}

We will describe some free modules which will be used in the following part
of this paper.

\bigskip

\textbf{Examples}: (i) If $\Gamma $ is a finite set then $C^{k}(\mathbb{R}%
)^{\Gamma }$ is a free module over $C^{k}(\mathbb{R})$ of rank $\left\vert
\Gamma \right\vert $ and a basis is given by 
\begin{equation*}
\left\{ \chi _{a}\right\} _{a\in \Gamma }
\end{equation*}%
If $u\in C^{k}(\mathbb{R})^{\Gamma },$ then 
\begin{equation*}
u(t,x)=\sum_{a\in \Gamma }c(t)\chi _{a}(x)
\end{equation*}%
Then we may think of $u$ as a real function defined in $\mathbb{R}\times
\Gamma $ or a function in $C^{k}(\mathbb{R},\Gamma )$ or in $C^{k}(\mathbb{R}%
)\otimes _{\mathbb{R}}\mathfrak{F}\left( \Gamma \right) .$

(ii) Let $W\subset \mathfrak{F}\left( \mathbb{R}^{N}\right) $ be a vector
space of finite dimension, then 
\begin{equation*}
C^{k}(\mathbb{R},W):=C^{k}(\mathbb{R})\otimes _{\mathbb{R}}W\subset 
\mathfrak{F}\left( \mathbb{R}^{N+1}\right)
\end{equation*}%
is a free $C^{k}$-module of rank equal to $\dim W$, namely every function $%
f\in C^{k}(\mathbb{R},W)$ can be written as follows:%
\begin{equation*}
f(t,x):=\sum_{k=1}^{\dim W}c_{k}(t)e_{k}(x)
\end{equation*}%
where $c_{k}\in C^{k}(\mathbb{R})$ and $\left\{ e_{k}\right\} $ is any basis
in $W$.

(iii) If $W=\lim_{\lambda \uparrow \Lambda }\ W_{\lambda }\subset \mathfrak{F%
}\left( \mathbb{R}^{N}\right) ^{\ast }$ is an internal vector space of
hyperfinite dimension, then by $C^{k}(\mathbb{E},W)$ we denote the internal $%
C^{k}(\mathbb{R})^{\ast }$-module defined by%
\begin{equation*}
C^{k}(\mathbb{E},W)=\lim_{\lambda \uparrow \Lambda }\ C^{k}(\mathbb{R}%
,W_{\lambda }).
\end{equation*}

Now we are ready to define the time-dependent ultrafunctions:

\begin{definition}
\label{lella3}The space of\textbf{\ time dependent ultrafunctions }of order $%
k\in \mathbb{N}$ is the free $C^{k}(\mathbb{R})^{\ast }$-module given by 
\begin{equation*}
C^{k}(\mathbb{E},V%
{{}^\circ}%
):=\left\{ u%
{{}^\circ}%
=u|_{\mathbb{E}\times \Gamma }\ |\ u\in C^{k}(\mathbb{E},V_{\Lambda
})\right\} .
\end{equation*}
\end{definition}

Every time-dependent ultrafunction can be represented by the following
hyperfinite sum:%
\begin{equation*}
u(t,x)=\sum_{a\in \Gamma }c(t)\chi _{a}(x)
\end{equation*}%
where $c(t)\in C^{k}(\mathbb{R})^{\ast }$ and $\left\{ \chi _{a}(x)\right\}
_{a\in \Gamma }$ is the canonical basis of $V%
{{}^\circ}%
$. The map%
\begin{equation*}
(%
{{}^\circ}%
):C^{k}(\mathbb{E},V_{\Lambda })\rightarrow C^{k}(\mathbb{E},V%
{{}^\circ}%
);\ \ u%
{{}^\circ}%
=u|_{\mathbb{E}\times \Gamma }
\end{equation*}%
is an isomorphism between free $C^{k}(\mathbb{R})^{\ast }$-modules: in fact,
by using (\ref{sigma}), we have that $\left\{ \sigma _{a}\right\} _{a\in
\Gamma }$ is a basis of $C^{k}(\mathbb{E},V_{\Lambda })$ and we have that%
\begin{equation}
\left( \sum_{a\in \Gamma }c(t)\sigma _{a}(x)\right) ^{\circ }=\sum_{a\in
\Gamma }c(t)\chi _{a}(x)  \label{lang}
\end{equation}%
The restriction map $(%
{{}^\circ}%
)$ can be extended to a $C^{k}(\mathbb{R})^{\ast }$-module homomorphism 
\begin{equation*}
(%
{{}^\circ}%
):C^{k}(\mathbb{R},\mathfrak{F}\left( \mathbb{R}^{N}\right) )^{\ast
}\rightarrow C^{k}(\mathbb{E},V%
{{}^\circ}%
)
\end{equation*}%
by setting%
\begin{equation*}
w%
{{}^\circ}%
(t,x)=\sum_{a\in \Gamma }w(t,a)\chi _{a}(x)
\end{equation*}%
and, of course, by (\ref{cicca}), 
\begin{equation*}
w_{\Lambda }(t,x)=\sum_{a\in \Gamma }w(t,a)\sigma _{a}(x)\in C^{k}(\mathbb{E}%
,V%
{{}^\circ}%
)
\end{equation*}%
In particular, if $f\in C^{k}(\mathbb{R},\mathfrak{F}\left( \mathbb{R}%
^{N}\right) ),$ we have that 
\begin{equation}
f%
{{}^\circ}%
(t,x)=\sum_{a\in \Gamma }f^{\ast }(t,a)\chi _{a}(x)\ \ \ \ \text{and}\ \ \ \
\ f_{\Lambda }(t,x)=\sum_{a\in \Gamma }f^{\ast }(t,a)\sigma _{a}(x)
\label{nina}
\end{equation}%
Observe that, unlikely of (\ref{ciccia2}), in this case $f_{\Lambda
}(t,x)\neq f^{\ast }(t,x)$ for some $t\in \mathbb{E}\backslash \mathbb{R}$;
this fact does not prevent the theory to work, but, in some circumstances,
it is necessary to be careful.

The notion of generalized derivative in the space variable is trivially
defined by linearity:%
\begin{equation}
D_{i}u(t,x)=D_{i}\left( \sum_{a\in \Gamma }c(t)\chi _{a}(x)\right)
=\sum_{a\in \Gamma }c(t)D_{i}\chi _{a}(x).  \label{dx}
\end{equation}%
It is not necessary to introduce a generalized time-derivative, since the
natural derivative%
\begin{equation*}
\partial _{t}^{\ast }:C^{k+1}(\mathbb{E},V%
{{}^\circ}%
)\rightarrow C^{k}(\mathbb{E},V%
{{}^\circ}%
),\ k\geq 0,
\end{equation*}%
is well defined by setting%
\begin{equation}
\partial _{t}^{\ast }u(t,x)=\partial _{t}^{\ast }\left( \sum_{a\in \Gamma
}c(t)\chi _{a}(x)\right) =\sum_{a\in \Gamma }\partial _{t}^{\ast }c(t)\chi
_{a}(x).  \label{dt}
\end{equation}%
In our applications, we do not need a generalized nor a weak time-derivative
for the functions in $C^{0}(\mathbb{E},V%
{{}^\circ}%
).$

\begin{theorem}
\label{TEV}If $f\in C^{k}(\mathbb{R},C_{c}^{m,1})$ and $m\geq 0$, then,$\
\partial _{t}^{\ast }f%
{{}^\circ}%
\in C^{k-1}(\mathbb{E},U^{m})$ and%
\begin{equation*}
\partial _{t}^{\ast }f%
{{}^\circ}%
(t,x)=\sum_{a\in \Gamma }\partial _{t}^{\ast }f^{\ast }(t,a)\chi _{a}(x).
\end{equation*}%
Moreover, for $i=1,...,N\ $and $m\geq 1$,$\ D_{i}f%
{{}^\circ}%
\in C^{k}(\mathbb{E},U^{m-1}),$ and 
\begin{equation*}
D_{i}f%
{{}^\circ}%
(t,x)=\sum_{a\in \Gamma }\partial _{i}^{\ast }f^{\ast }(t,a)\chi _{a}(x).
\end{equation*}
\end{theorem}

\textbf{Proof}. The first equality follows immediately from (\ref{nina}) and
(\ref{dt}).

Let us prove the second statement. By the definition \ref{ddd}, if $w\in
C^{k}(\mathbb{E},U^{m}),$ 
\begin{equation*}
w_{\Lambda }\left( t,x\right) =\sum_{a\in \Gamma }w(t,a)\sigma _{a}(x)\in
C^{k}(\mathbb{E},U_{\Lambda }^{m});\ 
\end{equation*}%
and by \ref{TU}.\ref{U3}, $\forall t\in \mathbb{E}$, $D_{i}w(t,\cdot )=\left[
\partial _{i}^{\ast }w_{\Lambda }(t,\cdot )\right] 
{{}^\circ}%
\in U^{m-1}$ and hence%
\begin{equation*}
D_{i}w(t,x)=\left[ \partial _{i}^{\ast }w_{\Lambda }(t,x)\right] 
{{}^\circ}%
=\sum_{a\in \Gamma }\partial _{i}^{\ast }w(t,a)\chi _{a}(x)
\end{equation*}%
In particular, by (\ref{nina}), if $f\in C^{k}(\mathbb{R},C_{c}^{m,1}),$ 
\begin{equation*}
f_{\Lambda }\left( t,x\right) =\sum_{a\in \Gamma }f^{\ast }(t,a)\sigma
_{a}(x)\in C^{k}(\mathbb{E},U_{\Lambda }^{m})
\end{equation*}%
and hence 
\begin{equation*}
D_{i}f%
{{}^\circ}%
(t,x)=\left[ \partial _{i}^{\ast }f^{\ast }(t,x)\right] 
{{}^\circ}%
=\sum_{a\in \Gamma }\partial _{i}^{\ast }f^{\ast }(t,a)\chi _{a}(x).
\end{equation*}

$\square $

\section{Basic properties of ultrafunctions\label{BP}}

In this section we analyze some properties of the fine ultrafunction that
seems interesting in themselves and/or relevant in the applications.

\subsection{Ultrafunctions and measures\label{UM}}

As we have seen, if $f\notin V$, in general%
\begin{equation*}
\doint f%
{{}^\circ}%
(x)\ dx\neq \int f(x)\ dx
\end{equation*}%
since $\doint f%
{{}^\circ}%
(x)dx$ takes account of the value of $f$ in any single point. Thus it is a
natural question to ask if there exists an ultrafunction $u$ some way
related to $f$ such that 
\begin{equation}
\doint u\ dx=\int f~dx.  \label{lilla}
\end{equation}%
This question has an easy answer if we think of $f$ as the density of a
measure $\mu _{f}$. In fact, the following definition appears quite natural:

\begin{definition}
\label{dualext}If $\mu \ $is a Radon measure, we define an ultrafunction $%
\mu 
{{}^\circ}%
$ as follows: for every $v\in V%
{{}^\circ}%
,$ we set%
\begin{equation*}
\doint v(x)\mu 
{{}^\circ}%
(x)dx=\lim_{\lambda \uparrow \Lambda }\int v_{\lambda }(x)d\mu
\end{equation*}
\end{definition}

Notice that the existence of $\mu 
{{}^\circ}%
$ is guaranteed by Th. \ref{riz}, since 
\begin{equation*}
\Phi (v):=\lim_{\lambda \uparrow \Lambda }\int v_{\lambda }(x)d\mu
\end{equation*}%
is an internal linear functional over $V%
{{}^\circ}%
$.

\textbf{Example:} 1 - If $\mu _{f}$ is a measure whose density is $f\in 
\mathcal{L}_{loc}^{1},$ then, $\forall v\in V%
{{}^\circ}%
$%
\begin{equation*}
\doint v\mu _{f}^{\circ }~dx=\lim_{\lambda \uparrow \Lambda }\int
f(x)v_{\lambda }(x)dx=\int^{\ast }f^{\ast }v_{\Lambda }dx.
\end{equation*}%
then, taking $f\in \mathcal{L}^{1},$ and $v=1%
{{}^\circ}%
$%
\begin{equation*}
\doint \mu _{f}^{\circ }dx=\int f~dx;
\end{equation*}%
then (\ref{lilla}) holds with $u=\mu _{f}^{\circ }.$

\textbf{Example:} 2 - If $\mathbf{\delta }_{a}$ is the Dirac measure, then $%
\mathbf{\delta }_{a}^{\circ }=\delta _{a}$ where $\delta _{a}$ is the Dirac
ultrafunction defined by (\ref{dirac2}).

\textbf{Example:} 3 - For all $u\in V%
{{}^\circ}%
,\ \partial _{i}u_{\Lambda }$ is a measure and we will denote by $\mu
_{\partial _{i}u}$ the related ultrafunction, namely, using the notation (%
\ref{nota}), we have that, $\forall v\in V%
{{}^\circ}%
$ 
\begin{equation*}
\doint v\mu _{\partial _{i}u}dx=\lim_{\lambda \uparrow \Lambda }\int
\partial _{i}u_{\lambda }v_{\lambda }d\mu
\end{equation*}%
and hence, by (\ref{A2}) 
\begin{equation}
\doint v\mu _{\partial _{i}u}dx=\int^{\ast }\left( \partial _{i}^{\ast
}u_{\Lambda }\right) v_{\Lambda }dx=\doint D_{i}u~v~dx  \label{pinta}
\end{equation}%
namely 
\begin{equation}
D_{i}u=\mu _{\partial _{i}u}.  \label{pla}
\end{equation}

\textbf{Example:} 4 - If $\Omega \subset \mathbb{R}^{N}$ is a set of finite
measure, and $\mu _{\Omega }$ is the measure whose density is $\chi _{\Omega
}$, then $\forall u\in V%
{{}^\circ}%
,$ we have that 
\begin{eqnarray}
\doint u(x)\mu _{\Omega }^{\circ }(x)dx &=&\lim_{\lambda \uparrow \Lambda }\
\int u_{\lambda }(x)d\mu _{\Omega }=\lim_{\lambda \uparrow \Lambda }\ \int
u_{\lambda }(x)\chi _{\Omega }dx  \notag \\
&=&\int^{\ast }u_{\Lambda }(x)\chi _{\Omega }^{\ast }~dx=\int_{\Omega ^{\ast
}}^{\ast }u_{\Lambda }(x)dx.  \label{polla}
\end{eqnarray}

In particular, if $f\in V,$%
\begin{equation}
\int_{\Omega }f\left( x\right) dx=\int^{\ast }f^{\ast }(x)\chi _{\Omega
}^{\ast }~dx=\doint f%
{{}^\circ}%
(x)\mu _{\Omega }^{\circ }(x)dx  \label{pollina}
\end{equation}

\textbf{Example:} 5 - If $\Omega \subset \mathbb{R}^{N}$ has a smooth
boundary $\partial \Omega $ whose $(N-1)$-dimensional measure is denoted by $%
S_{\partial \Omega },$ we have that%
\begin{equation}
\doint u(x)S_{\partial \Omega }^{\circ }(x)dx=\lim_{\lambda \uparrow \Lambda
}\ \int u_{\lambda }(x)dS_{\partial \Omega }=\int^{\ast }u_{\Lambda
}(x)~dS_{\partial \Omega }^{\ast }.  \label{pollaS}
\end{equation}

\textbf{Caveat: }It is absolutely natural to generalize equation (\ref{int3}%
) by the following definition: 
\begin{equation}
\doint_{E}u(x)dx:=\sum_{a\in E}u(a)d(a)  \label{u1}
\end{equation}%
for every set $E\subset \Gamma .$ At this point in is important to notice
that for a measurable set $\Omega \subset \mathbb{R}^{N}$ and $f\in C^{0}$%
\begin{equation}
\doint_{\Omega 
{{}^\circ}%
}f%
{{}^\circ}%
\left( x\right) dx\neq \doint f%
{{}^\circ}%
(x)\mu _{\Omega }^{\circ }(x)dx=\int_{\Omega }f\left( x\right) dx.
\label{sale+}
\end{equation}%
and%
\begin{equation*}
\doint_{\partial \Omega 
{{}^\circ}%
}f%
{{}^\circ}%
\left( x\right) dx\neq \doint f%
{{}^\circ}%
(x)S_{\Omega }^{\circ }(x)dx=\int_{\partial \Omega }f\left( x\right)
dS_{\partial \Omega }.
\end{equation*}%
\bigskip

The next proposition shows a useful way to represent $\mu _{f}^{\circ }$:

\begin{proposition}
\label{pilla}If $f\in \mathcal{L}_{loc}^{1}$, then%
\begin{equation}
\mu _{f}^{\circ }(x)=\sum_{a\in \Gamma }\left( \int^{\ast }f^{\ast
}(y)\sigma _{a}(y)dy\right) \delta _{a}(x)  \label{placido}
\end{equation}%
where $\left\{ \sigma _{a}(x)\right\} _{a\in \Gamma }$ is the $\sigma $%
-basis (see (\ref{sigma})).
\end{proposition}

\textbf{Proof: }Let $\sigma _{a,\lambda }\ $be the net such that $\sigma
_{a}=\ \lim_{\lambda \uparrow \Lambda }\ \sigma _{a,\lambda }.$ Then, by (%
\ref{dirac2}),%
\begin{eqnarray*}
\mu _{f}^{\circ }(a) &=&\doint \mu _{f}^{\circ }(y)\delta _{a}(y)dy=\frac{1}{%
d(a)}\doint \mu 
{{}^\circ}%
(y)\chi _{a}(y)dy=\frac{1}{d(a)}\lim_{\lambda \uparrow \Lambda }\int
f(y)\sigma _{a,\lambda }\left( y\right) dy \\
&=&\frac{1}{d(a)}\int^{\ast }f^{\ast }(y)\sigma _{a}(y)dy
\end{eqnarray*}%
Hence,%
\begin{eqnarray*}
\mu _{f}^{\circ }(x) &=&\sum_{a\in \Gamma }\mu _{f}^{\circ }(a)\chi
_{a}(x)=\sum_{a\in \Gamma }\left( \frac{1}{d(a)}\int^{\ast }f^{\ast
}(y)\sigma _{a}(y)dy\right) \chi _{a}(x) \\
&=&\sum_{a\in \Gamma }\left( \int^{\ast }f^{\ast }(y)\sigma _{a}(y)dy\right)
\ \delta _{a}(x)
\end{eqnarray*}

$\square $

\bigskip

Prop. \ref{pilla} suggests to generalize the operator $f^{\ast }\mapsto \mu
_{f}^{\circ }$ to an operator $w\mapsto \mu _{w}$ defined $\forall w\in
\left( L_{loc}^{1}\right) ^{\ast }$ by setting%
\begin{equation*}
\mu _{w}(x):=\sum_{a\in \Gamma }\left( \int^{\ast }w(y)\sigma
_{a}(y)dy\right) \delta _{a}(x)
\end{equation*}%
$\mu _{w}(x)$ can be considered as a sort of measure density defined on $%
\Gamma .$

Moreover, to simplify the notation, we set 
\begin{equation}
\mu _{E}(x)=\mu _{\chi _{E}}(x).  \label{polla+}
\end{equation}%
By the pointwise representation of $u$ (see (\ref{pu})), we get 
\begin{equation*}
\doint u(x)\mu _{E}(x)dx=\sum_{a\in \Gamma }u(a)\doint \chi _{a}(x)\mu
_{E}(x)dx;
\end{equation*}%
then setting $d_{E}(a)=\doint \chi _{a}(x)\mu _{E}(x)dx$,%
\begin{equation}
\doint u(x)\mu _{E}(x)dx=\sum_{a\in \Gamma }u(a)d_{E}(a)  \label{u2}
\end{equation}%
generalizing eq. (\ref{int3}). Notice the difference between the equalities (%
\ref{u1}) and (\ref{u2}). These equalities suggest the following notation:%
\begin{equation}
\doint u(x)\ d_{E}x=\doint u(x)\mu _{E}(x)dx  \label{polla++}
\end{equation}

In particular, if we take $f\in V,$ 
\begin{equation*}
\int_{\Omega }f(x)dx=\sum_{a\in \Gamma }f%
{{}^\circ}%
(a)d_{\Omega }(a)=\doint f%
{{}^\circ}%
(x)d_{\Omega }x
\end{equation*}%
\begin{equation*}
\int_{\partial \Omega }f(x)dx=\sum_{a\in \Gamma }f%
{{}^\circ}%
(a)d_{\partial \Omega }(a)=\doint f%
{{}^\circ}%
(x)d_{\partial \Omega }x.
\end{equation*}%
We end this section with the following

\begin{proposition}
\label{papa}$\mu _{w}(x)$ is the element of $V%
{{}^\circ}%
$ characterized by the following identity: $\forall u\in V%
{{}^\circ}%
,$ 
\begin{equation}
\doint u\left( x\right) \mu _{w}(x)dx=\int^{\ast }u_{\Lambda }\left(
x\right) w_{\Lambda }(x)dx=\doint w(x)\mu _{u}\left( x\right) dx
\label{domingo}
\end{equation}
\end{proposition}

\textbf{Proof. }We have that 
\begin{eqnarray*}
\doint \mu _{w}(x)u\left( x\right) dx &=&\doint \sum_{a\in \Gamma }\left(
\int^{\ast }w_{\Lambda }(y)\sigma _{a}(y)dy\right) \delta _{a}(x)u\left(
x\right) dx \\
&=&\sum_{a\in \Gamma }\left( \int^{\ast }w_{\Lambda }(y)\sigma
_{a}(y)dy\right) \doint \delta _{a}(x)u\left( x\right) dx \\
&=&\sum_{a\in \Gamma }\left( \int^{\ast }w_{\Lambda }(y)\sigma
_{a}(y)u\left( a\right) dy\right) =\int^{\ast }w_{\Lambda }(y)\left[
\sum_{a\in \Gamma }\sigma _{a}(y)u\left( a\right) \right] dy \\
&=&\int^{\ast }w_{\Lambda }(y)u_{\Lambda }\left( y\right) dy.
\end{eqnarray*}%
The last equality follows by symmetry.

$\square $

\bigskip

\subsection{The \textit{vicinity} of a set \label{pie}}

Given an open set $\Omega \subset \mathbb{R}^{N}$ and $f\in C^{1},$ then the
value of $\nabla f(x_{0})$ in a point $x_{0}\in \partial \Omega $ depends
only on the values which $f$ takes in $\Omega $ since 
\begin{equation*}
\nabla f(x_{0})=\lim_{\substack{ x\in \Omega ,  \\ x\rightarrow x_{0}}}\
\nabla f(x_{0})
\end{equation*}%
If $f$ is not continuous, $\nabla f$ is not defined, but $Df%
{{}^\circ}%
$ makes sense; however $Df%
{{}^\circ}%
(x_{0})$ in a point $x_{0}\in \partial \Omega 
{{}^\circ}%
$ depends on the values which $f$ takes in suitable points $y\sim x_{0}$
even if $y\notin \Omega 
{{}^\circ}%
.$ Roughly speaking, the vicinity of a set $E$ consists of the set of points
which influence the derivative of the points in $E.$ In order to make this
definition precise, we need to define a function $\theta _{E}\in V%
{{}^\circ}%
$ for every $E\subset \Gamma $. If $\Omega \subset \mathbb{R}^{N}$ is a set
such that $\bar{\chi}_{\Omega }\in V,$ we set%
\begin{equation}
\theta _{\Omega 
{{}^\circ}%
}=\bar{\chi}_{\Omega }^{\circ }  \label{pepe}
\end{equation}%
where the operator $u\mapsto \bar{u}$ has been defined by (\ref{due+}). We
want to generalize the above formula to every set $E\subset \Gamma .$ We put%
\begin{equation*}
\left( \chi _{E}\right) _{\Lambda }=\sum_{a\in E}\sigma _{a}(x)
\end{equation*}%
and 
\begin{equation*}
\theta _{E}=\left[ \overline{\left( \chi _{E}\right) _{\Lambda }}\right] 
{{}^\circ}%
\end{equation*}%
In the above formula, the operator $u\mapsto \bar{u}$ who has been defined
by (\ref{due+}) for $u\in \mathcal{L}_{loc}^{\infty },$ has been extended to 
$\left( \mathcal{L}_{loc}^{\infty }\right) ^{\ast }.\ \theta _{E}$ is
similar to the measure density $\mu _{E}$ but it is slightly different; for
example we have that $\forall f\in C^{1,1}$ and $\forall \Omega \subset 
\mathbb{R}^{N},$ bounded and open,%
\begin{equation*}
\doint f%
{{}^\circ}%
\mu _{\Omega 
{{}^\circ}%
}dx=\int f\ dx=\doint f%
{{}^\circ}%
\theta _{\Omega 
{{}^\circ}%
}dx
\end{equation*}%
since $f%
{{}^\circ}%
\theta _{\Omega 
{{}^\circ}%
}\in V;$ however, if you take $f=\theta _{\Omega 
{{}^\circ}%
}$ you have that 
\begin{equation*}
\doint \theta _{\Omega 
{{}^\circ}%
}\mu _{\Omega 
{{}^\circ}%
}dx=\int \bar{\chi}_{\Omega }\ dx=m\left( \Omega \right)
\end{equation*}%
and, if we assume that $\partial \Omega $ is smooth, by (\ref{chi}), we have
that 
\begin{equation*}
\doint \theta _{\Omega 
{{}^\circ}%
}\theta _{\Omega 
{{}^\circ}%
}dx=\doint \theta _{\Omega 
{{}^\circ}%
}\bar{\chi}_{\Omega }^{\circ }dx=\doint \theta _{\Omega 
{{}^\circ}%
}dx-\frac{1}{2}\doint \chi _{_{\partial \Omega 
{{}^\circ}%
}}dx=m\left( \Omega \right) -\frac{1}{2}\doint \chi _{_{\partial \Omega 
{{}^\circ}%
}}dx.
\end{equation*}%
Then $\theta _{\Omega 
{{}^\circ}%
}\mu _{\Omega 
{{}^\circ}%
}\neq \theta _{\Omega 
{{}^\circ}%
}\theta _{\Omega 
{{}^\circ}%
}$ and hence $\mu _{\Omega 
{{}^\circ}%
}\neq \theta _{\Omega 
{{}^\circ}%
}.$

Now we can state the following definitions:

\begin{definition}
\label{45}Given an internal set $E\subset \Gamma ,$ we define the \textbf{%
vicinity} of $E$ as follows%
\begin{equation*}
\mathfrak{vic}(E):=\mathfrak{supp}\left( \left\vert D\theta _{E}\right\vert
+\theta _{E}\right) .
\end{equation*}
\end{definition}

The operator $\Omega \longmapsto \mathfrak{vic}(\Omega )$ reminds the
closure operator $\Omega \longmapsto \overline{\Omega }$ but it is not a
closure operator in the topological sense; in fact, in general, $\mathfrak{%
vic}^{2}\left( E\right) :=\mathfrak{vic}\left( \mathfrak{vic}(E)\right)
\supset \mathfrak{vic}(E)$ in a strict sense. However, this similarity with
the closure operator suggests the following

\begin{definition}
\label{46}For any internal set $E\subset \Gamma ,$ we define%
\begin{equation*}
\mathfrak{bd}\left( E\right) :=\mathfrak{supp}\left( \left\vert D\theta
_{E}\right\vert \right)
\end{equation*}%
and we will call $\mathfrak{bd}\left( E\right) $ the \textbf{pointwise
boundary} of $E$ and 
\begin{equation*}
\mathfrak{int}(E):=\left\{ x\in \Gamma \ |\ D\theta _{E}(x)=0\ \text{and}\
\theta _{E}(x)>1\right\} =\mathfrak{vic}(E)\backslash \mathfrak{bd}(E)
\end{equation*}%
and we will call $\mathfrak{int}(E)$ the \textbf{pointwise interior} of $E$.
\end{definition}

If we take $E=\left( \mathbb{R}^{N}\right) 
{{}^\circ}%
=\Gamma $, then $\mathfrak{bd}\left( \Gamma \right) $ is called the \textbf{%
boundary at infinity}.

\subsection{The Gauss' divergence theorem}

In this section we want to generalize the Gauss' divergence theorem in the
framework of the ultrafunction; in particular it is interesting to analyze
the case when $\partial \Omega $ is not smooth.

First, we will examine the smooth case. If $\Omega \subset \mathbb{R}^{N}$
is a bounded set with smooth boundary and $\phi $ is a smooth vector field,
we have that%
\begin{equation*}
\int_{\Omega }\nabla \cdot \phi \ dx=\int_{\partial \Omega }\phi \cdot 
\mathbf{n}_{\Omega }\ dS_{\partial \Omega }
\end{equation*}%
where $S_{\partial \Omega }$ denotes the $(N-1)$-dimensional measure over $%
\partial \Omega $ and $\mathbf{n}_{\Omega }(x)$ is the exterior normal
derivative. For future purposes, we assume that $\mathbf{n}_{\Omega }(x)$ is
a $C_{c}^{1}$-function defined in all $\mathbb{R}^{N}$ which coincides with
the exterior normal derivative in $\partial \Omega .$ We have the following
result:

\begin{theorem}
\label{ruby}Let $\Omega \subset \mathbb{R}^{N}$ be a bounded open set such
that $\bar{\chi}_{\Omega }\in V$; then%
\begin{equation*}
S_{\partial \Omega }^{\circ }(x)=\left\vert D\theta _{\Omega 
{{}^\circ}%
}(x)\right\vert
\end{equation*}%
and $\forall x\in \mathfrak{bd}\left( \Omega \right) $%
\begin{equation*}
\mathbf{n}_{\Omega }(x)=-\frac{D\theta _{\Omega 
{{}^\circ}%
}(x)}{\left\vert D\theta _{\Omega 
{{}^\circ}%
}(x)\right\vert }.
\end{equation*}
\end{theorem}

\textbf{Proof}: The Gauss divergence theorem can be generalized to vector
fields $\phi \in V^{N}$ by writing%
\begin{equation*}
\int \nabla \cdot \phi \ \bar{\chi}_{\Omega }\ dx=\int_{\partial \Omega
}\phi \cdot \mathbf{n}_{\Omega }\ dS_{\partial \Omega }
\end{equation*}%
In this case $\nabla \cdot \phi $ is a measure and we have used the notation
(\ref{nota}). By the transfer principle, we have that $\forall \phi \in
\left( V^{N}\right) ^{\ast }$%
\begin{equation}
\int^{\ast }\nabla ^{\ast }\cdot \phi \ \bar{\chi}_{\Omega }^{\ast }\
dx=\int_{\partial \Omega ^{\ast }}^{\ast }\phi \cdot \mathbf{n}_{\Omega
}^{\ast }\ dS_{\partial \Omega }^{\ast }  \label{bo}
\end{equation}%
By (\ref{A2}), for every $v\in V_{\Lambda },$ 
\begin{equation*}
\int^{\ast }\nabla ^{\ast }\cdot \phi \ v\ dx=\sum_{i=1}^{N}\int^{\ast
}\partial _{i}\phi _{i}\ v\ dx=\sum_{i=1}^{N}\doint D_{i}\phi _{i}\ v%
{{}^\circ}%
\ dx=\doint D\cdot \phi \ v%
{{}^\circ}%
\ dx;
\end{equation*}%
then, by (\ref{pepe}) and (\ref{ip}), we have that 
\begin{equation}
\int^{\ast }\nabla ^{\ast }\cdot \phi \ \bar{\chi}_{\Omega }^{\ast
}dx=\doint D\cdot \phi ~\bar{\chi}_{\Omega }^{\circ }dx=\doint D\cdot \phi
~\theta _{\Omega 
{{}^\circ}%
}dx=-\doint \phi ~D\cdot \theta _{\Omega 
{{}^\circ}%
}dx  \label{ba}
\end{equation}%
Moreover, by (\ref{pollaS})%
\begin{equation}
\int_{\partial \Omega ^{\ast }}^{\ast }\phi ^{\ast }\cdot \mathbf{n}_{\Omega
}^{\ast }dS_{\partial \Omega }=\doint \phi \cdot \mathbf{n}_{\Omega }^{\circ
}\ S_{\partial \Omega }^{\circ }dx;  \label{bi}
\end{equation}%
Then, by (\ref{bo}), (\ref{ba}) and (\ref{bi}), $\forall \phi \in \left( V%
{{}^\circ}%
\right) ^{N}$%
\begin{equation*}
\doint \phi \cdot D\theta _{\Omega 
{{}^\circ}%
}dx=-\doint \phi \cdot \mathbf{n}_{\Omega }^{\circ }\ S_{\partial \Omega
}^{\circ }dx
\end{equation*}%
Now we take $\phi (x)=\delta _{a}(x)\mathbf{e}_{i}$ where $a\in \Gamma $ and 
$\left\{ \mathbf{e}_{i}\right\} _{i=1,..,N}$ is the canonical basis in $%
\mathbb{R}^{N}$ (and hence in $\mathbb{E^{N})}$, we replace $\phi $ in the
above formula:%
\begin{equation*}
\doint \delta _{a}(x)\mathbf{e}_{i}\cdot \mathbf{n}_{\Omega }^{\circ }(x)\
S_{\partial \Omega }(x)dx=-\doint \delta _{a}(x)\mathbf{e}_{i}\cdot D\theta
_{\Omega 
{{}^\circ}%
}(x)dx.
\end{equation*}%
and we get 
\begin{equation}
\left( \mathbf{n}_{\Omega }^{\circ }(a)\cdot \mathbf{e}_{i}\right)
S_{\partial \Omega }(a)=-\mathbf{e}_{i}\cdot D\theta _{\Omega 
{{}^\circ}%
}(a)=-D_{i}\theta _{\Omega 
{{}^\circ}%
}(a);  \label{fg}
\end{equation}%
then%
\begin{equation*}
S_{\partial \Omega }(a)=\sqrt{\dsum\limits_{i=1}^{N}\left[ \left( \mathbf{n}%
_{\Omega }^{\circ }(a)\cdot \mathbf{e}_{i}\right) S_{\partial \Omega }(a)%
\right] ^{2}}=\sqrt{\dsum\limits_{i=1}^{N}\left[ D_{i}\theta _{\Omega 
{{}^\circ}%
}(a)\right] ^{2}}=\left\vert D\theta _{\Omega 
{{}^\circ}%
}(a)\right\vert
\end{equation*}%
Moreover, using again (\ref{fg}) we have that 
\begin{equation*}
\left( \mathbf{n}_{\Omega }^{\circ }(a)\cdot \mathbf{e}_{i}\right)
\left\vert D_{i}\theta _{\Omega 
{{}^\circ}%
}(a)\right\vert =D_{i}\theta _{\Omega 
{{}^\circ}%
}(a);
\end{equation*}%
if $a\in \mathfrak{bd}\left( \Omega \right) ,$ we have that $\left\vert
D_{i}\theta _{\Omega 
{{}^\circ}%
}(a)\right\vert \neq 0$ and hence 
\begin{equation*}
\mathbf{n}_{\Omega }^{\circ }(a)=-\frac{D\theta _{\Omega 
{{}^\circ}%
}(a)}{\left\vert D\theta _{\Omega 
{{}^\circ}%
}(a)\right\vert }\cdot
\end{equation*}

$\square $

\bigskip

Theorem \ref{ruby} suggests the "right" generalization of the Gauss'
theorem; given any set $E\subset \Gamma ,$ we define 
\begin{equation}
\mathbf{n}_{E}(x)=\left\{ 
\begin{array}{cc}
-\frac{D\theta _{E}(x)}{\left\vert D\theta _{E}(x)\right\vert } & if\ \ x\in 
\mathfrak{bd}(E)\ \  \\ 
&  \\ 
0 & if\ \ x\notin \mathfrak{bd}(E)%
\end{array}%
\right.  \label{127}
\end{equation}%
It is amazing that $\mathbf{n}_{E}(x)$ makes sense even if $E$ consists of a
single point $x_{0}.$ Clearly in this case, $\mathfrak{supp}\left( \mathbf{n}%
_{\left\{ x_{0}\right\} }\right) \subset \mathfrak{mon}\left( x_{0}\right) .$

\begin{theorem}
\label{B}(\textbf{Generalized Gauss' divergence theorem}) Let $\phi :\Gamma
\rightarrow \left( V%
{{}^\circ}%
\right) ^{N}$ be a (ultrafunctions) vector field and let $E\subseteq \Gamma $
be an internal set; then%
\begin{equation*}
\doint D\cdot \phi \ d_{E}x=\doint \phi \cdot \mathbf{n}_{E}\ \left\vert
D\theta _{E}\right\vert \ dx
\end{equation*}
\end{theorem}

\textbf{Proof}: By (\ref{polla++}), and (\ref{ip}).%
\begin{eqnarray*}
\doint D\cdot \phi \ d_{E}x &=&\doint D\cdot \phi ~\theta _{E}\ dx=-\doint
\phi \cdot D\theta _{E}~dx \\
&=&-\doint \phi \cdot \frac{D\theta _{E}}{|D\theta _{E}|}~|D\theta
_{E}|~dx=\doint \phi \cdot \mathbf{n}_{E}\ \left\vert D\theta
_{E}\right\vert \ dx.
\end{eqnarray*}

$\square $

\bigskip

It is interesting to compare the above results with the notion of
Caccioppoli set:

\begin{definition}
A Caccioppoli set $\Omega $ is a Borel set such that $\chi _{\Omega }\in BV,$
namely such that $\nabla (\chi _{\Omega })$ (the distributional gradient of $%
\chi _{\Omega }$) is a finite Radon measure. If $\Omega $ is a Caccioppoli
set, then the measure $|\nabla (\chi _{\Omega })|$ is defined as follows: $%
\forall f\in C_{c}^{1},\ f\geq 0$, 
\begin{equation*}
\int f\ d\left( |\nabla (\chi _{\Omega })|\right) :=\sup \left\{
\int_{\Omega }\nabla \cdot \left( f\phi \right) dx\ |\ \phi \in \left(
C^{1}\right) ^{N},\ \left\Vert \phi \right\Vert _{L^{\infty }}\leq 1\right\}
\end{equation*}
\end{definition}

The number 
\begin{equation}
p(\Omega ):=\int d\left( |\nabla (\chi _{\Omega })|\right)  \label{per}
\end{equation}%
is called Caccioppoli perimeter of $\Omega .$ If the reduced boundary of $%
\Omega $ coincides with $\partial \Omega ,$ we have that (see \cite[Section
5.7]{eva-gar}) 
\begin{equation}
\int f(x)\ d\left( |\nabla (\chi _{\Omega })|\right) =\int_{\partial \Omega
}f(x)\,d\mathcal{H}^{N-1}  \label{25}
\end{equation}%
where $\mathcal{H}^{N-1}$ is the $(N-1)$-dimensional Hausdorff measure of $%
\partial \Omega .\ $The Gauss' theorem for a Caccioppoli set takes the
following form:%
\begin{equation*}
\int_{\partial \Omega }\phi \cdot \mathbf{n}_{\Omega }\,d\left( |\nabla
(\chi _{\Omega })|\right) =\int_{\partial \Omega ^{\ast }}\nabla \cdot \phi 
\mathbf{\ }dx
\end{equation*}%
By transfer, we have that for every $\ast $-Borellian vector field $\phi $,%
\begin{equation*}
\int_{\partial \Omega ^{\ast }}^{\ast }\phi \cdot \mathbf{n}_{\Omega }^{\ast
}\,d\left( |\nabla (\chi _{\Omega })|\right) ^{\ast }=\int_{\partial \Omega
^{\ast }}^{\ast }\nabla ^{\ast }\cdot \phi \mathbf{\ }dx
\end{equation*}%
and hence $\forall v\in V_{\Lambda },$ and by (\ref{pollaS}) and Th. \ref%
{ruby}, 
\begin{equation*}
\doint \phi 
{{}^\circ}%
\cdot \mathbf{n}_{\Omega }^{\circ }\ \mu _{|\nabla (\chi _{\Omega })|}\
dx=\doint \phi \cdot \mathbf{n}_{\Omega }\ \left\vert D\theta _{\Omega 
{{}^\circ}%
}\right\vert \ dx
\end{equation*}%
and by the arbitrariness of $\phi ,$ 
\begin{equation*}
\left\vert D\theta _{\Omega 
{{}^\circ}%
}(x)\right\vert =\mu _{|\nabla (\chi _{\Omega })|}
\end{equation*}

However, the measure $\left\vert D\mu _{\Omega }^{\circ }(x)\right\vert $ is
more general than $\mu _{|\nabla (\chi _{\Omega })|}$. For example if $%
\partial \Omega $ is a set with Hausdorff dimension $d>N-1$, $\mu _{|\nabla
(\chi _{\Omega })|}$ is not defined, while $\left\vert D\mu _{\Omega
}^{\circ }(x)\right\vert $ is well defined. Moreover, by the generalization
of (\ref{per}) 
\begin{equation*}
p(\Omega ):=\doint \left\vert D\theta _{\Omega 
{{}^\circ}%
}\right\vert \ dx
\end{equation*}%
$\left\vert D\mu _{\Omega }^{\circ }(x)\right\vert $ allows to define the
perimeter of every set $\Omega \subset \mathbb{R}^{N}\ $even when $\partial
\Omega $ is very wild such as e.g. $\Omega =\partial \Omega $. With this
respect, the theory of fine ultrafunctions is a good improvement of the
theory developed in \cite{BBG}.

\subsection{Ultrafunctions and distributions\label{ud}}

One of the most important properties of the ultrafunctions is that they can
be seen (in some sense that we will make precise in this section) as a
generalizations of the distributions.

\begin{definition}
We say that an ultrafunction $u$ is \textbf{distribution-like} ($DL$) if
there exist a distribution $T$ such that for any $\varphi \in \mathcal{D}$%
\begin{equation*}
\doint u(x)\varphi 
{{}^\circ}%
(x)dx=\left\langle T,\varphi \right\rangle
\end{equation*}%
We say that an ultrafunction $u$ is \textbf{almost} \textbf{distribution like%
} if there exist a distribution $T$ such that for any $\varphi \in \mathcal{D%
}$%
\begin{equation*}
\doint u(x)\varphi 
{{}^\circ}%
(x)dx\sim \left\langle T,\varphi \right\rangle
\end{equation*}
\end{definition}

\textbf{Example: }The measures\textbf{\ } $\left\vert D\mu _{\Omega }^{\circ
}(x)\right\vert \ $and$\ S_{\partial \Omega }^{\circ }(x)$ are
distribution-like since for any $\varphi \in \mathcal{D}$%
\begin{equation*}
\doint \left\vert D\mu _{\Omega }^{\circ }\right\vert \varphi 
{{}^\circ}%
dx=\doint S_{\partial \Omega }^{\circ }\varphi 
{{}^\circ}%
dx=\left\langle T_{S_{\partial \Omega }},\varphi \right\rangle
\end{equation*}%
where $T_{S_{\partial \Omega }}$ is the distribution related to the measure $%
S_{\partial \Omega }$.

It is easy to see that:

\begin{proposition}
\label{P1}For any distribution $T$ there is a distribution-like
ultrafunction $u_{T}$.
\end{proposition}

\textbf{Proof}: Let us consider any projection $P_{\lambda }:V_{\lambda
}\rightarrow V_{\lambda }\cap \mathcal{D}$ and set

\begin{equation*}
\doint u_{T}(x)v(x)dx=\lim_{\lambda \uparrow \Lambda }\ \left\langle
T,P_{\lambda }v_{\lambda }\right\rangle
\end{equation*}%
By Th. \ref{riz}, $u_{T}$ is well defined. So, $\forall \varphi \in \mathcal{%
D}$ 
\begin{equation*}
\doint u_{T}(x)\varphi 
{{}^\circ}%
(x)dx=\lim_{\lambda \uparrow \Lambda }\ \left\langle T,P_{\lambda }\varphi
\right\rangle =\lim_{\lambda \uparrow \Lambda }\ \left\langle T,\varphi
\right\rangle =\left\langle T,\varphi \right\rangle
\end{equation*}

$\square $

\bigskip

Clearly $u_{T}$ is not univocally defined since, in the proof of Prop. \ref%
{P1}, the projection $P_{\lambda }$ can be defined arbitrarily. So it make
sense to set%
\begin{equation*}
\left[ u\right] _{\mathcal{D}^{\prime }}=\left\{ v\in V%
{{}^\circ}%
\ |\ v\approx _{\mathcal{D}^{\prime }}u\right\}
\end{equation*}%
where%
\begin{equation*}
v\approx _{\mathcal{D}^{\prime }}u:\Leftrightarrow \forall \varphi \in 
\mathcal{D},\ \doint \left( u-v\right) \varphi 
{{}^\circ}%
dx=0
\end{equation*}%
Then there is a bijective map 
\begin{equation}
\Psi :\mathcal{D}^{\prime }\rightarrow V_{DL}^{\circ }/\approx _{\mathcal{D}%
^{\prime }}  \label{psi}
\end{equation}%
where $V_{DL}^{\circ }$ is the set of distribution like ultrafunction and%
\begin{equation*}
\Psi (T)=\left\{ u\in V%
{{}^\circ}%
\ |\ \forall \varphi \in \mathcal{D},\ \doint u\varphi 
{{}^\circ}%
dx=\left\langle T,\varphi \right\rangle \right\}
\end{equation*}%
The linear map is $\Psi $ consistent with the distributional derivative,
namely:

\begin{proposition}
If $\Psi (T)=\left[ u\right] _{\mathcal{D}^{\prime }}$ then $\Psi (\partial
_{i}T)=\left[ D_{i}u\right] _{\mathcal{D}^{\prime }}$
\end{proposition}

\textbf{Proof}: If $\Psi (T)=\left[ u\right] _{\mathcal{D}^{\prime }},$ then

\begin{equation*}
\doint D_{i}u\varphi 
{{}^\circ}%
~dx=-\doint uD_{i}\varphi 
{{}^\circ}%
dx
\end{equation*}

Since $\varphi 
{{}^\circ}%
\in U^{\infty },$ then by Th. \ref{TU}.\ref{U3}, $D_{i}\varphi 
{{}^\circ}%
=\left( \partial _{i}\varphi \right) 
{{}^\circ}%
$ and so%
\begin{eqnarray*}
\doint Du\varphi 
{{}^\circ}%
~dx &=&\doint uD\varphi 
{{}^\circ}%
dx=-\doint u\left( \partial _{i}\varphi \right) 
{{}^\circ}%
dx \\
&=&-\left\langle T,\partial _{i}\varphi \right\rangle =\left\langle \partial
_{i}T,\varphi \right\rangle
\end{eqnarray*}%
Hence $\left[ D_{i}u\right] _{\mathcal{D}^{\prime }}=\Psi (\partial _{i}T).$

$\square $

\bigskip

At this point it is a natural question to ask if there exists a linear map 
\begin{equation*}
\Phi :\mathcal{D}^{\prime }\rightarrow V%
{{}^\circ}%
\end{equation*}%
which selects in any equivalence class $\left[ u\right] _{\mathcal{D}%
^{\prime }}$ a distribution-like ultrafunction $\phi (u)$ in a way
consistent with the distributional derivative, namely%
\begin{equation}
\Phi \left( \partial _{i}T\right) =D_{i}\Phi \left( T\right)  \label{mer}
\end{equation}

Actually this goal can be achieved in several ways. We will describe one of
them by using Th. \ref{riz}:

\begin{definition}
\label{cina}For every $T\in \mathcal{D}^{\prime }$, we denote by $T%
{{}^\circ}%
$ the unique ultrafunction in $U^{\infty }$ (see(\ref{Uinf})) such that $%
\forall \psi \in U^{\infty }$ 
\begin{equation}
\doint T%
{{}^\circ}%
(x)\psi (x)dx=\lim_{\lambda \uparrow \Lambda }\ \left\langle T,\psi
_{\lambda }\right\rangle =\left\langle T^{\ast },\psi _{\Lambda
}\right\rangle ^{\ast }  \label{mer1}
\end{equation}
\end{definition}

Clearly, $T%
{{}^\circ}%
$ is a $DL$-ultrafunction since $\forall \varphi \in \mathcal{D}$, $\varphi 
{{}^\circ}%
\in U^{\infty }$ and 
\begin{equation*}
\doint T%
{{}^\circ}%
(x)\varphi 
{{}^\circ}%
dx=\lim_{\lambda \uparrow \Lambda }\ \left\langle T,\varphi \right\rangle
=\left\langle T,\varphi \right\rangle .
\end{equation*}

\begin{theorem}
The map $T\mapsto T%
{{}^\circ}%
$ defined by (\ref{mer1}) satisfies (\ref{mer}), namely%
\begin{equation*}
\left( \partial _{i}T\right) 
{{}^\circ}%
=D_{i}T%
{{}^\circ}%
\end{equation*}
\end{theorem}

\textbf{Proof}: By (\ref{ciccia}) and Th. \ref{TU}.\ref{U4}, we have that $%
\forall \psi \in U^{\infty },$ 
\begin{eqnarray*}
\doint D_{i}T%
{{}^\circ}%
\psi ~dx &=&-\doint T%
{{}^\circ}%
D_{i}\psi ~dx=-\doint T%
{{}^\circ}%
\partial _{i}^{\ast }\psi ~dx \\
&=&-\left\langle T^{\ast },\partial _{i}^{\ast }\psi _{\Lambda
}\right\rangle =\left\langle \partial _{i}^{\ast }T^{\ast },\psi _{\Lambda
}\right\rangle \\
&=&\doint \left( \partial _{i}T\right) 
{{}^\circ}%
\psi ~dx
\end{eqnarray*}%
$\square $

\bigskip

Every function $f\in \mathcal{L}_{loc}^{1}$ defines a distribution $T_{f};$
then, given $f$, we can define three ultrafunctions: $f%
{{}^\circ}%
,$ $\mu _{f}^{\circ }$ and $T_{f}^{\circ }.$ What is the relation between
them? By Prop. \ref{pilla}, we have that $\mu _{f}^{\circ }$ is a projection
of $f^{\ast }$ over $V%
{{}^\circ}%
$, namely%
\begin{equation*}
\mu _{f}^{\circ }(x)=\sum_{a\in \Gamma }\left( \int^{\ast }f^{\ast
}(x)\sigma _{a}(y)dy\right) \delta _{a}(x)
\end{equation*}%
Similarly also $T_{f}^{\circ }$ is a projection, but over a smaller space as
the following proposition shows:

\begin{proposition}
If $f\in \mathcal{L}_{loc}^{1},$ then%
\begin{equation*}
T_{f}^{\circ }=\Pi _{\infty }f%
{{}^\circ}%
.
\end{equation*}%
where $\Pi _{\infty }$ has been defined by (\ref{181}).
\end{proposition}

\textbf{Proof}: By (\ref{ciccia}), we have that $\psi \in U^{\infty },$ 
\begin{equation*}
\doint \Pi _{\infty }f%
{{}^\circ}%
\psi ~dx=\lim_{\lambda \uparrow \Lambda }\int f\psi _{\lambda }\
dx=\lim_{\lambda \uparrow \Lambda }\ \left\langle T,\psi _{\lambda
}\right\rangle =\doint T_{f}^{\circ }\psi ~dx
\end{equation*}%
Since both $T_{f}^{\circ }$ and $\Pi _{\infty }f%
{{}^\circ}%
\in U^{\infty }$, the conclusion follows.

$\square $

\bigskip

So, if $f\in \mathcal{L}_{loc}^{1},$ $T_{f}^{\circ },$ similarly to $\mu
_{f}^{\circ }$, destroys some information contained in $f;$ namely $%
T_{f}^{\circ }$ (resp. $\mu _{f}^{\circ }$) cannot be distiguished by $%
T_{g}^{\circ }$ (resp. $\mu _{f}^{\circ }$) if $f$ and $g$ agree almost
everywere. Similarly, if $\mu $ is any Radon measure and $T_{\mu }$ is the
corresponding distribution, then $T_{\mu }^{\circ }$ destroys some
information contained in $\mu 
{{}^\circ}%
$ since 
\begin{equation*}
T_{\mu }^{\circ }=\Pi _{\infty }\mu 
{{}^\circ}%
\end{equation*}

\textbf{Example}: The $\delta _{a}$ ultrafunction is distribution like since
for every $\varphi \in \mathcal{D}$, we have%
\begin{equation*}
\doint \delta _{a}\varphi 
{{}^\circ}%
(x)dx=\varphi (a)=\left\langle \mathbf{\delta }_{a},\varphi \right\rangle ;
\end{equation*}%
(here we have used the boldface to distinguish the ultrafunction $\delta
_{a} $ from the distribution $\mathbf{\delta }_{a}$). However 
\begin{equation*}
\delta _{a}\neq \mathbf{\delta }_{a}^{\circ }
\end{equation*}%
Actually, according to Def. \ref{sm},%
\begin{equation*}
\delta _{a}=\mathbf{\delta }_{a}^{\circ }+\Pi _{\infty }^{\perp }\delta _{a},
\end{equation*}%
namely $\mathbf{\delta }_{a}^{\circ }\ $is the smooth part of $\delta _{a}.$

\section{Some applications\label{SA}}

In this section we will sketch how the theory of fine ultrafunctions can be
used in the study of Partial Differential Equations. In the framework of
ultrafunctions, a very large class of problems is well posed and has
solutions. Very often, hard \textit{a priori }estimates are not necessary in
proving the existence, but only in understanding the properties of a
solution (qualitative analysis). In particular, if you have a problem from
Physics or from Geometry, it is interesting to investigate whether the
generalized solutions describe the Physical or the Geometric phenomenon. We
refer to \cite{ultra},\cite{BBG},\cite{benciISO},\cite{belu2012},..,\cite%
{bls} where such kind of problems have been treated in the framework of
ultrafunction. \textit{A fortiori}, these problems can be treated using fine
ultrafunctions. In this section, we limit ourselves to give some new
examples just to illustrate the use of fine ultrafunctions with a particular
emphasis in the study of ill posed problems. Obviously, each example is
treated superficially. A deep analysis of each case, probably, would deserve
a full paper.

\subsection{Second order equations in divergence form\label{soe}}

Let $\Omega \subseteq \mathbb{R}^{N}$ be a bounded open set with regular
boundary and let us consider the following boundary value problem:%
\begin{equation}
-\nabla \cdot \left[ k(x,u)\nabla u\right] +f(x,u)=0\ \ \ in\ \ \ \Omega
\label{ab}
\end{equation}%
\begin{equation}
u(x)=0\ \ \ for\ \ \ x\in \partial \Omega  \label{b}
\end{equation}%
where $f$ is a function and $k(x,u)$ is a function or a $\left( N\times
N\right) $-matrix.

A function which satisfies (\ref{ab}) and (\ref{b}) is called classical
solution if $u\in C^{2}\left( \overline{\Omega }\right) .$ The natural
"translation" of this problem in the world of ultrafunctions is the
following:%
\begin{equation}
-D\cdot \left[ k^{\ast }(x,u)Du\right] +f^{\ast }(x,u)=0\ \ \ in\ \ \ \Omega 
{{}^\circ}
\label{b'}
\end{equation}%
\begin{equation}
u(x)=0\ \text{\ }\ if\ \ \ x\in \partial \Omega 
{{}^\circ}
\label{ab'}
\end{equation}

We must be very careful in the interpretation of eq. (\ref{b'}). In fact $Du$
makes sense if $u$ is defined in $\mathfrak{vic}\left( \Omega \right) \ $and 
$D\cdot \left[ k^{\ast }(x,u)Du\right] $ makes sense if $u$ is defined in $%
\mathfrak{vic}^{2}\left( \Omega \right) :=\mathfrak{vic}\left( \mathfrak{vic}%
\left( \Omega \right) \right) .$ Thus the problem (\ref{b'}),(\ref{ab}) is
well posed if 
\begin{equation}
f(x,\cdot )\ \text{and\ }k(x,\cdot )\ \text{are defined in a neighborhood of 
}\Omega .  \label{basic}
\end{equation}

In particular we have that $f^{\ast }(x,\cdot )$, $k^{\ast }(x,\cdot )$ and
consequently $u$ are definite in 
\begin{equation*}
\Omega ^{+}:=\mathfrak{vic}^{2}\left( \Omega \right) .
\end{equation*}%
In conclusion, we are lead to the following definition:

\begin{definition}
A solution of (\ref{b'}), (\ref{ab'}), is the restriction to $\Omega 
{{}^\circ}%
$ of an ultrafunction $u\in V%
{{}^\circ}%
.$ $u|_{\Omega 
{{}^\circ}%
}$ will be called ultrafunction solution of the problem (\ref{ab}), (\ref{b}%
).
\end{definition}

This is the right definition as the following theorem shows:

\begin{theorem}
\label{122}If $w$ is a classical solution of (\ref{ab}),(\ref{b}) then $%
u=w^{\ast }|_{\Omega 
{{}^\circ}%
}$ is a ultrafunction solution.
\end{theorem}

\textbf{Proof}: By our definition of $C^{2}\left( \overline{\Omega }\right)
, $ we have that $w$ can be extended in a neighborhood $\mathcal{N}%
_{\varepsilon }\left( \Omega \right) $ of $\Omega ;\ $consequently also $%
f(x,\cdot )$ and $k(x,\cdot )$ can be extended in $\mathcal{N}_{\varepsilon
}\left( \Omega \right) $ and hence, by the transfer principle, $w^{\ast }$
satisfies the equation%
\begin{equation*}
-\nabla ^{\ast }\cdot \left[ k^{\ast }(x,w^{\ast })\nabla ^{\ast }w^{\ast }%
\right] +f^{\ast }(x,w^{\ast })=0\ \ \ \forall x\in \mathcal{N}\left( \Omega
\right) ^{\ast }
\end{equation*}%
where, with some abuse of notation we have denoted with $f(x,\cdot )$ and $%
k(x,\cdot )$ the extension of the homonymous functions. Then, $u=w^{\ast
}|_{\Omega 
{{}^\circ}%
}$ satisfies (\ref{b'}),(\ref{ab}).

$\square $

\bigskip

If a problem does not have a classical solution, we can look for weak
solutions in some Sobolev space or in a space of distribution. However if
there are not weak solutions, we can find a ultrafunction solution
exploiting the following theorem:

\begin{theorem}
\label{12}Assume (\ref{basic}). If the $\exists M,R\in \mathbb{E}^{+}$ such
that 
\begin{eqnarray}
\left\Vert u\right\Vert &\geq &R\Rightarrow \doint_{\Omega ^{+}}\left[
k^{\ast }(x,u)Du\cdot Dv+f^{\ast }(x,u)u\right] \ dx  \label{pip} \\
&\geq &M\cdot \sqrt{\doint_{\Omega ^{+}}\left\vert u\right\vert ^{2}dx}, 
\notag
\end{eqnarray}%
then problem (\ref{b'}),(\ref{ab'}) has at least one solution.
\end{theorem}

\textbf{Proof}: First of all we set%
\begin{equation}
V_{\text{\textsc{dir}}}^{\circ }\left( \Omega 
{{}^\circ}%
\right) :=\left\{ u\in V%
{{}^\circ}%
\ |\ \forall x\in \partial \Omega 
{{}^\circ}%
,\ u(x)=0\right\} ;  \label{VVV}
\end{equation}%
$V_{\text{\textsc{dir}}}^{\circ }\left( \Omega 
{{}^\circ}%
\right) $ is a hyperfinite Hilbert space equipped with the norm 
\begin{equation*}
\left\Vert u\right\Vert :=\sqrt{\doint \left\vert u\right\vert ^{2}dx}
\end{equation*}%
Now let $\Psi $ be a continuous function such that $\Psi (x)=1$ if $x\in 
\mathcal{N}_{\varepsilon /2}\left( \Omega \right) ,$ $\Psi (x)=0$ if $%
x\notin \mathcal{N}_{\varepsilon }\left( \Omega \right) \ $and $0\leq \Psi
\leq 1.$ We put 
\begin{equation*}
\hat{f}(x,u)=f^{\ast }(x,u)\Psi ^{\ast }(x)
\end{equation*}%
\begin{equation*}
\hat{k}(x,u)=k^{\ast }(x,u)\Psi ^{\ast }(x)+\left[ 1-\Psi ^{\ast }(x)\right]
\end{equation*}%
and we define an operator 
\begin{equation*}
\mathcal{A}:V_{\text{\textsc{dir}}}^{\circ }\left( \Omega 
{{}^\circ}%
\right) \rightarrow V_{\text{\textsc{dir}}}^{\circ }\left( \Omega 
{{}^\circ}%
\right)
\end{equation*}%
by setting, $\forall v\in V_{\text{\textsc{dir}}}^{\circ }\left( \Omega 
{{}^\circ}%
\right) ,$ 
\begin{equation*}
\doint \mathcal{A}\left[ u\right] v\ dx=\doint \left[ \hat{k}(x,u)Du\cdot Dv+%
\hat{f}(x,u)v\right] \ dx.
\end{equation*}%
By (\ref{pip}), we have that $\exists M,R\in \mathbb{E}^{+}$ such that%
\begin{equation*}
\left\Vert u\right\Vert \geq R\Rightarrow \doint \mathcal{A}\left[ u\right]
u\ dx\geq M\left\Vert u\right\Vert ^{2}.
\end{equation*}

Then, by the Brower fixed point theorem and the fact that $V^{\circ }\left(
\Omega 
{{}^\circ}%
\right) $ has hyperfinite dimension the equation $\mathcal{A}\left[ u\right]
=0\ $has at least a solution\ $\bar{u}\in V_{\text{\textsc{dir}}}^{\circ
}\left( \Omega 
{{}^\circ}%
\right) $, and hence $\forall v\in V_{\text{\textsc{dir}}}^{\circ }\left(
\Omega 
{{}^\circ}%
\right) ,$ 
\begin{eqnarray}
0 &=&\doint \mathcal{A}\left[ \bar{u}\right] v\ dx=\doint \left[ \hat{k}(x,%
\bar{u})D\bar{u}\cdot Dv+\hat{f}(x,\bar{u})v\right] \ dx  \notag \\
&=&\doint \left[ -D\cdot \left( \hat{k}(x,\bar{u})D\bar{u}\right) +\hat{f}(x,%
\bar{u})\right] v\ dx  \label{ddc}
\end{eqnarray}%
and hence, for every $x_{0}\in \Omega 
{{}^\circ}%
,$ taking $v=\delta _{x_{0}},$ 
\begin{equation*}
-D\cdot \left( k^{\ast }(x_{0},\bar{u})D\bar{u}\right) +f^{\ast }(x_{0},\bar{%
u})=0.
\end{equation*}%
since for $x_{0}\in \Omega 
{{}^\circ}%
,$ $\hat{k}=k^{\ast }$ and $\hat{f}=f^{\ast }.$ Then $\bar{u}|_{\Omega 
{{}^\circ}%
}$ solves (\ref{b'}),(\ref{ab'}).

$\square $

\begin{remark}
By the above theorem we have that $\bar{u}|_{\Omega 
{{}^\circ}%
}$ might depend on the values that $k^{\ast }$ and $f^{\ast }$ assume in $%
\Omega ^{+}$ and hence, assumption (\ref{basic}) seems essential. However
Th. \ref{122} shows that when $\bar{u}|_{\Omega 
{{}^\circ}%
}$ is sufficiently regular, this dependence disappears.
\end{remark}

The solution of (\ref{b'}),(\ref{ab'}) is the restriction $u|_{\Omega 
{{}^\circ}%
}$ of an ultrafunction $u\in V%
{{}^\circ}%
.$ It is interesting to examine the nature of $u$ near $\partial \Omega 
{{}^\circ}%
.$ By the definition of $V_{\text{\textsc{dir}}}^{\circ }\left( \Omega 
{{}^\circ}%
\right) ,$ we have that 
\begin{equation*}
v\in V%
{{}^\circ}%
\Longleftrightarrow v\chi _{\partial \Omega 
{{}^\circ}%
}\in V_{\text{\textsc{dir}}}^{\circ }\left( \Omega 
{{}^\circ}%
\right)
\end{equation*}%
So, the equation (\ref{ddc}) is equivalent to%
\begin{equation*}
\forall v\in V%
{{}^\circ}%
,\ \doint \left[ \hat{k}(x,u)Du\cdot D\left( v\chi _{\partial \Omega 
{{}^\circ}%
}\right) +f^{\ast }(x,u)\left( v\chi _{\partial \Omega 
{{}^\circ}%
}\right) \right] dx=0
\end{equation*}%
Then, we have that, $\forall v\in V%
{{}^\circ}%
,$

\begin{equation}
\doint \left[ -D\cdot \left[ \hat{k}(x,u)Du\right] \chi _{\partial \Omega 
{{}^\circ}%
}+f^{\ast }(x,u)\chi _{\partial \Omega 
{{}^\circ}%
}\right] vdx=0  \label{pasta}
\end{equation}%
and hence 
\begin{equation}
\left( -D\cdot \left[ \hat{k}(x,u)Du\right] +f^{\ast }(x,u)\right) \chi
_{\partial \Omega 
{{}^\circ}%
}=0.  \label{gigia}
\end{equation}%
If we set%
\begin{eqnarray*}
\left[ V_{\text{\textsc{dir}}}^{\circ }\left( \Omega 
{{}^\circ}%
\right) \right] ^{\perp } &:&=\left\{ \Phi \in V%
{{}^\circ}%
\ |\ \forall v\in V_{\text{\textsc{dir}}}^{\circ }\left( \Omega 
{{}^\circ}%
\right) ,\ \doint \Phi (x)v(x)dx=0\right\} \\
&=&\left\{ \Phi \in V%
{{}^\circ}%
\ |\ \mathfrak{supp}\left( \Phi \right) \subseteq \partial \Omega 
{{}^\circ}%
\right\}
\end{eqnarray*}%
eq. (\ref{gigia}) can be rewritten as follows: 
\begin{equation}
-D\cdot \left[ \hat{k}(x,u)Du\right] +\hat{f}(x,u)=\Phi \left( x\right) \ \
\ with\ \ \ \Phi \in \left[ V_{\text{\textsc{dir}}}^{\circ }\left( \Omega 
{{}^\circ}%
\right) \right] ^{\perp }.  \label{morri}
\end{equation}%
This equation describes in terms of infinitesimal analysis what happens near 
$\partial \Omega 
{{}^\circ}%
.$ $\Phi \left( x\right) $ can be regarded as a force which constrains $u$
to vanish on $\partial \Omega 
{{}^\circ}%
.$ Essentially $\Phi $ describes the reaction of a constraint.

\begin{remark}
\label{RRR}If we want to consider equation (\ref{ab}) with non-homogeneous
boundary condition, i.e. 
\begin{equation*}
u(x)=g(x)\ \ \ for\ \ \ x\in \partial \Omega
\end{equation*}%
we can adopt the standard trick to set 
\begin{equation*}
w(x)=u(x)-\bar{g}(x)
\end{equation*}%
where $\bar{g}(x)$ is any function which extends $g$ in $\Omega $ and to
solve the resulting equation in $w$ with the homogeneous boundary conditions.
\end{remark}

\bigskip

\begin{remark}
The equation (\ref{ab}) with the homogeneous Neumann boundary conditions,
i.e. 
\begin{equation}
\frac{du}{d\mathbf{n}}(x)=\nabla u\cdot \mathbf{n}_{\Omega }=0\ \ \ for\ \ \
x\in \partial \Omega  \label{b+}
\end{equation}%
can be treated in a very similar way. The boundary condition, in the
framework of the ultrafunctions becomes%
\begin{equation}
Du\cdot \mathbf{n}_{E}=0\ \ for\ x\in \partial \Omega 
{{}^\circ}
\label{bb}
\end{equation}%
where $\mathbf{n}_{E}$ has been defined by (\ref{127}). The space $V_{\text{%
\textsc{dir}}}^{\circ }\left( \Omega 
{{}^\circ}%
\right) $ must be replaced by 
\begin{equation*}
V_{\text{\textsc{neu}}}^{\circ }\left( \Omega 
{{}^\circ}%
\right) :=\left\{ u\in V%
{{}^\circ}%
\ |\ \forall x\in \partial \Omega 
{{}^\circ}%
,\ Du\cdot \mathbf{n}_{E}=0\right\} .
\end{equation*}%
Everything else proceeds in a similar way.
\end{remark}

\subsection{Examples}

\textbf{Example 1}: Let us consider the following problem:%
\begin{equation}
-\nabla \cdot \left[ k(x,u)\nabla u\right] +u=f(x)\ \ \ in\ \ \Omega
\label{b1}
\end{equation}%
\begin{equation}
u(x)=0,\ \ for\ \ x\in \partial \Omega  \label{c1}
\end{equation}%
where $k\in C^{1}(\mathbb{R)}$. We set 
\begin{equation*}
k_{0}:=\ \inf \ \left\{ k(x,s)\ |\ \left( x,s\right) \in \Omega \times 
\mathbb{R}\right\}
\end{equation*}%
and%
\begin{equation}
k_{0}>0,  \label{rita}
\end{equation}%
then, if 
\begin{equation}
k\ \text{does\ not\ depend\ on\ }u,  \label{rita1}
\end{equation}%
$A=-\nabla \cdot \left[ k\left( x\right) \nabla u\right] +u$ is a strictly
monotone operator and it is immediate to check that eq. (\ref{b1}) has a
unique weak solution in $H_{0}^{1}(\Omega )\cap L^{4}(\Omega )$ for every $%
f\in H^{-1}(\Omega )+L^{4/3}(\Omega ).$ Moreover, if $f$ and $\partial
\Omega \ $are smooth$,$ by the usual regularity results, problem (\ref{b1}),(%
\ref{c1}) has a classical solution. If (\ref{rita}) or (\ref{rita1}) is not
satisfied, this problem is more delicate. In particular, if for some (but
not all) value of $u$%
\begin{equation*}
k(x,u)<0,
\end{equation*}%
the problem is not well posed and, in general, it has no solution in any
distribution space. Nevertheless we have the following result:

\begin{theorem}
\label{te}If (\ref{basic}) holds and 
\begin{equation}
\underset{\left\vert u\right\vert \rightarrow \infty }{\min \lim }\
k(x,u)\geq 0  \label{bol1}
\end{equation}%
then problem (\ref{b1}), (\ref{c1}) has at least a ultrafunction solution,
namely, $\forall f\in V%
{{}^\circ}%
,$ there exists $u\in V_{\text{\textsc{dir}}}^{\circ }(\Omega ),$ such that%
\begin{equation}
-D\cdot \left[ k^{\ast }(x,u)Du\right] +u=f,\ \ \ \forall x\in \Omega 
{{}^\circ}%
.  \label{can}
\end{equation}
\end{theorem}

\textbf{Proof: }By Th. \ref{12} we have to prove that the operator is (\ref%
{can}) coercive, namely that (\ref{pip}) holds. By (\ref{bol1}) and the
continuity of $k$, there exists a constant $M>0$ such that for $\left\vert
u\right\vert >M,$ 
\begin{equation*}
k(x,u)\geq -\frac{1}{2\left\Vert D\right\Vert ^{2}}\ \ ;\ \ \left\Vert
D\right\Vert =\ \underset{u\neq 0}{\max }\frac{\left\Vert Du\right\Vert }{%
\left\Vert u\right\Vert }
\end{equation*}

Then, setting%
\begin{equation*}
k_{\infty }=\sup \left\{ -k(x,u)\ |\ \left\vert u\right\vert <M\right\}
\end{equation*}%
for any $M\geq 0,$ 
\begin{eqnarray*}
\doint_{\Omega ^{+}}\left[ k(x,u)\left\vert Du\right\vert ^{2}+\left\vert
u\right\vert ^{2}\right] dx &\geq &\doint_{\Omega ^{+},\left\vert
u\right\vert \geq M}k(x,u)\left\vert Du\right\vert ^{2}dx+\doint_{\Omega
^{+},\left\vert u\right\vert <M}k(x,u)\left\vert Du\right\vert
^{2}dx+\left\Vert u\right\Vert ^{2} \\
&\geq &-\frac{1}{2\left\Vert D\right\Vert ^{2}}\doint \left\vert
Du\right\vert ^{2}dx-k_{\infty }\doint_{\Omega ^{+},\left\vert u\right\vert
<M}\left\vert Du\right\vert ^{2}dx+\left\Vert u\right\Vert ^{2} \\
&\geq &-\frac{1}{2\left\Vert D\right\Vert ^{2}}\doint \left\vert
Du\right\vert ^{2}dx-k_{\infty }\underset{u\neq 0,\left\vert u\right\vert
\leq M}{\max }\left\Vert Du\right\Vert +\left\Vert u\right\Vert ^{2} \\
&\geq &-\frac{\left\Vert Du\right\Vert ^{2}}{2\left\Vert D\right\Vert ^{2}}%
-C+\left\Vert u\right\Vert ^{2}\geq \frac{1}{2}\left\Vert u\right\Vert ^{2}-C
\end{eqnarray*}%
$\square $

\bigskip

If $k$ is not positive, our problem might have infinitely many solutions and
they can be quite wild. For example, if we consider the problem 
\begin{equation}
u\in V_{\text{\textsc{dir}}}^{\circ }\left( \Omega 
{{}^\circ}%
\right) ,\ -D\cdot \left[ (u^{3}-u)Du\right] =0.  \label{tcha}
\end{equation}%
we can check directly that, for any internal set $E\subset \mathfrak{int}%
\left( \Omega 
{{}^\circ}%
\right) $ the function $u(x)=\chi _{E}(x)$ is a legitimate solution in the
frame of ultrafunction.

If (\ref{rita}) and (\ref{rita1}) hold, problem (\ref{b1}), (\ref{c1}) has a
classical solution provided that $f$ is sufficiently good. Nevertheless, our
problem has a unique generalized solution even if $f$ is a wild function
(e.g. $f\notin H^{-1}(\Omega )+L^{4/3}(\Omega )$ or $f$ not measurable or $f$
is a distribution not in $H^{-1}$). For example if $N\geq 2$, then $\mathbf{%
\delta }_{a}\notin H^{-1}(\Omega );$ in this case the ultrafunction
solution, for any delta-like $f$, concentrates in $\mathfrak{mon}\left(
a\right) .$

\textbf{Example 2}: Let us consider the following problem:%
\begin{equation*}
u\in C^{2}\left( \overline{\Omega }\right) :
\end{equation*}%
\begin{equation}
-\nabla \cdot \left[ k(u)\nabla u\right] +u^{3}=f(x)\ \ \ in\ \ \Omega
\label{mar1}
\end{equation}%
\begin{equation*}
u(x)=0,\ \ for\ \ x\in \partial \Omega
\end{equation*}%
where $k\ $is a matrix such that%
\begin{equation}
k(u)\xi \cdot \xi \geq -k_{0}\left\vert \xi \right\vert ^{2},\ k_{0}\geq 0.
\label{mar}
\end{equation}

\begin{theorem}
If (\ref{mar}) holds, problem (\ref{b1}), (\ref{c1}) has at least a
ultrafunction solution, namely, $\forall f\in V%
{{}^\circ}%
,$ there exists $u\in V_{\text{\textsc{dir}}}^{\circ }\left( \Omega 
{{}^\circ}%
\right) ,$ such that%
\begin{equation}
-D\cdot \left[ k^{\ast }(u)Du\right] +u^{3}=f\ \ \ \forall x\in \Omega 
{{}^\circ}%
.  \label{can+}
\end{equation}
\end{theorem}

\textbf{Proof: }By Th. \ref{12} we have to prove that the operator is (\ref%
{can}) coercive namely that it satisfies (\ref{pip}). By (\ref{norma}), all
the norms on $V_{\text{\textsc{dir}}}^{\circ }\left( \Omega 
{{}^\circ}%
\right) $ are equivalent; hence, we have that%
\begin{eqnarray*}
\doint \left( -D\cdot \left[ k(u)Du\right] +u^{3}\right) u\ dx
&=&\doint_{\Omega 
{{}^\circ}%
}\left[ k(u)Du\cdot Du+u^{4}\right] dx \\
&\geq &-k_{0}\left\Vert Du\right\Vert ^{2}+\left\Vert u\right\Vert ^{4} \\
&\geq &-k_{0}\left\Vert D\right\Vert ^{2}\left\Vert u\right\Vert
^{2}+\left\Vert u\right\Vert ^{4}.
\end{eqnarray*}%
an hence the operator (\ref{can}) is coercive.

$\square $

\bigskip

As a particular case, we can take 
\begin{equation*}
k=\left[ 
\begin{array}{cc}
-1 & 0 \\ 
0 & c^{2}%
\end{array}%
\right] ;\ 
\end{equation*}%
so equation (\ref{mar1}) reduces to the nonlinear wave equation:%
\begin{equation*}
\square u+u^{3}=f;\ \ \square =\partial _{1}^{2}-c^{2}\partial _{2}^{2}
\end{equation*}%
If we take $\Omega =\left[ 0,T\right] \times \left[ 0,1\right] $ and we
impose periodic boundary conditions in $x_{1}$ this problem reduces to the
classical problem relative to the existence of periodic solution of the
nonlinear wave equation. In general, this problem does not have classical
solutions because of the presence of small divisors which prevent the
approximate solutions to converge. So this is a problem that can be studied
in the framework of ultrafunctions where the existence is guaranteed.

\subsection{Regular weak solutions}

The expression \textit{regular} \textit{weak }sounds like an oxymoron,
nevertheless it well describes the notion we are going to present now.

For example, let us consider the equation (\ref{tcha}). It is possible that
the function of the form $\chi _{E}$ are not acceptable as solutions of a
some physical model described by (\ref{tcha}). Then we may ask if there
exist "approximate" solutions of equation (\ref{tcha}) which exhibit some
form of regularity.

More in general, given the second order operator (\ref{ab}) and (\ref{b}),
we might be interested in regular solutions $u\in U^{m}$ for some $m\in 
\mathbb{N}\cup \left\{ \infty \right\} $ (see section \ref{RU}).

For example, let us consider problem (\ref{ab}),(\ref{b}). We set 
\begin{equation}
U_{\text{\textsc{dir}}}^{1}(\Omega 
{{}^\circ}%
):=\left\{ u\in U^{1}\ |\ \forall x\in \partial \Omega 
{{}^\circ}%
,\ u(x)=0\right\}  \label{Udir}
\end{equation}%
we translate problem (\ref{ab}),(\ref{b}) as follows 
\begin{equation}
u\in U_{\text{\textsc{dir}}}^{1}(\Omega 
{{}^\circ}%
)  \label{w1}
\end{equation}%
such that $\forall v\in U_{\text{\textsc{dir}}}^{1}(\Omega 
{{}^\circ}%
),$%
\begin{equation}
\doint_{\Omega ^{+}}\left[ k^{\ast }(x,u)Du\cdot Dv+f^{\ast }(x,u)v\right]
dx=0  \label{w2}
\end{equation}%
Arguing as in theorem \ref{12}, we have the following result:

\begin{theorem}
If the operator 
\begin{equation*}
\mathcal{A}:U_{\text{\textsc{dir}}}^{1}(\Omega 
{{}^\circ}%
)\rightarrow U_{\text{\textsc{dir}}}^{1}(\Omega 
{{}^\circ}%
)
\end{equation*}%
defined by%
\begin{equation*}
\doint \mathcal{A}\left[ u\right] v~dx=\doint_{\Omega 
{{}^\circ}%
}\left[ k^{\ast }(x,u)Du\cdot Dv+f^{\ast }(x,u)v\right] dx\ \ \forall v\in
U^{1}(\Omega 
{{}^\circ}%
),
\end{equation*}%
is coercive, the equation (\ref{w2}) has at least one solution.
\end{theorem}

In general a solution of (\ref{w1}), (\ref{w2}) does not satisfy equation (%
\ref{b'}), but the equation 
\begin{equation}
-D\cdot \left[ k^{\ast }(x,u)Du\right] +f^{\ast }(x,u)=\psi (x)\ \ \ in\ \ \
\Omega 
{{}^\circ}
\label{w3}
\end{equation}%
where 
\begin{equation}
\psi (x)\in U_{\text{\textsc{dir}}}^{1}(\Omega 
{{}^\circ}%
)^{\perp }=\left\{ w\in U_{\text{\textsc{dir}}}^{1}(\Omega 
{{}^\circ}%
)\ |\ \forall v\in U_{\text{\textsc{dir}}}^{1}(\Omega 
{{}^\circ}%
),\ \doint wv\ dx=0\right\}  \label{w5}
\end{equation}%
can be considered as a sort of error. The error $\psi $ can be considered
negligible since%
\begin{equation*}
\forall v\in U^{1}(\Omega 
{{}^\circ}%
),\ \doint \psi v\ dx=0
\end{equation*}

Probably, in many situation, the regular weak solutions are more relevant
than the solution of type (\ref{b'}).

\begin{remark}
If we adopt the strategy of using regular weak solutions, the more
appropriate functional framework is the use of the quotient space%
\begin{equation*}
\widetilde{U^{1}}:=V%
{{}^\circ}%
/I
\end{equation*}%
where 
\begin{equation*}
I:=U_{\text{\textsc{dir}}}^{1}(\Omega 
{{}^\circ}%
)^{\perp }=\left\{ \psi \in V%
{{}^\circ}%
\ |\ \forall v\in U^{1},\ \doint \psi v\ dx=0\right\}
\end{equation*}%
The space $\widetilde{U^{1}}$ is the analogous, in the world of
ultrafunctions, of the Sobolev space $H^{1};$ in fact in both cases the
functions are not defined pointwise, but they are classes of equivalence of
functions defined up to negligible functions.
\end{remark}

\begin{remark}
\label{bullo}For suitable choices of $k,$ $f$ and $\Omega ,$ it is possible
that the regular weak solutions coincide with the ultrafunctions solutions,
however for ill posed problem in general this fact does not happens.
Clearly, in most cases the regular weak solutions are infinitely close to
ultrafunctions solutions in the appropriate topology (e.g. in the topology
of $H^{1}\left( \Omega \right) ).$ However, if $\partial \Omega $ is very
wild they can differ from each other in a relevant way. Then, in dealing
with a problem relative to a physical phenomenon, the choice of the space in
which to work might be very relevant (see also Remark \ref{reg}).
\end{remark}

\subsection{Calculus of variations}

Let us consider the minimization problem of the functional%
\begin{equation}
J(u)=\int_{\Omega }F(x,u,\nabla u)\ dx  \label{J}
\end{equation}%
If we assume the Dirichlet boundary condition, the natural space where to
work is $C^{1}(\Omega )\cap C_{0}^{0}(\overline{\Omega })$. If we translate
this problem in the framework of ultrafunctions the natural space is $V_{%
\text{\textsc{dir}}}^{1}(\Omega 
{{}^\circ}%
)$ defined by (\ref{VVV}) and the condition (\ref{basic}) is translated in
the assumption that $F(x,u,\xi )$ be defined in $\mathcal{N}_{\varepsilon
}\left( \Omega \right) \times \mathbb{R}\times \mathbb{R}^{N}.$ Then the
functional (\ref{J}) becomes%
\begin{equation}
J%
{{}^\circ}%
(u):=\doint F^{\ast }(x,u,Du)\ d_{\Omega }x.  \label{JJ}
\end{equation}%
and we have the following result:

\begin{theorem}
\label{nadia}If (\ref{JJ}) is coercive, then $J%
{{}^\circ}%
(u)$ has a minimizer $\bar{u}$ in $V_{\text{\textsc{dir}}}^{\circ }\left(
\Omega 
{{}^\circ}%
\right) .$ Moreover, if $J(u)$ has a minimizer $w\in C^{1}(\overline{\Omega }%
),$ then, $\forall x\in \overline{\Omega },$ 
\begin{equation*}
\bar{u}\left( x\right) =w^{\ast }(x).
\end{equation*}
\end{theorem}

\textbf{Proof:} Trivial.

$\square $

Also in this case it is interesting to see the form assumed by The
Euler-Lagrange equations. For simplicity, we assume that 
\begin{equation*}
F^{\ast }(x,u,Du)=\frac{1}{2}k^{\ast }(x,u)\left\vert Du\right\vert
^{2}+h^{\ast }(x,u)
\end{equation*}%
then, 
\begin{equation*}
dJ%
{{}^\circ}%
(u)\left[ v\right] :=\doint \left[ k^{\ast }(x,u)Du\cdot Dv+f(x,u)v\right] \
d_{\Omega }x;\ f(x,u)=\frac{\partial h(x,u)}{\partial u}.
\end{equation*}%
By (\ref{polla}), we have that, $\forall v\in V_{\text{\textsc{dir}}%
}^{1}(\Omega 
{{}^\circ}%
)$ 
\begin{equation*}
dJ%
{{}^\circ}%
(u)\left[ v\right] :=\doint \left[ k^{\ast }(x,u)Du\cdot Dv+f(x,u)v\right]
\mu _{\Omega }^{\circ }(x)\ dx.
\end{equation*}%
and arguing as we did at the end of section \ref{soe} (see eq. (\ref{morri}%
)) we get the equation%
\begin{equation}
-D\cdot \left[ \mu _{\Omega }^{\circ }(x)k^{\ast }(x,u)Du\right] +f^{\ast
}(x,u)\mu _{\Omega }^{\circ }(x)=\Phi \left( x\right) \ \ \ with\ \ \ \Phi
\in \left[ V_{\text{\textsc{dir}}}^{\circ }\left( \Omega 
{{}^\circ}%
\right) \right] ^{\perp }.  \label{cone}
\end{equation}

It is interesting to note that eqs. (\ref{morri}) and (\ref{cone}) are
slightly different; they coincide for the $x\in \Omega 
{{}^\circ}%
$ when $mon(x)\subset \Omega 
{{}^\circ}%
,$ but they differ for some $x\sim \partial \Omega 
{{}^\circ}%
.$ This fact depend on the fact that in the framework of the ultrafunction a
weak solution formulated with the measure $\mu _{\Omega }^{\circ }$ does not
coincide with a solution defined by (\ref{b'}),(\ref{ab'}). In any case, we
have seen that the regular solutions are the same.

The minimization problem in $C^{1}(\overline{\Omega }),$ with no constraint
on the boundary, in the regular case leads to equation (\ref{ab}) with the
Neumann boundary conditions; in the world of ultrafunctions, this problem
gives to the equation%
\begin{equation}
-D\cdot \left[ \mu _{\Omega }^{\circ }(x)k^{\ast }(x,u)Du\right] +f^{\ast
}(x,u)\mu _{\Omega }^{\circ }(x)=0.  \label{ung}
\end{equation}%
Thus the boundary conditions are included in this equation.

\bigskip

When we deal with the calculus of variations the notion of regular weak
solutions arises in a natural way. In fact, it makes perfect sense to
minimize the functional (\ref{JJ}) in the spaces%
\begin{equation*}
U_{\text{\textsc{dir}}}^{m}\left( \Omega 
{{}^\circ}%
\right) :=V_{\text{\textsc{dir}}}^{\circ }\left( \Omega 
{{}^\circ}%
\right) \cap U^{m};\ m\geq 0.
\end{equation*}%
If $J%
{{}^\circ}%
$ is coercive, the minimizer exists in each space and if the minimizer in $%
V_{\text{\textsc{dir}}}^{\circ }\left( \Omega 
{{}^\circ}%
\right) $ is not $m$-regular it is different from the minimizer in $U_{\text{%
\textsc{dir}}}^{m}\left( \Omega 
{{}^\circ}%
\right) $; this phenomenon is typical of the ultrafunctions and it does not
have any analogous thing in the framework of $C^{m}$-functions or in the
Sobolev spaces.

\bigskip

\textbf{Example}: Let us consider the functional defined in $C^{1}(0,1)\cap
C_{0}^{0}(\left[ 0,1\right] )$ 
\begin{equation}
J(u)=\int_{\Omega }\left[ \left( \left\vert \partial u\right\vert
^{2}-1\right) ^{2}+u^{2}\right] \ dx  \label{JL}
\end{equation}%
which presents the well known Lavrentiev phenomenon, namely every minimizing
sequence $u_{n}$ converges uniformly to $0$ but%
\begin{equation*}
0=\ \underset{n\rightarrow \infty }{\lim }J(u_{n})\neq J(0)=1.
\end{equation*}%
If we use Th. \ref{nadia}, the minimizing sequence $u_{\lambda }\in
U_{\lambda }^{1}(\Omega )$ has an infinitesimal $\Lambda $-limit $\bar{u}\in
U_{\text{\textsc{dir}}}^{1}\left( \Omega 
{{}^\circ}%
\right) $ and $J%
{{}^\circ}%
(\bar{u})$ is a positive infinitesimal which satisfies the following
Euler-Lagrange equations in $\mathfrak{int}(\Omega 
{{}^\circ}%
)$:%
\begin{equation}
D\left[ \left( \left\vert D\bar{u}\right\vert ^{2}-1\right) D\bar{u}\right] -%
\frac{1}{2}\bar{u}=\psi  \label{glu}
\end{equation}%
where%
\begin{equation*}
\forall v\in U^{1}(\Omega )^{\perp },\ \doint \psi v\ dx=0.
\end{equation*}%
However, it is possible to minimize the functional%
\begin{equation*}
J%
{{}^\circ}%
(u):=\int \left[ \left( \left\vert Du\right\vert ^{2}-1\right) ^{2}+u^{2}%
\right] \ dx
\end{equation*}%
in $V_{\text{\textsc{dir}}}^{\circ }\left( \Omega 
{{}^\circ}%
\right) ;$ in this case we get a minimizer $\hat{u}$ such that $J%
{{}^\circ}%
(\hat{u})<J%
{{}^\circ}%
(\bar{u})$ and in $\mathfrak{int}(\Omega 
{{}^\circ}%
)$ we have: 
\begin{equation*}
D\left[ \left( \left\vert D\hat{u}\right\vert ^{2}-1\right) D\hat{u}\right] -%
\frac{1}{2}\hat{u}=0.
\end{equation*}%
The function $\psi $ in equation (\ref{glu}) can be interpreted ad a
structural "force" which prevents $u$ to form "angles". The presence of this
force, increases the energy level and we have that 
\begin{equation*}
J%
{{}^\circ}%
(\bar{u})>J%
{{}^\circ}%
(\hat{u})>0
\end{equation*}%
even if $J%
{{}^\circ}%
(\bar{u})$ remains infinitesimal.

Also this example shows how the choice of the space where to work changes
the solution of the problem and hence, as observed in Remarks \ref{bullo}
and \ref{reg}, the choice of a particular space depends on the phenomenon
which we want to describe.

\begin{remark}
We have assumed $F(x,u,\xi )$ to be continuous in $u\ $and $\xi $; however
also the case of discontinuous functions can be easily analyzed (see e.g. 
\cite{BBG}).
\end{remark}

\subsection{Evolution problems\label{ep}}

Let $\Omega \subset \mathbb{R}^{N}$ be an open set and let%
\begin{equation*}
A(x,\partial _{i}):D_{A}\left( \Omega \right) \rightarrow C\left( \Omega
\right) ;\ 
\end{equation*}%
be a differential operator. Here $D_{A}\left( \Omega \right) $denotes the
domain of $A(x,\partial _{i}).$

We consider the following Cauchy problem: given $u_{0}(x)\in C^{0}\left( 
\overline{\Omega }\right) ,$ find%
\begin{equation}
u\in C^{1}(I,D_{A}\left( \Omega \right) )\cap C^{0}(I,C\left( \Omega \right)
),\ 0\in I\subseteq \mathbb{R}:  \label{ae}
\end{equation}%
\begin{equation}
\partial _{t}u=A(x,\partial _{i})\left[ u\right]  \label{be}
\end{equation}%
\begin{equation}
u(0,x)=u_{0}(x).  \label{ce}
\end{equation}

A function which satisfies (\ref{ae}), (\ref{be}), (\ref{ce}) is called
classical solution. We want to translate this problem in the world of
ultrafunctions. Because of the nature of the Cauchy problem, it is not
convenient to translate this problem in the space $V%
{{}^\circ}%
\left( \mathbb{R}\times \Omega \right) .$ It is better to use the internal
space of functions $C^{1}(\mathbb{E},V^{\circ })$ defined in section \ref%
{TDU}. We assume that the boundary conditions which are "contained" in the
definition of the domain $D_{A}\left( \Omega \right) $ can be translated in
a domain $D_{A}^{\circ }\left( \Omega \right) \subset V%
{{}^\circ}%
$ (in sections \ref{soe} we have seen how this can be done for second order
operators with Dirichlet and Neumann boundary conditions). Then, setting 
\begin{equation*}
C^{1}(I^{\ast },D_{A}^{\circ }\left( \Omega \right) )=\left\{ u\in C^{1}(%
\mathbb{E},V%
{{}^\circ}%
)\ |\ \forall t\in I^{\ast },\ u(t,x)\in D_{A}^{\circ }\left( \Omega \right)
\right\}
\end{equation*}

The problem (\ref{ae}), (\ref{be}), (\ref{ce}) translates in the following
one: 
\begin{equation}
u\in C^{1}(I,D_{A}^{\circ }\left( \Omega \right) )  \label{ae'}
\end{equation}%
\begin{equation}
\partial _{t}^{\ast }u=A^{\ast }(x,D_{i})\left[ u\right] ,\ t\in I^{\ast }
\label{be'}
\end{equation}%
\begin{equation}
u(0,x)=u_{0}(x)  \label{ce'}
\end{equation}

Remember that the time derivative $\partial _{t}^{\ast }$ is not the
generalized derivative, but the natural extension of $\partial _{t}$ defined
by (\ref{dt}).

A solution of (\ref{ae'}), (\ref{be'}), (\ref{ce'}) will be called
ultrafunction solution of the problem (\ref{ae}), (\ref{be}), (\ref{ce}).

Using this notation, we can state the following fact:

\begin{theorem}
\label{bott}If $w$ is a classical solution of (\ref{ae}), (\ref{be}), (\ref%
{ce}) and we assume that $w$ be extended continuously in a neighborhood $%
\mathcal{N}\left( \Omega \right) $; then $u=w%
{{}^\circ}%
$ is a ultrafunction solution of (\ref{ae'}), (\ref{be'}), (\ref{ce'}).
\end{theorem}

\textbf{Proof}: It is the same than the proof of Th. \ref{122}.

$\square $

\bigskip

Also in the evolution case, the conditions which guarantee the existence of
a ultrafunction solution are very weak.

\begin{theorem}
\label{evol}Assume that $A(x,\partial _{i})\left[ u\right] $ restricted to $%
V_{\lambda }\left( \Omega \right) $ is locally Lipschitz continuous in $u;$
then there exists $T_{\Lambda }$ such that the problem (\ref{ae'}), (\ref%
{be'}), (\ref{ce'}) has a unique ultrafunction solution for $t\in \left[
0,T_{\Lambda }\right) .$ Moreover, if there is an a prori bound for such a
solution , then there exists a unique ultrafunction solution in $%
C^{1}(I^{\ast },D_{A}^{\circ }\left( \Omega \right) )$.
\end{theorem}

\textbf{Proof}. For every $\lambda $ consider the following system of $%
ODE^{\prime }s$ in $\mathfrak{F}\left( \Gamma _{\lambda }\cap \Omega \right) 
$%
\begin{equation}
\partial _{t}u_{\lambda }(t,a_{\lambda })=A(x,D_{i\lambda })\left[
u_{\lambda }\right] (t,a_{\lambda }),\ a_{\lambda }\in \Gamma _{\lambda }
\label{ru}
\end{equation}%
Since $\mathfrak{F}\left( \Gamma _{\lambda }\right) $ is a finite
dimensional vector space and $A(x,D_{i\lambda })\left[ u\right] $ restricted
to $V_{\lambda }\left( \Omega \right) $ is locally Lipschitz continuous the
above system has a local solution for $t\in \left[ 0,T_{\lambda }\right) .$
Then, setting 
\begin{equation*}
T_{\Lambda }=~\lim_{\lambda \uparrow \Lambda }~T_{\lambda }
\end{equation*}%
and taking the $\Lambda $-limit in (\ref{ru}), $\forall t\in \left[
0,T_{\Lambda }\right) ,$ we get 
\begin{equation*}
\partial _{t}^{\ast }u(t,a)=A^{\ast }(x,D_{i})\left[ u\right] (t,a)
\end{equation*}%
Moreover, if there is an \textit{a priori} bound for such a solution, then
there is a bound for the approximate solution (\ref{ru}) in $\left[
0,T_{\lambda }\right) $ (which might also depend on $\lambda $); then, it is
well known that $\left[ 0,T_{\lambda }\right) =I$ and hence $\left[
0,T_{\Lambda }\right) =I^{\ast }$

$\square $

\subsection{Some examples of evolution problems}

\textbf{Example 1}: Let $\Omega $ be a bounded open set and let us consider
the following problem:%
\begin{equation}
u\in C^{1}(I,C_{0}^{0}\left( \overline{\Omega }\right) )\cap
C^{0}(I,C^{2}\left( \Omega \right) ):
\end{equation}%
\begin{equation}
\partial _{t}u=\nabla \cdot \left[ k(u)\nabla u\right]  \label{liliana}
\end{equation}%
\begin{equation}
u(0,x)=u_{0}(x)
\end{equation}

If $k(u)$ satisfies (\ref{rita}), then (\ref{liliana}) is a parabolic
equation and the problem, under suitable condition, has a classical
solution. If $k(u)<0$ for some $u\in \mathbb{R}$, then the problem is
ill-posed, and in general classical solutions do not exist.

Let us translate this problem in the world of ultrafunctions; taking account
of the results of section \ref{soe}, we get:%
\begin{equation}
u\in C^{1}(I^{\ast },V_{\text{\textsc{dir}}}^{\circ }\left( \Omega 
{{}^\circ}%
\right) )  \label{aa1}
\end{equation}%
\begin{equation}
\partial _{t}^{\ast }u=D\cdot \left[ k^{\ast }(u)Du\right] \ \ for\ \ x\in
\Omega 
{{}^\circ}
\label{aa2}
\end{equation}%
\begin{equation}
u(0,x)=u_{0}^{\circ }(x)  \label{aa3}
\end{equation}%
So we can apply the theorems \ref{evol} and we get the existence of a unique
(local in time) solution. It is not difficult to get sufficient conditions
which guarantee the existence of a solution for every $t\in I^{\ast }.$ For
example

\begin{theorem}
If the set%
\begin{equation}
B=\left\{ r\in \mathbb{R\ }|\ k(r)<0\right\}  \label{billo}
\end{equation}%
is bonded, problem (\ref{aa1}),(\ref{aa2}),(\ref{aa3}) has a global solution.
\end{theorem}

\textbf{Proof} - We set $k(r)=k^{+}(r)-k^{-}(r)$ and with some abuse of
notation we write, for $r\in \mathbb{E}$, $k^{\ast }(r)=k^{+}(r)-k^{-}(r).$
Then, we have that%
\begin{eqnarray*}
\partial _{t}^{\ast }\doint u^{2}d_{\Omega }x &=&2\doint u\partial
_{t}^{\ast }u\ d_{\Omega }x=2\doint uD\cdot \left[ k^{\ast }(u)Du\right]
d_{\Omega }x \\
&=&-2\doint k^{\ast }(u)\left\vert Du\right\vert ^{2}d_{\Omega }x\leq
2\doint_{\Omega 
{{}^\circ}%
}k^{-}(u)\left\vert Du\right\vert ^{2}d_{\Omega }x \\
&\leq &2\left\Vert D\right\Vert ^{2}\doint_{\Omega 
{{}^\circ}%
}k^{-}(u)\left\vert u\right\vert ^{2}d_{\Omega }x
\end{eqnarray*}%
If we set 
\begin{equation*}
M=\ \underset{r\in \mathbb{R}}{\max }\ k^{-}(r)\left\vert r\right\vert ^{2}
\end{equation*}%
then, we have that%
\begin{equation*}
\partial _{t}^{\ast }\doint u^{2}d_{\Omega }x\leq 2M\left\Vert D\right\Vert
^{2}\cdot \doint d_{\Omega }x
\end{equation*}%
and this implies the existence of a global solution.

$\square $

\bigskip

In many applications of eq. (\ref{liliana}) $u$ represent a density and $%
k(u)\nabla u$ its flow. Then, thanks to the generalized Gauss' theorem \ref%
{B}, it is easy to prove that the "mass" $\doint_{\Omega }u(x)dx$ is
preserved up to the flow crossing $\partial \Omega $:%
\begin{eqnarray*}
\partial _{t}^{\ast }\doint u(x)d_{\Omega }x &=&\doint \partial _{t}^{\ast
}u(x)d_{\Omega }x \\
&=&\doint D\cdot \left[ k^{\ast }(u)Du\right] d_{\Omega }x \\
&=&\doint k^{\ast }(u)Du\cdot \mathbf{n}_{\Omega }\ \left\vert D\theta
_{E}\right\vert \ dx
\end{eqnarray*}%
If we want to model a situation where the flow of $u$ cannot cross $\partial
\Omega $, $\Omega $ bounded, in a classical context, the Neumann boundary
conditions are imposed: 
\begin{equation}
\forall x\in \partial \Omega ,\ \frac{du}{d\mathbf{n}}\left( x\right) =0;
\label{ritina}
\end{equation}%
in the framework of ultrafunction this situation can be easily described
using the analog of eq. (\ref{ung}) for equation (\ref{liliana}):%
\begin{equation}
\partial _{t}^{\ast }u=D\cdot \left[ \mu _{\Omega }^{\circ }(x)k^{\ast
}(x,u)Du\right]
\end{equation}

In this case, we have the flow $\mu _{\Omega }^{\circ }(x)k^{\ast }(u)Du$
vanishes out of $\mathfrak{vic}^{2}\left( \Omega 
{{}^\circ}%
\right) $ and we have the conservation of the "mass" $\doint u(x)dx$. We can
prove this fact directly, in fact, since $\Omega $ is bounded, $\forall x\in 
\mathfrak{vic}\left( \Omega 
{{}^\circ}%
\right) ,\ D1%
{{}^\circ}%
=0;$ hence 
\begin{eqnarray*}
\partial _{t}^{\ast }\doint u(x)dx &=&\doint D\cdot \left[ \mu _{\Omega
}^{\circ }k^{\ast }(u)Du\right] dx \\
&=&\doint D\cdot \left[ \mu _{\Omega }^{\circ }k^{\ast }(u)Du\right] 1%
{{}^\circ}%
dx \\
&=&\doint k^{\ast }(u)\mu _{\Omega }^{\circ }Du\cdot D1%
{{}^\circ}%
dx=0.
\end{eqnarray*}

\begin{remark}
This problem has been studied with nonstandard methods by Bottazzi in the
framework of grid functions \cite{Bott}. One of the main differences is that
in the context of \cite{Bott}, theorem \ref{bott} does not hold.
\end{remark}

\textbf{Example 2.} Now let us consider the following "conservation law": 
\begin{equation}
u\in C^{1}(I,C^{0}\left( \Omega \right) )\cap C^{0}(I,C^{1}\left( \Omega
\right) ):  \label{x}
\end{equation}%
\begin{equation}
\partial _{t}u=\nabla \cdot F(x,u)  \label{y}
\end{equation}%
\begin{equation}
u(0,x)=u_{0}(x)  \label{z}
\end{equation}%
where $F:\mathbb{R}^{N+1}\rightarrow \mathbb{R}^{N}$ is a smooth function
with support in $\overline{\Omega }$ as well as the initial data $u_{0}(x).$
This problem is not well posed and when $N>1$ and very little is known.
Nevertheless this problem is well posed in the frame of ultrafunctions:%
\begin{equation}
u\in C^{1}(I,V%
{{}^\circ}%
\left( \mathfrak{vic}\left( \Omega 
{{}^\circ}%
\right) \right) ):
\end{equation}%
\begin{equation}
\partial _{t}^{\ast }u=D\cdot F^{\ast }(x,u)  \label{hyden}
\end{equation}%
\begin{equation}
u(0,x)=u_{0}^{\circ }(x)
\end{equation}

\begin{theorem}
If 
\begin{equation}
\left\vert \frac{\partial F}{\partial u}\right\vert \leq
c_{1}+c_{2}\left\vert u\right\vert ,  \label{cielo}
\end{equation}
problem (\ref{x}),(\ref{y}),(\ref{z}) has a unique (global in time)
ultrafunction solution $u(t,x)$ which satisfies the following properties: 
\begin{equation}
\mathfrak{supp}\left( u(t,x)\right) \subset \mathfrak{vic}\left( \Omega 
{{}^\circ}%
\right)  \label{dolores}
\end{equation}%
\begin{equation}
\partial _{t}^{\ast }\doint u\ dx=0  \label{elisa}
\end{equation}
\end{theorem}

\textbf{Proof}: The existence of a unique solution follows from Th. \ref%
{evol}. (\ref{dolores}) follows from the fact that for every $x\notin 
\mathfrak{vic}\left( \Omega 
{{}^\circ}%
\right) ,$ $D\cdot F(x,u)=0.$ (\ref{elisa}) is an immediate consequence of
the Gauss' theorem \ref{B}.

$\square $

\bigskip

A particular case of eq. (\ref{hyden}) is the Burger's equation%
\begin{equation*}
\partial _{t}u=-u\partial _{x}u
\end{equation*}%
\begin{equation*}
u(0,x)=u_{0}(x)\geq 0.
\end{equation*}%
It is well known that this equation has infinitely many weak solutions which
preserve the mass. One of them, the entropy solution, describes the
phenomena occurring in fluid mechanics. It is convenient to write the
Burger's equation in the framework of ultrafunctions as follows:%
\begin{equation}
u\in C^{1}(\mathbb{E},V%
{{}^\circ}%
)  \label{blua}
\end{equation}%
\begin{equation}
\partial _{t}^{\ast }u=-D_{x}\left( \frac{u\left\vert u\right\vert }{2}%
\right)  \label{blue}
\end{equation}%
\begin{equation}
u(0,x)=u_{0}^{\circ }(x)  \label{bluf}
\end{equation}

Since the right hand side of this equation does not satisfy (\ref{cielo}),
the formulation (\ref{blue}) grants \textit{a priori} bounds and hence the
existence of a global solutions. In fact, by (\ref{A2}), 
\begin{equation*}
\partial _{t}^{\ast }\doint \left\vert u\right\vert ^{3}dx=\frac{1}{3}\doint
\left( u\left\vert u\right\vert \right) \partial _{t}u\ dx=-\frac{1}{6}%
\doint \left( u\left\vert u\right\vert \right) D_{x}\left( u\left\vert
u\right\vert \right) dx=0
\end{equation*}%
The solutions of eq. (\ref{blue}) are different from the entropy solution
since they preserves also the quantity $\doint \left\vert u\right\vert
^{3}dx.$ The viscosity solution can be modelled in the frame of
ultrafunctions by the equation%
\begin{equation*}
\partial _{t}^{\ast }u=-D_{x}\left( \frac{u^{2}}{2}\right) +\nu D_{x}^{2}u
\end{equation*}%
where $\nu $ is a suitable infinitesimal (see \cite{blbur}).

Finally we remark that in the frame of ultrafunctions the equation (\ref%
{blue}) is different from the equation%
\begin{equation}
\partial _{t}^{\ast }u=-uD_{x}u  \label{bue+}
\end{equation}%
even if $u\geq 0;$ in fact, in the point where $u$ is singular $D_{x}\left(
u^{2}\right) \neq 2uD_{x}u$ and \textit{a priori }eq. (\ref{bue+}) is not a
conservation law since it does not have the form (\ref{hyden}). Nevertheless
the mass is preserved since, by (\ref{A2}), 
\begin{equation*}
\partial _{t}\doint u\ dx=-\doint uD_{x}u\ dx=0
\end{equation*}%
It is is immediate to see that $u_{0}^{\circ }(x)\geq 0$ implies that, $%
\forall t,\ u(t,x)\geq 0$ and hence, by the fact that $\doint \left\vert
u\right\vert \ dx$ is constant, we get an a priori bound and the existence
of a global solution.

Moreover, if $\mathfrak{supp}\left( u_{0}\right) \subset \left[ a,b\right] 
{{}^\circ}%
$ then $\mathfrak{supp}\left( u\left( t,\cdot \right) \right) \subset 
\mathfrak{vic}\left( \left[ a,b\right] 
{{}^\circ}%
\right) $ since for every $x\notin \mathfrak{vic}\left( \left[ a,b\right] 
{{}^\circ}%
\right) ,$ $\partial _{t}u(x)=0\ $for every $t\in \mathbb{E}$. So a new
phenomenon occurs: the mass concentrates in the front of the shock waves.
For example consider the initial condition%
\begin{equation*}
u_{0}(x)=\left\{ 
\begin{array}{cc}
x & if\ \ x\in \left[ 0,1\right] 
{{}^\circ}
\\ 
&  \\ 
0 & if\ \ x\in \Gamma \backslash \left[ 0,1\right] 
{{}^\circ}%
\end{array}%
\right.
\end{equation*}%
In this case the solution, for $t\geq 0,$ is 
\begin{equation*}
u(t,x)=\left\{ 
\begin{array}{cc}
\frac{x}{t+1} & if\ \ x\in \left[ 0,1\right] 
{{}^\circ}%
\backslash \mathfrak{mon}\left( 1\right) \\ 
&  \\ 
0 & if\ \ x\leq 0\ \ or\ x>0.\ 
\end{array}%
\right.
\end{equation*}%
and hence, as $t\rightarrow \infty ,$ all the mass concentrates in $%
\mathfrak{mon}\left( 1\right) .$

In conclusion, the translation of the Burger's equation in the frame of
ultrafuctions, leads to several different situation which might reproduce
different physical models.

\bigskip

\textbf{Example 3.} Let us consider the following problem:%
\begin{equation*}
u\in C^{2}(I,C^{0})\cap C^{0}(I,C^{2})
\end{equation*}%
\begin{equation*}
\square u+\left\vert u\right\vert ^{p-2}u=0\ \ \ in\ \ \ \mathbb{R}^{N};\ p>2
\end{equation*}%
\begin{eqnarray*}
u(0,x) &=&u_{0}(x), \\
\ \partial _{t}u(0,x) &=&u_{1}(x);
\end{eqnarray*}%
where%
\begin{equation*}
\square u=\partial _{t}^{2}u-\Delta u.
\end{equation*}%
For simplicity we assume that $u_{0}(x)$ and $u_{1}(x)$ have compact support.

It is well known that this problem has a unique weak solution in suitable
Sobolev spaces, provided that, for $N\geq 3,$%
\begin{equation*}
p\leq 2^{\ast }=\frac{2N}{N-2}
\end{equation*}%
and, in this case, the energy%
\begin{equation*}
E(u(t,\cdot ))=\int \left( \frac{1}{2}\left\vert \partial _{t}^{\ast
}u\right\vert ^{2}+\frac{1}{2}\left\vert \nabla u\right\vert ^{2}+\frac{1}{p}%
\left\vert u\right\vert ^{p}\right) dx
\end{equation*}%
is a constant of motion.

If $p>2^{\ast }$, the problem has a weak solution but it is an open question
if it is unique and if the energy is preserved. This problem can be
translated and generalized in the framework of ultrafunctions:%
\begin{equation}
u\in C^{2}(I,V^{\circ })  \label{20}
\end{equation}%
\begin{equation}
\partial _{t}^{\ast 2}u-D^{2}u+\left\vert u\right\vert ^{p-2}u=0\ \ \ for\ \
x\in \Gamma  \label{21}
\end{equation}%
\begin{equation}
u(0,x)=u_{0}^{\circ }(x);\ \partial _{t}^{\ast }u(0,x)=u_{1}^{\circ }(x).
\label{23}
\end{equation}

\begin{theorem}
\label{T1}The problem (\ref{20}),(\ref{21}),(\ref{23}) has a unique
solution. Moreover the energy equality holds, namely%
\begin{equation*}
\partial _{t}^{\ast }E(u(t,\cdot ))=\partial _{t}^{\ast }\doint \left( \frac{%
1}{2}\left\vert \partial _{t}^{\ast }u\right\vert ^{2}+\frac{1}{2}\left\vert
Du\right\vert ^{2}+\frac{1}{p}\left\vert u\right\vert ^{p}\right) dx=0
\end{equation*}
\end{theorem}

\textbf{Proof}: Although the solution of this problem could be easily proved
directly, we will reduce this problem to the form (\ref{ae'}), (\ref{be'}) (%
\ref{ce'}) by setting%
\begin{eqnarray*}
\partial _{t}^{\ast }u &=&w \\
\partial _{t}^{\ast }w &=&D^{2}u-\left\vert u\right\vert ^{p-2}u
\end{eqnarray*}%
so problem (\ref{20}),(\ref{21}),(\ref{23}) can be reformulated as follows: 
\begin{equation*}
\mathbf{u}\in C^{1}(I^{\ast },\left( V^{\circ }\right) ^{2}):
\end{equation*}%
\begin{equation*}
A%
{{}^\circ}%
(D)\left[ \mathbf{u}\right] =0
\end{equation*}%
\begin{equation*}
u(0,x)=u_{0}(x);\ w(0,x)=u_{1}(x)
\end{equation*}%
where%
\begin{equation*}
\mathbf{u}=\left[ 
\begin{array}{c}
u \\ 
w%
\end{array}%
\right] \ \ \ \ \text{and}\ \ \ \ A%
{{}^\circ}%
(D)\left[ \mathbf{u}\right] =\left[ 
\begin{array}{c}
w \\ 
D^{2}u-\left\vert u\right\vert ^{p-2}u%
\end{array}%
\right] ;
\end{equation*}

Then the existence of a unique solution is guaranteed by Th. \ref{evol}
provided that we get an \textit{a priori} estimate. This is given by the
conservation of the energy: 
\begin{eqnarray*}
\partial _{t}^{\ast }E(\mathbf{u}(t,\cdot )) &=&\partial _{t}^{\ast }\doint
\left( \frac{1}{2}\left\vert w\right\vert ^{2}+\frac{1}{2}\left\vert
Du\right\vert ^{2}+\frac{1}{p}\left\vert u\right\vert ^{p}\right) dx \\
&=&\doint \left( w\partial _{t}^{\ast }w+Du\cdot \partial _{t}^{\ast
}Du+\left\vert u\right\vert ^{p-1}\partial _{t}^{\ast }u\right) dx \\
&=&\doint \left( w\partial _{t}^{\ast }w+Du\cdot Dw+\left\vert u\right\vert
^{p-1}w\right) dx \\
&=&\doint \left( \partial _{t}^{\ast }w-D^{2}u+\left\vert u\right\vert
^{p-1}\right) w\ dx=0
\end{eqnarray*}

Notice that $\partial _{t}^{\ast }$ and $D$ commute, since $\partial
_{t}^{\ast }$ "behaves" as an ordinary derivative and $D$ "behaves" as a
matrix in a finite dimensional space with the coefficients independent of
time.

$\square $

\bigskip

\subsection{Linear problems}

Let us consider the following linear boundary value problem: 
\begin{equation}
u\in C^{2}\left( \Omega \right) \cap C_{0}^{0}\left( \overline{\Omega }%
\right)
\end{equation}%
\begin{equation}
-\nabla \cdot \left[ k(x)\nabla u\right] +\lambda u=f(x)\ \ in\ \ \ \Omega
,\ f\in C^{0}\left( \overline{\Omega }\right) .
\end{equation}

For what we have discussed in section \ref{soe}, it is convenient to
translate it in the framework of ultrafunctions as follows:%
\begin{equation}
u\in V_{\text{\textsc{dir}}}^{\circ }\left( \Omega 
{{}^\circ}%
\right)  \label{p1}
\end{equation}%
\begin{equation}
-D\cdot \left[ k^{\ast }(x)Du\right] +\lambda u=f^{\ast }(x),  \label{p3}
\end{equation}

Since the operator $-D\cdot \left[ k^{\ast }(x)Du\right] $ is symmetric it
has a hyperfinite spectrum $\Sigma =\left\{ \lambda _{k}\right\} _{k\in K}$
with an orthonormal basis of eigenvalues $\left\{ e_{k}\right\} _{k\in K}.$
Then the Fredholm alternative holds and if $\lambda \notin \Sigma ,$ problem
(\ref{p1}), (\ref{p3}) has a unique solution given by%
\begin{equation*}
u(x)=\dsum\limits_{k\in K}\frac{f_{k}}{\lambda _{k}+\lambda }e_{k}\ \ \ 
\text{where}\ \ \ \ f(x)=\dsum\limits_{k\in K}f_{k}e_{k}(x).
\end{equation*}

The operator $L%
{{}^\circ}%
u:=-D\cdot \left[ k^{\ast }(x)Du\right] $ can be regarded as a sort of
selfadjoint realization of $L$ with respect to the scalar product $%
(u,v)\mapsto \doint uv~d_{\Omega }x.$

Let us examine the spectrum of $L%
{{}^\circ}%
$ in some cases. If $k(x)$ is a strictly positive smooth function and $%
\Omega $ is bounded, then the classical operator 
\begin{equation*}
Lu=-\nabla \cdot \left[ k(x)\nabla u\right]
\end{equation*}%
has a discrete spectrum of positive eigenvalues which correspond to an
orthonormal basis $\left\{ e_{k}\right\} _{k\in \mathbb{N}}$ of smooth
functions. Then the spectrum of $L%
{{}^\circ}%
u$ has an orthonormal basis $\left\{ h_{k}\right\} _{k\leq \dim \left( V_{%
\text{\textsc{dir}}}^{\circ }\left( \Omega 
{{}^\circ}%
\right) \right) }.$ For some infinite number $\bar{k}<\dim (V%
{{}^\circ}%
\left( \Omega 
{{}^\circ}%
\right) ),$ $\left\{ h_{k}\right\} $ coincides with the spectrum of $L^{\ast
}$; in particular, if $k\in \mathbb{N}$, the eigenvalues of $L%
{{}^\circ}%
$ coincide with the eigenvalues of $L$ and $h_{k}=e_{k}^{\circ }$.

If $k(x)$ is negative in a subset of $\Omega $ of positive measure, then $Lu$
has a continuum unbounded spectrum. In this case the eigenvalues of $L%
{{}^\circ}%
$ are infinitely close to each other.

\bigskip

\textbf{Example}: Let us consider the following ill posed problem relative
to the Tricomi equation:%
\begin{equation*}
\partial _{1}^{2}u+x_{1}\partial _{2}^{2}u=0\ in\ \Omega ;\ \ 0\in \Omega .
\end{equation*}%
\begin{equation*}
u=g(x)\ \ for\ x\in \partial \Omega
\end{equation*}%
In this case 
\begin{equation*}
L%
{{}^\circ}%
u:=-\left( D_{1}^{2}u+x_{1}D_{2}^{2}u\right) \theta _{\Omega }^{\circ }(x).
\end{equation*}%
It is not difficult to prove that for a "generic" open set $\Omega ,$ $0$ is
not in the spectrum of $L%
{{}^\circ}%
u.$ In this case, using Remark \ref{RRR}, this problem has a unique
ultrafunction solutions.

\section{A model for ultrafunctions\label{CSU}}

This section is devoted to prove Th. \ref{main} namely to the construction a
space of fine ultrafunctions. The construction presented here is the
simplest which we have been able to find. Nevertheless it is quite involved
but we do not know if a substantially simpler one exists. The main
difficulty relies in the fact that \textit{all} the properties of Def. \ref%
{A} need to be satisfied simultaneously; in particular, properties (\ref%
{lana}) and (\ref{loc}) are quite difficult to be obtained simultaneously

Our model of fine ultrafunctions combines the theory of ultrafunctions with
the techniques related to step functions. Roughly speaking we can say that,
in this model, the fine ultrafunctions, as well as the standard continuous
functions, can be well approximated by step functions and this fact is a
cornerstone of our construction.

\subsection{A construction of the Euclidean numbers\label{cen}}

As we have seen, a basic tool in the theory of ultrafunctions is the field
of Euclidean numbers and the notion of $\Lambda $-limit. Althout a
construction of a field which satisfies Axiom \ref{EU} can be found in
several papers, we repeat this construction here for the sake of the reader
(see e.g. \cite{ultra,belu2013}, see also \cite{BL2021,ben2022} for richer
models).

\begin{theorem}
Let $\mathfrak{L}$ be defined by (\ref{elle}). There exists a field $\mathbb{%
E\supset R}$ ($\mathbb{E\neq R}$) and a surjective field homomorphism, 
\begin{equation*}
J:\mathfrak{F}\left( \mathfrak{L},\mathbb{R}\right) \rightarrow \mathbb{E}
\end{equation*}%
namely a map with the following properties:%
\begin{equation*}
J\left( \varphi +\psi \right) =J\left( \varphi \right) +J\left( \psi \right)
,
\end{equation*}%
\begin{equation*}
J\left( \varphi \cdot \psi \right) =J\left( \varphi \right) \cdot J\left(
\psi \right) .
\end{equation*}
\end{theorem}

\textbf{Proof}: Let $\mathcal{U}$ be a fine ultrafilter on $\mathfrak{L}$,
namely a filter of sets such that

\begin{itemize}
\item {Maximality:} $Q\in \mathcal{U}\Leftrightarrow \mathfrak{L}\backslash
Q\notin \mathcal{U}$;

\item {Finess:} $\forall \lambda \in \mathfrak{L},\ Q\left[ \lambda \right]
\in \mathcal{U}$, where 
\begin{equation}
Q\left[ \lambda \right] :=\left\{ \mu \in \mathfrak{L}\ |\ \mu \supseteq
\lambda \right\} .  \label{finess}
\end{equation}
\end{itemize}

The existence of $\mathcal{U}$ is a well known and easy consequence of
Zorn's Lemma. We use $\mathcal{U}$ to introduce an equivalence relation on
nets, by letting for all $\psi ,\varphi \in \mathfrak{F}\left( \mathfrak{L},%
\mathbb{R}\right) $ 
\begin{equation*}
\varphi \approx _{\mathcal{U}}\psi \Longleftrightarrow \exists Q\in \mathcal{%
U},\ \forall \lambda \in Q,\ \varphi \left( \lambda \right) =\psi \left(
\lambda \right) .
\end{equation*}

We set 
\begin{equation*}
\mathbb{E}:=\mathfrak{F}\left( \mathfrak{L},\mathbb{R}\right) /\approx _{%
\mathcal{U}}
\end{equation*}%
and we denote by $\left[ \varphi \right] _{\mathcal{U}}$ the equivalence
classes. The operations on $\mathbb{E}$ can be easily defined by letting 
\begin{equation*}
\left[ \varphi \right] _{\mathcal{U}}+\left[ \psi \right] _{\mathcal{U}}=%
\left[ \varphi +\psi \right] _{\mathcal{U}};\ \ \ \left[ \varphi \right] _{%
\mathcal{U}}\cdot \left[ \psi \right] _{\mathcal{U}}=\left[ \varphi \cdot
\psi \right] _{\mathcal{U}}.
\end{equation*}

It is very well known (see e.g. \cite{keisler76}) and simple to show that,
thanks to $\mathcal{U}$ being an ultrafilter, $\mathbb{E}$ endowed with the
above operations is a field. The operator $J$ is defined by the canonical
projection 
\begin{equation*}
J\left( \varphi \right) :=\left[ \varphi \right] _{\mathcal{U}}.
\end{equation*}

$\square $

\begin{proposition}
\label{prop}The set 
\begin{equation*}
\mathbb{E}^{+}=\left\{ J(\varphi )\ |\ \forall \lambda \in \mathfrak{L},\
\varphi (\lambda )>0\right\}
\end{equation*}%
provides $\mathbb{E}$ of the linear order structure.
\end{proposition}

\textbf{Proof}: We need to prove that 
\begin{equation*}
\mathbb{E}=\mathbb{E}^{+}\cup \left\{ 0\right\} \cup \mathbb{E}^{-}
\end{equation*}%
namely, if we take $\xi \in \mathbb{E}\backslash \left\{ 0\right\} $, then $%
\exists \varphi \in \mathbb{E}^{+}$ such that 
\begin{equation*}
\xi =J(\varphi )\ \ or\ \ \xi =J(-\varphi ).
\end{equation*}%
By the surjectivity of $J$ there exists $\psi $ such that 
\begin{equation*}
\xi =J(\psi ).
\end{equation*}%
If we set 
\begin{eqnarray*}
R^{+} &=&\left\{ \lambda \in \mathfrak{L}\ |\ \psi (\lambda )>0\right\} \\
R^{-} &=&\left\{ \lambda \in \mathfrak{L}\ |\ \psi (\lambda )<0\right\} \\
R^{0} &=&\left\{ \lambda \in \mathfrak{L}\ |\ \psi (\lambda )=0\right\}
\end{eqnarray*}%
then 
\begin{equation*}
\mathfrak{L}=R^{+}\cup R^{-}\cup R^{0}
\end{equation*}%
Since $\mathcal{U}$ is an ultrafilter, only one of the sets $%
R^{+},R^{-},R^{0}$ is in $\mathcal{U}$. It is not possible that $R^{0}\in 
\mathcal{U}$ since 
\begin{equation*}
0\neq \xi =J\left( \psi \right) =\left[ \psi \right] _{\mathcal{U}}.
\end{equation*}%
If $R^{+}\in \mathcal{U}$, then $\xi \in \mathbb{E}^{+}$ since 
\begin{equation*}
\xi =J\left( \psi \right) =J\left( \psi ^{+}\right)
\end{equation*}%
where%
\begin{equation*}
\psi ^{+}(\lambda ):=\left\{ 
\begin{array}{cc}
\psi (\lambda ) & if\ \ \lambda \in R^{+} \\ 
&  \\ 
1 & if\ \ \lambda \in R^{-}\cup R^{0}%
\end{array}%
\right.
\end{equation*}

If $R^{-}\in \mathcal{U}$, then, arguing in a similar way, $\xi \in \mathbb{E%
}^{-}.$

$\square $

\bigskip

Now, we can define the notion of $\Lambda $-limit. For every net $\varphi :%
\mathfrak{L}\rightarrow V_{0}(\mathbb{R}),$ we set 
\begin{equation}
\lim_{\lambda \uparrow \Lambda }\varphi (\lambda )=J\left( \varphi \right)
\label{liz}
\end{equation}%
and for a net $\varphi :\mathfrak{L}\rightarrow V_{n}(\mathbb{R}),$ $n\geq
0, $ we define the $\Lambda $-limit by induction. If $n=0,$ $\lim_{\lambda
\uparrow \Lambda }$ $\varphi \left( \lambda \right) $ has been defined
above; if $n>0,$ we set%
\begin{equation}
\lim_{\lambda \uparrow \Lambda }\ \varphi (\lambda ):=\left\{ \lim_{\lambda
\uparrow \Lambda }\ \psi (\lambda )\ |\ \forall \lambda \in \mathfrak{L},\
\psi (\lambda )\in \varphi (\lambda )\right\} .  \label{vivaldi}
\end{equation}

\begin{theorem}
The $\Lambda $-limit defined by (\ref{liz}) and (\ref{vivaldi}) satisfies
the requests of Axiom \ref{EU}.
\end{theorem}

\textbf{Proof}: Let us check each request separately.

\ref{EU}.\ref{EU1} - If eventually\textit{\ }%
\begin{equation*}
\varphi (\lambda )=\psi (\lambda ),
\end{equation*}%
the above relation is satisfied by every $\lambda \in Q\left[ \lambda _{0}%
\right] $ for a suitable $\lambda _{0}.$ By (\ref{finess}), $Q\left[ \lambda
_{0}\right] \in \mathcal{U}$ and hence $J\left( \varphi \right) =J\left(
\psi \right) .$

\ref{EU}.\ref{A3} - By (\ref{vivaldi}), it is immediate to see that 
\begin{equation*}
\lim_{\lambda \uparrow \Lambda }\left\{ \varphi _{1}(\lambda ),...,\varphi
_{n}(\lambda )\right\} \subseteq \left\{ \lim_{\lambda \uparrow \Lambda
}\varphi _{1}(\lambda ),...,\lim_{\lambda \uparrow \Lambda }\varphi
_{n}(\lambda )\right\}
\end{equation*}%
It is more delicate to show that 
\begin{equation*}
\left\{ \lim_{\lambda \uparrow \Lambda }\varphi _{1}(\lambda
),...,\lim_{\lambda \uparrow \Lambda }\varphi _{n}(\lambda )\right\}
\subseteq \lim_{\lambda \uparrow \Lambda }\left\{ \varphi _{1}(\lambda
),...,\varphi _{n}(\lambda )\right\}
\end{equation*}%
namely that for every $\xi \in \lim_{\lambda \uparrow \Lambda }\left\{
\varphi _{1}(\lambda ),...,\varphi _{n}(\lambda )\right\} ,$ $\exists \bar{k}
$ such that%
\begin{equation*}
\xi =\lim_{\lambda \uparrow \Lambda }\varphi _{\bar{k}}(\lambda )
\end{equation*}%
Take 
\begin{equation*}
\xi \in \lim_{\lambda \uparrow \Lambda }\left\{ \varphi _{1}(\lambda
),...,\varphi _{n}(\lambda )\right\}
\end{equation*}%
then%
\begin{equation*}
\xi =\lim_{\lambda \uparrow \Lambda }\psi (\lambda )
\end{equation*}%
where $\forall \lambda \in \mathfrak{L},$ $\exists k\leq n,$ such that 
\begin{equation*}
\psi (\lambda )=\varphi _{k}(\lambda )
\end{equation*}%
Now, if we set, for $k\leq n,$%
\begin{equation*}
R_{k}=\left\{ \lambda \in \mathfrak{L}\ |\ \psi (\lambda )=\varphi
_{k}(\lambda )\right\}
\end{equation*}%
it turns out that 
\begin{equation*}
R_{1}\cup R_{2}\cup ....\cup R_{n}=\mathfrak{L.}
\end{equation*}%
Since $\mathcal{U}$ is a ultrafilter, then one and only one of the $R_{k}$'s
is in $\mathcal{U}$. If $R_{\bar{k}}\in \mathcal{U}$, then%
\begin{equation*}
\xi =\lim_{\lambda \uparrow \Lambda }\psi (\lambda )=J\left( \psi \right)
=J\left( \varphi _{\bar{k}}\right) =\lim_{\lambda \uparrow \Lambda }\varphi
_{\bar{k}}(\lambda )
\end{equation*}

\ref{EU}.\ref{A4} - It is nothing else but the definition (\ref{vivaldi}).

\ref{EU}.\ref{EU2} - It follows directly by the definition of $\mathbb{E}$
and Prop. \ref{prop}.

$\square $

\subsection{Hyperfinite step functions\label{SF}}

The step functions are easy to handle ad hence they are largely used many
branches of mathematics and in particular in nonstandard analysis. Let us
describe the kind of step functions that will be considered in this paper.

Given an infinite hypercube 
\begin{equation*}
Q=\left[ -L,L\right) _{\mathbb{E}}^{N}\subset \mathbb{E}^{N},\ \ \left[
-L,L\right) _{\mathbb{E}}=\left\{ \xi \in \mathbb{E}\ |\ -L\leq \xi
<L\right\}
\end{equation*}%
we say that a partition $\left\{ S_{a}\right\} _{a\in \Gamma }$ of $Q$ is 
\textbf{fine} if

\begin{itemize}
\item $\Gamma $ is an hyperfinite set,

\item the index $a$ is a point in $int(S_{a})$ and if $S_{a}\cap \mathbb{R}%
^{N}=\left\{ r\right\} ,$ then $a=r.$

\item $\forall a\in \Gamma ,$ $\bar{\chi}_{S_{a}}\in V^{\ast }.$

\item the \textbf{size }of the partitition, namely the number%
\begin{equation*}
\eta =\ \underset{a\in \Gamma }{\max }\ diam\left( S_{a}\right) ,
\end{equation*}%
is infinitesimal; here $diam\left( S_{a}\right) $ denotes the diameter of $%
S_{a}.$
\end{itemize}

\begin{definition}
Given a partition $\left\{ S_{a}\right\} _{a\in \Gamma },$ a \textbf{step
function} $u:\mathbb{E}^{N}\rightarrow \mathbb{E}$ is defined as follows:%
\begin{equation*}
u(x)=\sum_{a\in \Gamma }u_{a}\chi _{S_{a}}(x).
\end{equation*}%
where $u_{(\cdot )}:\Gamma \rightarrow \mathbb{E}$ is an internal function
namely a grid function.
\end{definition}

The integral of a step function $u$ is given by%
\begin{equation*}
\int^{\ast }u(x)dx=\sum_{a\in \Gamma }u(a)m(S_{a});
\end{equation*}

Given a fine partition $\left\{ S_{a}\right\} _{a\in \Gamma },$ any internal
function $u\in \mathfrak{F}(\mathbb{E}^{N},\mathbb{E})$ can be approximated
by a step function $\breve{u}$ by setting%
\begin{equation}
\breve{u}(x)=\sum_{a\in \Gamma }u(a)\chi _{S_{a}}(x).  \label{step}
\end{equation}

The next theorem shows a very interesting property of the infinitesimal
partition in the contest of the Euclidean numbers.

\begin{theorem}
\label{super}Let $\left\{ S_{a}\right\} _{a\in \Gamma }$ be a fine
partition: then, if $f\in \mathcal{L}^{1},$%
\begin{equation*}
\int f(x)dx\sim \int^{\ast }\breve{f}(x)dx
\end{equation*}
\end{theorem}

First we need the following Lemma.

\begin{lemma}
\label{le}Let $\left\{ S_{a}\right\} _{a\in \Gamma }$ be as in Th. \ref%
{super}; if $f$ is a measurable bounded function with compact support, then 
\begin{equation*}
\int f(x)dx\sim \int^{\ast }\breve{f}(x)dx
\end{equation*}
\end{lemma}

\textbf{Proof: }Since $\left\{ S_{a}\right\} _{a\in \Gamma }$ is an internal
set, then 
\begin{equation*}
\left\{ S_{a}\right\} _{a\in \Gamma }=\lim_{\lambda \uparrow \Lambda }\
\left\{ S_{a}^{\lambda }\right\} _{a\in \Gamma _{\lambda }}
\end{equation*}%
where $\left\{ S_{a}^{\lambda }\right\} _{a\in \Gamma _{\lambda }}$ is a net
of finite measurable partitions. Then 
\begin{equation*}
\breve{f}=\lim_{\lambda \uparrow \Lambda }f_{\lambda }
\end{equation*}%
where%
\begin{equation*}
f_{\lambda }(x)=\sum_{a\in \Gamma _{\lambda }}f(a)\ m(S_{a}^{\lambda }),
\end{equation*}%
are standard measurable step functions. Since $f(a)$ is bounded and with
compact support, the $f_{\lambda }$'s are uniformly bounded and with compact
support. Moreover, since $\mathbb{R}^{N}\subset \Gamma =\lim_{\lambda
\uparrow \Lambda }\Gamma _{\lambda },$ if you take $x\in \mathbb{R}^{N},$ we
have that, for $\lambda $ sufficiently large, $x\in \Gamma _{\lambda }$ and
hence 
\begin{equation*}
f_{\lambda }(x)=f(x)
\end{equation*}%
namely the $f_{\lambda }$'s converge pointwise to $f$: 
\begin{equation*}
\lim_{\lambda \rightarrow \Lambda }f_{\lambda }(x)=f(x)
\end{equation*}%
Then by the Lebesgue Dominated Convergence Theorem 
\begin{equation*}
\int \lim_{\lambda \rightarrow \Lambda }f_{\lambda }(x)dx=\int f(x)dx.
\end{equation*}%
On the other hand, by (\ref{bn}), 
\begin{eqnarray*}
\lim_{\lambda \rightarrow \Lambda }\int f_{\lambda }(x)dx &=&st\left(
\lim_{\lambda \uparrow \Lambda }\int f_{\lambda }(x)dx\right) \\
&=&st\left( \int^{\ast }\lim_{\lambda \uparrow \Lambda }f_{\lambda
}(x)dx\right) =st\left( \int^{\ast }\breve{f}(x)dx\right)
\end{eqnarray*}

$\square $

\textbf{Proof of Th. \ref{super}}. Since $\mathcal{L}_{c}^{\infty }$ is
dense in $\mathcal{L}^{1}$ there is a sequence $f_{n}$ of measurable bounded
functions with compact support such that 
\begin{equation*}
f(x)=\sum_{n=0}^{\infty }f_{n}(x);\ \int f(x)\ dx=\sum_{n=0}^{\infty }\int
f_{n}(x)\ dx
\end{equation*}%
Then, by transfer%
\begin{equation*}
f^{\ast }(x)=\sum_{n=0}^{\infty ^{\ast }}f_{n}^{\ast }(x);\ \int^{\ast
}f^{\ast }(x)\ dx=\sum_{n=0}^{\infty ^{\ast }}\int^{\ast }f_{n}^{\ast }(x)\
dx
\end{equation*}%
Since $\left\{ S_{a}\right\} _{a\in \Gamma }$ has been fixed we have that%
\begin{equation*}
\breve{f}(x)=\sum_{n=0}^{\infty ^{\ast }}\breve{f}_{n}(x);\ \int^{\ast }%
\breve{f}(x)\ dx=\sum_{n=0}^{\infty ^{\ast }}\int^{\ast }\breve{f}_{n}\ dx
\end{equation*}%
Hence, for every $N\in \mathbb{N},$%
\begin{eqnarray}
\left\vert \int^{\ast }f^{\ast }(x)\ dx-\int^{\ast }\breve{f}(x)\
dx\right\vert &\leq &\left\vert \sum_{n=0}^{\infty ^{\ast }}\int^{\ast
}f_{n}^{\ast }(x)\ dx-\sum_{n=0}^{\infty ^{\ast }}\int^{\ast }\breve{f}_{n}\
dx\right\vert  \notag \\
&\leq &\sum_{n=0}^{\infty ^{\ast }}\left\vert \int^{\ast }f_{n}^{\ast }(x)\
dx-\int^{\ast }\breve{f}_{n}\ dx\right\vert  \notag \\
&=&\sum_{n=0}^{N}\left\vert \int^{\ast }f_{n}^{\ast }(x)\ dx-\int^{\ast }%
\breve{f}_{n}\ dx\right\vert  \label{pal} \\
&&+\sum_{n=N+1}^{\infty ^{\ast }}\left\vert \int^{\ast }f_{n}^{\ast }(x)\
dx-\int^{\ast }\breve{f}_{n}\ dx\right\vert  \notag \\
&\leq &\sum_{n=0}^{N}\left\vert \int^{\ast }f_{n}^{\ast }(x)\ dx-\int^{\ast }%
\breve{f}_{n}\ dx\right\vert  \notag \\
&&+\sum_{n=N+1}^{\infty }\left\vert \int^{\ast }f_{n}^{\ast }(x)\
dx\right\vert +\sum_{n=N+1}^{\infty }\left\vert \int^{\ast }\breve{f}_{n}\
dx\right\vert  \notag
\end{eqnarray}

Now choose $\varepsilon \in \mathbb{R}^{+}$ arbitrarily and $N\in \mathbb{N}$
such that 
\begin{equation*}
\sum_{n=N+1}^{\infty ^{\ast }}\left\vert \int f_{n}(x)\ dx\right\vert
=\sum_{n=N+1}^{\infty ^{\ast }}\left\vert \int^{\ast }f^{\ast }(x)\
dx\right\vert <\varepsilon
\end{equation*}%
consequently, 
\begin{equation*}
\sum_{n=N+1}^{\infty }\left\vert \int^{\ast }\breve{f}_{n}\ dx\right\vert
<2\varepsilon
\end{equation*}

Since $N$ is a finite number, by lemma \ref{le}, we have that 
\begin{equation*}
\sum_{n=0}^{N}\left\vert \int^{\ast }f_{n}^{\ast }(x)\ dx-\int^{\ast }\breve{%
f}_{n}\ dx\right\vert \sim 0
\end{equation*}%
Then, by (\ref{pal}), we have that 
\begin{equation*}
\left\vert \int^{\ast }f^{\ast }(x)\ dx-\int^{\ast }\breve{f}(x)\
dx\right\vert \leq 4\varepsilon .
\end{equation*}%
The conclusion follows from the arbitrariness of $\varepsilon $.

$\square $

\subsection{$\protect\sigma $-bases}

\begin{definition}
Let $W\subset \mathfrak{F}\left( \mathbb{R}^{N}\right) $ be a function space
of finite dimension; we say that a family of functions $\left\{ \sigma
_{a}\right\} _{a\in \mathfrak{K}},$ $\mathfrak{K}\subset \mathbb{R}^{N}$ is
a $\sigma $-basis for $W$ if every function $u$ in $W$ can be written as
follows: 
\begin{equation*}
u(x)=\sum_{a\in \mathfrak{K}}u(a)\sigma _{a}(x).
\end{equation*}%
in a unique way.
\end{definition}

Given $W\subset \mathfrak{F}\left( \mathbb{R}^{N}\right) $ and a set of
points $\mathfrak{K=}\left\{ a_{1},...,a_{k}\right\} \subset \mathbb{R}^{N},$
we can define "restriction" map%
\begin{equation}
\Phi :W\rightarrow \mathfrak{F}\left( \mathfrak{K}\right)  \label{map}
\end{equation}%
\begin{equation*}
\Phi (f)=\left( f\left( a_{1}\right) ,...,f\left( a_{k}\right) \right)
\end{equation*}%
If $\left\{ e_{1},,,,e_{m}\right\} $ is a basis of $W,$ then $\Phi $ can be
"represented" by the matrix%
\begin{equation}
\left\{ e_{n}\left( a_{l}\right) \right\} _{n\leq m,l\leq k}  \label{M}
\end{equation}

\begin{definition}
\label{icr}Let $W\subset \mathfrak{F}\left( \mathbb{R}^{N}\right) $ be a
function space of finite dimension; we say that a set $\mathfrak{K}=\left\{
a_{1},...,a_{k}\right\} \subset \mathbb{R}^{N}$ is:

\begin{itemize}
\item \textbf{independent} in $W$ if the map $\Phi $ is surjective and hence
for any $k$-ple of points $\left( c_{1},...,c_{k}\right) \in \mathbb{R}^{N}$%
, there exists $f\in W$ such that%
\begin{equation}
f\left( a_{l}\right) =c_{l};\ \ l=1,....,k.  \label{elina}
\end{equation}%
in this case the matrix (\ref{M}) has rank $k.$

\item \textbf{complete }in $W$ if the map $\Phi $ is bijective and hence
there exists a \textbf{unique} $f\in W$ which satisfies (\ref{elina}); in
this case 
\begin{equation}
\det \left[ e_{k}\left( a_{l}\right) \right] \neq 0  \label{det}
\end{equation}

\item \textbf{redundant }in $W$ if the map $\Phi $ is injective and hence 
\begin{equation*}
\left( \forall f\in W,\forall a\in \mathfrak{K},\ f\left( a\right) =0\right)
\Rightarrow \left( f=0\right) ;
\end{equation*}%
in this case the matrix (\ref{M}) has rank $m.$
\end{itemize}
\end{definition}

Notice that a set of points is complete in $W$ if and only if it is
independent and redundant.

\begin{lemma}
\label{a}Let $W\subset \mathfrak{F}\left( \mathbb{R}^{N}\right) $ be a
function space of finite dimension and let $\mathfrak{K}=\left\{
a_{1},...,a_{m}\right\} \subset \mathbb{R}^{N}$. Then it is complete if and
only if there exists a $\sigma $-basis $\left\{ \sigma _{a}(x)\right\}
_{a\in \mathfrak{K}}$.
\end{lemma}

\textbf{Proof:} Let $\mathfrak{K}=\left\{ a_{1},...,a_{n}\right\} \subset 
\mathbb{R}^{N}$ be a of complete set points and let $\left\{ \zeta
_{a}\right\} _{a\in \mathfrak{K}}$ be the "canonical" basis in $\mathfrak{F}%
\left( \mathfrak{K}\right) $ namely%
\begin{equation*}
\zeta _{a}(b)=\delta _{ab}
\end{equation*}%
Then, $\left\{ \Phi ^{-1}\left( \zeta _{a}\right) \right\} _{a\in \mathfrak{K%
}}$ is a $\sigma $-basis. If $\left\{ \sigma _{a}(x)\right\} _{a\in 
\mathfrak{K}}$ is a $\sigma $-basis, then the map defined by $\Phi \left(
\sigma _{a}\right) =\zeta _{a}$ is bijective.

$\square $

\bigskip

It is evident that every finite dimensional vector space has a redundant set
of points. The existence of a complete set is a consequence of the following
theorem.

\begin{theorem}
\label{puta4}Let $\mathfrak{R}$ be redundant in $W$. Then there exists a set 
$\mathfrak{C}\subseteq \mathfrak{R},$ complete in $W$.
\end{theorem}

\textbf{Proof}: The matrix (\ref{M}) has rank $m=\dim \left( W\right) $. Let 
$M_{m}$ be a nonsingular $m\times m$ submatrix of the matrix $M$ and let $%
\mathfrak{C}$ be the image of the operator defined by $M_{m}.$ So, the set $%
\left\{ \sigma _{a}(x)\right\} _{a\in \mathfrak{C}}$ is a $\sigma $-basis
and by lemma \ref{a}, $\mathfrak{C}$ is complete.

$\square $

\begin{corollary}
\label{beppa}Every finite or hyperfinite dimensional vector space $W$ has a
complete system of points.
\end{corollary}

\textbf{Proof}: If $W$ is finite dimensional, take any redundant set of
points in $W$ and apply Th. \ref{puta4}. If the dimension of $W$ is
hyperfinite, take the $\Lambda $-limit.

$\square $

\begin{corollary}
\label{beppa1}Every finite or hyperfinite dimensional vector space $W$ has a 
$\sigma $-basis.
\end{corollary}

\textbf{Proof}: It is an immediate consequence of Corollary \ref{beppa} and
Lemma \ref{a}.

$\square $

\bigskip

We end this section with a technical lemma which will be used in the proof
of Lemma \ref{pluto}.

\begin{lemma}
\label{ben}Let $W$ be a finite or hyperfinite dimensional vector space and
let $r\in \mathbb{R}^{N}$ be a point such that 
\begin{equation*}
\exists w\in W,\ w(r)\neq 0.
\end{equation*}%
Then $W$ has a complete set of points $\mathfrak{K}$ with $r\in \mathfrak{K.}
$
\end{lemma}

\textbf{Proof}: Let $\left\{ a_{1},...,a_{m}\right\} $ be complete set of
points and let $\left\{ e_{n}\left( a_{l}\right) \right\} _{n,l\leq m}$ be
the relative matrix defined by (\ref{M}). Since $w(r)\neq 0,\ $the vector $%
\left\{ e_{n}\left( r\right) \right\} _{n\leq m}$ is $\neq 0$. Then there
exists a point $a_{\bar{l}}$ such that the vector $\left\{ e_{n}\left( a_{%
\bar{l}}\right) \right\} _{n\leq m}$ is a linear combination of the vectors $%
\left\{ e_{n}\left( a_{l}\right) \right\} _{l\neq \bar{l}}\cup \left\{
e_{n}\left( r\right) \right\} _{n\leq m},$ and hence the set%
\begin{equation*}
\mathfrak{K}:=\left\{ a_{l}\ |\ l\neq \bar{l}\right\} \cup \left\{ r\right\}
\end{equation*}%
is complete.

$\square $

\subsection{The space $V_{\Lambda }$}

This and the next sections are devoted to construct the space $V_{\Lambda }$
so that it is possible to define a pointwise integral and a generalized
derivative which satisfy the requests of Def. \ref{A}. This will be done in
three steps building three spaces $W_{\Lambda }^{0},W_{\Lambda },V_{\Lambda
} $ with 
\begin{equation}
V^{\circledcirc }\subseteq W_{\Lambda }^{0}\subseteq W_{\Lambda }\subseteq
V_{\Lambda }.  \label{inc}
\end{equation}

We fix 
\begin{equation}
Q=\left[ L,L\right) ^{N}  \label{Q}
\end{equation}%
where $L\in \mathbb{E}$ is an infinite number such that for every $%
u_{\lambda }\in V\cap \lambda ,$ we have that 
\begin{equation}
supp\left( u\right) \subset \left( -L,L\right) ^{N}.  \label{sid}
\end{equation}%
We set 
\begin{equation}
W_{\Lambda }^{0}:=span\left\{ \overline{uv}\ |\ u,v\in V^{\circledcirc }\cup
\left\{ \bar{\chi}_{Q}\right\} \right\}  \label{F3}
\end{equation}%
We recall that the operator $u\mapsto \bar{u}$ has been defined by (\ref%
{due+}).

By Cor. \ref{beppa1}, the space $W_{\Lambda }^{0}$ has a $\sigma $-basis;
however we need to have a $\sigma $-basis which satisfies suitable
properties. For this reason we need the following lemma.

\begin{lemma}
\label{pluto}There exists a hyperfinite dimensional vector space $W_{\Lambda
}\supseteq W_{\Lambda }^{0}$ which admits a $\sigma $-basis $\left\{ \tau
_{a}\right\} _{a\in \Gamma _{W}}$ such that%
\begin{equation}
\mathbb{R}^{N}\subset \Gamma _{W};  \label{tt1}
\end{equation}%
and%
\begin{equation}
supp\left( \tau _{a}\right) \subset B_{2\varepsilon }\left( x_{j}\right)
\label{tt2}
\end{equation}%
where $B_{2\varepsilon }\left( x_{j}\right) \subset \mathbb{E}^{N}$ is a
ball of radius $2\varepsilon \sim 0.$
\end{lemma}

\textbf{Proof}: Let $\left\{ E_{k}\right\} _{k\in \mathfrak{K}}$ be a fine
partition of $Q$ such that, for every $k,\ diam\left( E_{k}\right)
<\varepsilon \sim 0.$ Then, $\forall u\in W_{\Lambda }^{0},$ we have that%
\begin{equation*}
u(x)=\sum_{k\in \mathfrak{K}}\overline{u(x)\chi _{E_{k}}\left( x\right) }
\end{equation*}%
Now, we set 
\begin{equation*}
W(E_{k})=\left\{ \overline{u(x)\chi _{E_{k}}\left( x\right) }|\ u\in
W_{\Lambda }^{0}\right\}
\end{equation*}%
and 
\begin{equation*}
W_{\Lambda }=W(E_{1})\oplus ....\oplus \ W(E_{r})
\end{equation*}%
Clearly $W_{\Lambda }^{0}\subset W_{\Lambda }.$ Now, for every $k\in 
\mathfrak{K},$ take a $\sigma $-basis $\left\{ \tau _{a}\right\} _{a\in
\Gamma _{k}}$ of $W(E_{k})$ which exists by Cor. \ref{beppa}. If $\mathbb{R}%
^{N}\cap E_{k}=\left\{ r_{k}\right\} \neq \varnothing ,$ we can take $\Gamma
_{k}$ such that $r_{k}\in \Gamma _{k};$ this is possible by lemma \ref{ben}.
Now, if we set 
\begin{equation*}
\Gamma _{W}=\bigcup_{k\in \mathfrak{K}}\Gamma _{k}
\end{equation*}%
we have that $\left\{ \tau _{a}\right\} _{a\in \Gamma _{W}}$ is a $\sigma $%
-basis which satisfies our requests.

$\square $

\bigskip

Now, we fix once for ever a positive infinitesimal%
\begin{equation}
\gamma <\left\vert \Gamma _{W}\right\vert ^{-1}  \label{gamma}
\end{equation}%
and we set 
\begin{equation*}
\Omega _{W}:=\bigcup_{a\in \Gamma _{W}}B_{\varepsilon }(a)
\end{equation*}%
where, for $a$ in $\Gamma _{W},$ $B_{\varepsilon }(a)$ is the ball of center 
$a$ and radius $\varepsilon .$

\begin{lemma}
\label{muti}We can take $\varepsilon $ so small that $\forall a,b\in \Gamma
_{W},$ $a\neq b$%
\begin{equation*}
B_{\varepsilon }(a)\cap B_{\varepsilon }(b)=\varnothing
\end{equation*}%
\begin{equation}
\left( 1-\gamma ^{2}\right) m\left( B_{\varepsilon }\right) \leq
\int_{B_{\varepsilon }(a)}^{\ast }\tau _{a}(x)dx\leq \left( 1-\gamma
^{2}\right) m\left( B_{\varepsilon }\right)  \label{mozart1}
\end{equation}%
\begin{equation}
\int_{\Omega _{W}\backslash B_{\varepsilon }(a)}^{\ast }\left\vert \tau
_{a}(x)\right\vert dx\leq \gamma ^{2}m\left( B_{\varepsilon }\right)
\label{mozart2}
\end{equation}%
where $m\left( B_{\varepsilon }\right) :=m(B_{\varepsilon }(0))$ denotes the
Lebesgue measure of $B_{\varepsilon }(0).$
\end{lemma}

\textbf{Proof}: Since the functions $\tau _{a}$ are epilogic the conclusion
follows from the fact that $\tau _{a}(a)=1$ and, for $b\neq a,$ $\tau
_{a}(b)=0.$

$\square $

\bigskip

Now, we set 
\begin{equation*}
\Omega _{Z}:=Q\backslash \Omega _{W}
\end{equation*}%
and we denote by $\left\{ S_{a}\right\} _{a\in \Gamma _{Z}^{\eta }}$ a fine
partition of $\Omega _{Z}$ of size $\eta \leq \varepsilon ;$ moreover we set%
\begin{equation*}
Z_{\eta }:=span\left\{ \bar{\chi}_{S_{a}}\ |\ a\in \Gamma _{Z}^{\eta
}\right\} .
\end{equation*}

\begin{lemma}
\label{ravenna}$W_{\Lambda }\cap Z_{\eta }=\left\{ 0\right\} $
\end{lemma}

\textbf{Proof}: Let $u\in W_{\Lambda }\cap Z_{\eta }.$ Since $u\in Z_{\eta
}, $ 
\begin{equation*}
\forall a\in \Gamma _{W},\ u\left( a\right) =0.
\end{equation*}%
$\Gamma _{W}$ is a complete set for $W_{\Lambda }$ (see Def. \ref{icr}) and $%
u\in W_{\Lambda }$, then $u=0.$

$\square $\bigskip

We set%
\begin{equation}
V_{\Lambda }^{\eta }:=W_{\Lambda }\oplus Z_{\eta };\ \Gamma ^{\eta }:=\Gamma
_{W}\cup \Gamma _{Z}^{\eta }  \label{mary}
\end{equation}%
Notice that $\left\{ \tau _{a}\right\} _{a\in \Gamma _{W}}\cup \left\{ \bar{%
\chi}_{S_{a}}\right\} _{a\in \Gamma _{Z}}$ is a basis of $V_{\Lambda }$ but
it is not a $\sigma $-basis even if$\ \left\{ \tau _{a}\right\} _{a\in
\Gamma _{W}}$ is a $\sigma $-basis of $W_{\Lambda }\ $and$\ \left\{ \bar{\chi%
}_{S_{a}}\right\} _{a\in \Gamma _{Z}^{\eta }}$ is a $\sigma $-basis of $%
Z_{\eta }^{1}.$ If we put 
\begin{equation*}
S_{a}=B_{\varepsilon }(a)
\end{equation*}%
we have that 
\begin{equation*}
\left\{ S_{a}\right\} _{a\in \Gamma ^{\eta }}=\left\{ B_{\varepsilon
}(a)\right\} _{a\in \Gamma _{W}}\cup \left\{ S_{a}\right\} _{a\in \Gamma
_{Z}^{\eta }};\ \Gamma =\Gamma _{W}\cup \Gamma _{Z}^{\eta }
\end{equation*}%
is a fine partition of $Q$ of size $2\varepsilon $.

Using the notation (\ref{step}), we set, 
\begin{equation*}
\breve{\tau}_{a}(x)=\sum_{s\in \Gamma ^{\eta }}\tau _{a}(s)\chi _{S_{s}}(x)
\end{equation*}

\begin{lemma}
\label{pap1}If $\eta $ is sufficiently small, $\forall a\in \Gamma _{W},\ $%
\begin{equation*}
\int_{\Omega _{Z}}^{\ast }\left\vert \tau _{a}(x)-\breve{\tau}%
_{a}(x)\right\vert \leq \gamma ^{2}m\left( B_{\varepsilon }\right) .
\end{equation*}
\end{lemma}

\textbf{Proof}: Since the $\tau _{a}$'s are a hyperfinite number of Riemann
integrable functions with compact support, then we can choose $\eta $ so
small that our request be satisfied.

$\square $

\begin{remark}
The proof of lemma \ref{pap1} is the only point in which it is required that
the functions in $V$ be Riemann-integrable.
\end{remark}

From now on, we will fix $\eta $ in such a way that Lemma \ref{pap1} be
satisfied and we will write $V_{\Lambda },$ $Z,\ \Gamma \ $and$\ \Gamma _{Z}$
instead of $V_{\Lambda }^{\eta },$ $Z_{\eta },\ \Gamma ^{\eta }\ $and$\
\Gamma _{Z}^{\eta }.$

\begin{lemma}
\label{rosa}We set%
\begin{equation}
\sigma _{a}(x)=\left\{ 
\begin{array}{cc}
\bar{\chi}_{S_{a}}(x) & if\ \ a\in \Gamma _{Z} \\ 
&  \\ 
\tau _{a}(x)-\overline{\breve{\tau}_{a}(x)\chi _{\Omega _{Z^{1}}}(x)} & if\
\ a\in \Gamma _{W}%
\end{array}%
\right. ,  \label{sig}
\end{equation}%
then $\left\{ \sigma _{a}\right\} _{a\in \Gamma }$ is a $\sigma $-basis of $%
V_{\Lambda }$ such that 
\begin{equation}
supp\left( \sigma _{a}\right) \subset B_{2\varepsilon }(a).  \label{lola}
\end{equation}
\end{lemma}

\textbf{Proof}: First of all let us check that 
\begin{equation}
\sigma _{a}\left( c\right) =\delta _{ac}  \label{bocelli}
\end{equation}%
for every $c\in \Gamma _{W}\cup \Gamma _{Z}.$ If $a\in \Gamma _{Z}$, then, $%
\sigma _{a}\left( c\right) =\bar{\chi}_{S_{a}}\left( c\right) =\delta _{ac}.$

Let us see the case in which $a\in \Gamma _{W}$. In this case, (\ref{sig})
takes the following form:%
\begin{equation*}
\sigma _{a}(x)=\tau _{a}(x)-\dsum\limits_{b\in \Gamma _{Z}}\tau _{a}(b)\bar{%
\chi}_{S_{b}}(x);
\end{equation*}

\begin{itemize}
\item if $c\in \Gamma _{Z}$ then 
\begin{equation*}
\sigma _{a}\left( c\right) =\tau _{a}\left( c\right) -\dsum\limits_{b\in
\Gamma _{Z}}\tau _{a}(b)\bar{\chi}_{S_{b}}(x)=\tau _{a}\left( c\right) -\tau
_{a}(c)=0
\end{equation*}

\item if $c\in \Gamma _{W}$, since $\left\{ \tau _{a}\right\} _{a\in \Gamma
_{W}}$ is a $\sigma $-basis, then%
\begin{equation*}
\sigma _{a}\left( c\right) =\tau _{a}\left( c\right) -\dsum\limits_{b\in
\Gamma _{Z}}\tau _{a}(b)\bar{\chi}_{S_{b}}(x)=\delta _{ac}-0=\delta _{ac}
\end{equation*}
\end{itemize}

So the $\sigma _{a}$'s are $\left\vert \Gamma \right\vert =\left\vert \Gamma
_{W}\right\vert +\left\vert \Gamma _{Z}\right\vert $ elements independent in 
$V_{\Lambda }$ and hence they span all $V_{\Lambda }.$ So, they form a $%
\sigma $-basis.

Now let us prove (\ref{lola}). If $a\in \Gamma _{Z},\ $%
\begin{equation*}
supp\left[ \sigma _{a}\right] =\overline{Q_{a}\backslash \Omega _{W}}\subset
B_{2\eta }(a)\subset B_{2\varepsilon }(a);
\end{equation*}%
if $a\in \Gamma _{Z}^{\eta },$ by (\ref{tt2}) 
\begin{equation*}
supp\left[ \sigma _{a}\right] \subseteq supp\left[ \tau _{a}\right] \subset
B_{2\varepsilon }(a).
\end{equation*}

$\square $

\begin{lemma}
\label{ganzone}For every $a\in \Gamma ^{\eta },$%
\begin{equation}
\left( 1-3\gamma ^{2}\right) m\left( S_{a}\right) \leq \int^{\ast }\sigma
_{a}dx\leq \left( 1+3\gamma ^{2}\right) m\left( S_{a}\right) ;  \label{moz3}
\end{equation}
\end{lemma}

\textbf{Proof}: If $a\in \Gamma _{Z},\ \sigma _{a}=\bar{\chi}_{S_{a}}$ and
hence%
\begin{equation*}
\int^{\ast }\sigma _{a}dx=\int^{\ast }\bar{\chi}_{S_{a}}dx=m\left(
S_{a}\right)
\end{equation*}

Now let us consider the case $a\in \Gamma _{W};$ we have that%
\begin{equation}
\tau _{a}(x)=\overline{\tau _{a}(x)\chi _{B_{\varepsilon }(a)}}+\overline{%
\tau _{a}(x)\chi _{\Omega _{W}\backslash B_{\varepsilon }(a)}}+\overline{%
\tau _{a}(x)\chi _{\Omega _{Z}}}  \label{32}
\end{equation}%
then, by Lemma \ref{muti} and Lemma \ref{pap1}, 
\begin{eqnarray*}
\int^{\ast }\sigma _{a}dx &=&\int^{\ast }\left[ \tau _{a}(x)-\breve{\tau}%
_{a}(x)\cdot \chi _{\Omega _{Z}}(x)\right] dx \\
&=&\int_{_{B_{\varepsilon }(a)}}^{\ast }\tau _{a}dx+\int_{_{\Omega
_{W}\backslash B_{\varepsilon }(a)}}^{\ast }\tau _{a}dx+\int_{_{\Omega
_{Z}}}^{\ast }\left[ \tau _{a}-\breve{\tau}_{a}\right] dx
\end{eqnarray*}%
By Lemma \ref{muti} and Lemma \ref{pap1}, we have that 
\begin{equation*}
\left\vert \int_{_{\Omega _{W}\backslash B_{\varepsilon }(a)}}^{\ast }\tau
_{a}dx+\int_{_{\Omega _{Z}}}^{\ast }\left[ \tau _{a}-\breve{\tau}_{a}\right]
dx\right\vert \leq 2\gamma ^{2}m\left( B_{\varepsilon }\right)
\end{equation*}%
Then, using (\ref{mozart1}) 
\begin{equation*}
\int^{\ast }\sigma _{a}dx\geq \left( 1-\gamma ^{2}\right) m\left(
B_{\varepsilon }\right) -2\gamma ^{2}m\left( B_{\varepsilon }\right) =\left(
1-3\gamma ^{2}\right) m\left( B_{\varepsilon }\right)
\end{equation*}%
and%
\begin{equation*}
\int^{\ast }\sigma _{a}dx\leq \left( 1+\gamma ^{2}\right) m\left(
B_{\varepsilon }\right) +2\gamma ^{2}m\left( B_{\varepsilon }\right) =\left(
1+3\gamma ^{2}\right) m\left( B_{\varepsilon }\right)
\end{equation*}

$\square $

\bigskip

\begin{corollary}
\label{123}If $u\in V_{\Lambda }$, then%
\begin{equation*}
\int^{\ast }\left\vert u(x)-\breve{u}(x)\right\vert dx\leq 3\gamma
\left\Vert u\right\Vert _{L^{\infty }}.
\end{equation*}
\end{corollary}

\textbf{Proof}: By Lemma \ref{ganzone} and (\ref{gamma}) we have that 
\begin{eqnarray*}
\int^{\ast }\left\vert u(x)-\breve{u}(x)\right\vert dx &=&\int^{\ast
}\dsum\limits_{a\in \Gamma }\left[ u(a)\left( \sigma _{a}(x)-\bar{\chi}%
_{S_{a}}(x)\right) \right] dx \\
&\leq &\dsum\limits_{a\in \Gamma _{W}}\left\vert u(a)\right\vert \int^{\ast
}\left\vert \sigma _{a}(x)-\bar{\chi}_{S_{a}}(x)\right\vert dx \\
&\leq &\dsum\limits_{a\in \Gamma _{W}}\left( \left\Vert u\right\Vert
_{L^{\infty }}\cdot 3\gamma ^{2}\right) =3\gamma ^{2}\left\Vert u\right\Vert
_{L^{\infty }}\left\vert \Gamma _{W}\right\vert \\
&\leq &3\gamma \left\Vert u\right\Vert _{L^{\infty }}
\end{eqnarray*}

\bigskip $\square $

\subsection{The pointwise integral and the generalized derivative\label{PI}}

Finally we set

\begin{equation}
V%
{{}^\circ}%
=\left\{ w|_{\Gamma }\ |\ w\in V_{\Lambda }\right\} .  \label{ann}
\end{equation}%
and we define the pointwise integral for the functions in $V%
{{}^\circ}%
\ $as follows:%
\begin{equation}
\doint u\ dx=\int^{\ast }u_{\Lambda }dx,  \label{chicca}
\end{equation}%
where $u_{\Lambda },$ as usual, is defined by (\ref{ciccia1}). Since $%
V_{\Lambda }$ is a hyperfinite space, then 
\begin{equation*}
V_{\Lambda }=\lim_{\lambda \uparrow \Lambda }\ V_{\lambda }
\end{equation*}%
for a suitable net of finite dimensional spaces $V_{\lambda }\subset V$ and
hence this definition agrees with (\ref{int}). In particular we have that%
\begin{equation*}
\doint u\ dx=\dsum\limits_{a\in \Gamma }u(a)d(a);\ d(a):=\int^{\ast }\sigma
_{a}(x)dx.
\end{equation*}%
So, if $w\in V_{\Lambda },$ then 
\begin{equation*}
\doint w%
{{}^\circ}%
\ dx=\int^{\ast }w\ dx,
\end{equation*}%
Moreover we have also the following interesting result:

\begin{theorem}
\label{leb}If $f\in \mathcal{L}^{1},$ $\doint f%
{{}^\circ}%
\ dx\sim \int f\ dx.$
\end{theorem}

\textbf{Proof}: First let us assume that $f$ be bounded. By corollary \ref%
{123}, we have that%
\begin{equation*}
\int^{\ast }\left\vert f%
{{}^\circ}%
(x)-\breve{f}(x)\right\vert dx\leq 3\gamma \left\Vert f\right\Vert
_{L^{\infty }}\sim 0
\end{equation*}
Then, by Th. \ref{super}, 
\begin{equation*}
\left\vert \doint f%
{{}^\circ}%
\ dx-\int f\ dx\right\vert \leq \left\vert \doint f%
{{}^\circ}%
\ dx-\int^{\ast }\breve{f}\ dx\right\vert +\left\vert \int^{\ast }\breve{f}\
dx-\int f\ dx\right\vert \sim 0.
\end{equation*}%
If $f$ is not bounded we argue as in the proof of Th. \ref{super}.

$\square $\bigskip

For every $u\in V_{\Lambda }$ we define the "generalized derivative"%
\begin{equation*}
D_{i}:V_{\Lambda }\rightarrow V_{\Lambda }
\end{equation*}%
as follows: given $u\in V_{\Lambda },\ D_{i}u$ is the only element in $%
V_{\Lambda }$ such that 
\begin{equation}
\doint D_{i}uv~dx=\int^{\ast }\partial _{i}^{\ast }uv\ dx\ \ \forall v\in
V_{\Lambda }  \label{deri+}
\end{equation}%
namely 
\begin{equation*}
D_{i}u=P\partial _{i}^{\ast }u
\end{equation*}%
where 
\begin{equation}
P:\mathfrak{M}^{\ast }\rightarrow V_{\Lambda }  \label{ppp}
\end{equation}%
is the "projection" of the natural extension of the space Radon measures $%
\mathfrak{M}$ over $V_{\Lambda }$ with respect to the duality 
\begin{equation*}
\left\langle w,v\right\rangle =\doint vw\ dx.
\end{equation*}

\begin{theorem}
\label{main}$V%
{{}^\circ}%
$ is a fine space of ultrafunctions.
\end{theorem}

\textbf{Proof}: We will check that $V%
{{}^\circ}%
$ verifies the requests of Def. \ref{A}.

\ref{A}.\ref{uno} follows from the construction of $V%
{{}^\circ}%
;$ in fact, by (\ref{F3}), we have that%
\begin{equation*}
u,v\in V^{\circledcirc }\Rightarrow \overline{uv}\in W_{\Lambda }^{0}\subset
W_{\Lambda }\subset V_{\Lambda }.
\end{equation*}

\ref{A}.\ref{11} follows from (\ref{chicca}).

\ref{A}.\ref{22} follows from Lemma \ref{ganzone} since $\chi _{a}=\left(
\sigma _{a}\right) 
{{}^\circ}%
$.

\ref{A}.\ref{YY} follows from Th. \ref{leb}.

\ref{A}.\ref{3} follows from (\ref{deri+}).

\ref{A}.\ref{2+} follows from the definition (\ref{deri+}) of $D_{i}$ since
also $\partial _{i}^{\ast }$ is a local operator.

$\square $

\section{Conclusive remarks}

The definition \ref{A} of fine ultrafunctions extends the notion of real
functions in such a way that:

\begin{itemize}
\item (a) "almost" all the partial differential equations (and other
functional problems) have a solution;

\item (b) some of the main properties of the smooth functions are preserved.
\end{itemize}

The assumptions included in definition \ref{A} seems quite natural; however
they do not define a unique model and there is a lot of room to require
others properties that allow to prove more facts about the solutions of a
given problem.

Let's illustrate this point with an example: let us consider equation (\ref%
{21}); we know that it has a unique solution which preserves the energy;
however there are a lot of questions that are relevant for the physical
interpretation such as:

\begin{itemize}
\item If $f(u)$ grows more than $\left\vert u\right\vert ^{\left( N+2\right)
/(N-2)}$ most likely, this solution has singular regions $S\subset \Gamma $
where the density of energy is an infinite number; is there a precise
theorem? What can we say about the properties of $S?$

\item In general the momentum is not "exactly" preserved since the "space" $%
\Gamma $ is not invariant for translation (and hence we cannot apply the
Noether's theorem); then the natural question is to know the initial
conditions for which the momentum is preserved and/or the initial conditions
for which the momentum is preserved up to an infinitesimal.

\item Under which assumptions the solutions converge to $0$ as $t\rightarrow
\infty ?$

\item And so on....
\end{itemize}

Some of these questions have an answer in the context of ultrafunctions
provided that you work enough on a single question. However there are
questions which do not have a YES/NO answer since there are models of
ultrafunctions in which the answer is YES and other in which the answer is
NO. A possible development of the theory is to add to the definition of
ultrafunctions other properties which allows a more detailed description of
the physical phenomenon described by the mathematical model. For example, in
Def. \ref{A} we have not included \textit{all} the properties of the
generalized derivative which can be deduced by the choice of $V%
{{}^\circ}%
$ and Def. (\ref{deri+}). The main difficulty, if you want to add a new
property to the ultrafunctions, is the proof of its consistency, namely the
construction of a model. For example it is easy to prove that the Leibniz
rule is not consistent with an algebra of functions which includes
idempotent functions (see the discussion in section \ref{RU}) but in
general, it is difficult to guess if a "reasonable property" is consistent
with all the others (for example the identity $D_{i}D_{j}=D_{j}D_{i}$).

At this points, it is interesting compare the ultrafunction theory (and more
in general a theory which includes infinitesimals) and a traditional theory
based on real valued function spaces. In both cases there exist questions
which do not have a YES/NO answer because of the Goedel incompleteness
theorem, and, in both cases, we can add new axioms which allow to solve the
problem. The difficult issue is the proof of the consistency of these new
axioms. However, this issue is much easier in the world of Non-Archimedean
mathematics, since the mathematical universe is wider and there is more room
to construct a model.

Let us consider an example that clarifies this point. We do not know if the
Navier-Stokes equations have a smooth solution; however they have a unique
ultrafunction solution. Probably this fact is not relevant for an applied
mathematician or for an engineer since it does not help to discover new fact
relative to the motion of a real fluid. However, it is possible to add new
axioms which allow to prove properties of the solutions consistent with the
experiments and to end with a richer mathematical model. This model can be
used to get new practical results even if the problem of the existence of a
smooth solutions remains unsolved.

For this reasons we think that it is worthwhile to investigate the
potentialities of the Non-Archimedean mathematics.

\bigskip

\bigskip

\noindent \textbf{Acknowledgement}: I want to express my gratitude to the
referee whose remarks have been very useful.

\end{document}